\ifdefined\journalstyle
  \documentclass[numbers,webpdf,imanum]{ima-authoring-template}%
\else
  \documentclass[11pt]{article}
  \usepackage[a4paper]{geometry}
  \usepackage{latexsym,amsfonts,amsmath,theorem,amssymb,url}
  \usepackage{graphicx}
  \usepackage{hyperref}
  \usepackage[usenames,dvipsnames,svgnames]{xcolor}
  \definecolor{myblue}{rgb}{0,0,0.6}
  \hypersetup{colorlinks=true, linkcolor=myblue,citecolor=myblue,filecolor=myblue,urlcolor=myblue}
\fi

\usepackage{mleftright}
\usepackage{pdfsync}

\usepackage{tikz}
\usetikzlibrary{patterns.meta}

\ifdefined\journalstyle
  \theoremstyle{thmstyleone}
  \theoremstyle{thmstyletwo}
  \newtheorem{theorem}{Theorem}[section]         
  \newtheorem{proposition}[theorem]{Proposition} 
  \newtheorem{lemma}[theorem]{Lemma}             
  \newtheorem{corollary}[theorem]{Corollary}     
  \newtheorem{example}[theorem]{Example}  
  \newtheorem{remark}[theorem]{Remark}    
  \newtheorem{definition}[theorem]{Definition}   
  \newtheorem{notation}[theorem]{Notation}       
  \newtheorem{assumption}[theorem]{Assumption}   
  \theoremstyle{thmstylethree}
  \numberwithin{equation}{section}

\else
  \newtheorem{theorem}{Theorem}[section]
  \newtheorem{assumption}[theorem]{Assumption}
  \newtheorem{lemma}[theorem]{Lemma}
  \newtheorem{corollary}[theorem]{Corollary}
  \newtheorem{proposition}[theorem]{Proposition}
  
  \newtheorem{definition}[theorem]{Definition}
  
  \newtheorem{remark}[theorem]{Remark}
  \newtheorem{example}[theorem]{Example}
  \newtheorem{notation}[theorem]{Notation}
  \newenvironment{proof}{\begin{trivlist}
  \item[\hskip\labelsep{\it Proof.}]}{$\hfill\Box$\end{trivlist}}
  \numberwithin{equation}{section}
\fi

\newcommand{\beq}{\begin{equation}}
\newcommand{\eeq}{\end{equation}}
\newcommand{\beqs}{\begin{equation*}}
\newcommand{\eeqs}{\end{equation*}}
\newcommand{\bit}{\begin{itemize}}
\newcommand{\eit}{\end{itemize}}
\newcommand{\ben}{\begin{enumerate}}
\newcommand{\een}{\end{enumerate}}
\newcommand{\bal}{\begin{align}}
\newcommand{\eal}{\end{align}}
\newcommand{\bals}{\begin{align*}}
\newcommand{\eals}{\end{align*}}
\newcommand{\bpr}{\begin{proposition}}
\newcommand{\epr}{\end{proposition}}
\newcommand{\bre}{\begin{remark}}
\newcommand{\ere}{\end{remark}}
\newcommand{\bpf}{\begin{proof}}
\newcommand{\epf}{\end{proof}}
\newcommand{\ble}{\begin{lemma}}
\newcommand{\ele}{\end{lemma}}
\newcommand{\bco}{\begin{corollary}}
\newcommand{\eco}{\end{corollary}}
\newcommand{\bex}{\begin{example}}
\newcommand{\eex}{\end{example}}
\newcommand{\bth}{\begin{theorem}}
\newcommand{\enth}{\end{theorem}}
\newcommand{\bcon}{\begin{condition}}
\newcommand{\econ}{\end{condition}}
\newcommand{\bas}{\begin{assumption}}
\newcommand{\eas}{\end{assumption}}
\newcommand{\bde}{\begin{definition}}
\newcommand{\ede}{\end{definition}}
\newcommand{\ton}{\text{ on }}
\newcommand{\tin}{\text{ in }}
\newcommand{\tfa}{\text{ for all }}

\newcommand{\tand}{\text{ and }}

\newcommand{\CC}{\mathbb C}
\newcommand{\grad}{\nabla}
\newcommand{\N}[1]{\mleft\|#1\mright\|}
\newcommand{\pdiff}[2]{\frac{\partial #1}{\partial #2}}
\newcommand{\tendi}{\rightarrow \infty}
\newcommand{\DtN}{{\rm DtN}_k}
\newcommand{\de}{:=}
\newcommand{\vbar}{\overline{v}}
\newcommand{\minispace}{\;\!}

\newcommand{\NTA}[1]{\mathrm{NT}_{\mathrm{mat}}\mleft(#1\mright)}
\newcommand{\NTn}[1]{\mathrm{NT}_{\mathrm{scal}}\mleft(#1\mright)}
\DeclareMathOperator*{\esssup}{ess\,sup}

\newcommand{\cD}{\mathcal D}

\newcommand{\cG}{\mathcal G}
\newcommand{\cH}{\mathcal H}

\newcommand{\aPML}{a_{\mathrm{PML}}}
\newcommand{\uPML}{u_{\mathrm{PML}}}
\newcommand{\Etrunc}{E^{\rm trunc}}
\newcommand{\EQMC}{E^{\rm QMC}}
\newcommand{\EFEM}{E^{\rm FEM}}

\newcommand{\Gdir}{\partial D}

\newcommand{\bsalpha}{{\boldsymbol{\alpha}}}
\newcommand{\bsDelta}{{\boldsymbol{\Delta}}}

\newcommand{\bsb}{{\boldsymbol{b}}}

\newcommand{\bse}{{\boldsymbol{e}}}

\newcommand{\bsgamma}{{\boldsymbol{\gamma}}}

\newcommand{\bstau}{{\boldsymbol{\tau}}}
\newcommand{\bsnu}{{\boldsymbol{\nu}}}
\newcommand{\bszeta}{{\boldsymbol{\zeta}}}

\newcommand{\bst}{{\boldsymbol{t}}}

\newcommand{\bsx}{{\boldsymbol{x}}}

\newcommand{\bsy}{{\boldsymbol{y}}}
\newcommand{\bsz}{{\boldsymbol{z}}}

\newcommand{\bsm}{{\boldsymbol{m}}}

\newcommand{\bszero}{{\boldsymbol{0}}}

\newcommand{\ri}{{\mathrm{i}}}
\newcommand{\rd}{{\mathrm{d}}}

\newcommand{\bshalf}{{\boldsymbol{\tfrac{1}{2}}}}

\newcommand{\Honezero}{H^1_{0,D}(D_R)}
\newcommand{\bbN}{{\mathbb{N}}}
\newcommand{\bbR}{{\mathbb{R}}}

\newcommand{\bbE}{{\mathbb{E}}}

\newcommand{\calG}{{\mathcal{G}}}

\newcommand{\calO}{{\mathcal{O}}}
\newcommand{\calW}{{\mathcal{W}}}

\newcommand{\supp}{{\mathrm{supp}}}
\newcommand{\setu}{{\mathfrak{u}}}
\newcommand{\setv}{{\mathfrak{v}}}

\newcommand{\ualt}{u^{\rm alt}}
\newcommand{\falt}{f^{\rm alt}}
\newcommand{\Coscil}{C_{\rm osc}}
\newcommand{\Dminus}{D}
\newcommand{\Cstab}{C_{\mathrm{stab}}}

\begin{document}

\ifdefined\journalstyle
  \DOI{DOI HERE}
  \copyrightyear{2025}
  \vol{00}
  \pubyear{2025}
  \access{Advance Access Publication Date: Day Month Year}
  \appnotes{Paper}
  \copyrightstatement{Published by Oxford University Press on behalf of the Institute of Mathematics and its Applications. All rights reserved.}
  \firstpage{1}


  \title[Quasi-Monte Carlo methods for wave propagation and scattering]{Quasi-Monte Carlo methods for uncertainty quantification of wave propagation and scattering problems modelled by the Helmholtz equation}

  \author{\textsc{Ivan G.~Graham*}
          \address{Department of Mathematical Sciences, University of Bath, United Kingdom}}
  \author{\textsc{Frances~Y.~Kuo}
          \address{School of Mathematics and Statistics, UNSW Sydney, Australia}} 
  \author{\textsc{Dirk Nuyens}
          \address{Department of Computer Science, KU Leuven, Belgium}} 
  \author{\textsc{Ian H.~Sloan}
          \address{School of Mathematics and Statistics, UNSW Sydney, Australia}}
  \author{\textsc{Euan A.~Spence}
         \address{Department of Mathematical Sciences, University of Bath, United Kingdom}}

  \authormark{Graham et al.}

  \corresp[*]{Corresponding author. \href{email:i.g.graham@bath.ac.uk}{i.g.graham@bath.ac.uk}}

  \received{Date}{0}{2025}
  \revised{Date}{0}{2025}
  \accepted{Date}{0}{2025}


  \abstract{We analyse and implement a quasi-Monte Carlo (QMC) finite element method
(FEM) for the forward problem of uncertainty quantification (UQ) for the
Helmholtz equation with random coefficients, both in the second-order and
zero-order terms of the equation, thus modelling wave scattering in random
media. The problem is formulated on the infinite propagation domain, after
scattering by the heterogeneity, and also (possibly) a bounded
impenetrable scatterer. The spatial discretization scheme includes
truncation to a bounded domain via a perfectly matched layer (PML)
technique and then FEM approximation. A special case is the problem of an
incident plane wave being scattered by a bounded sound-soft impenetrable
obstacle surrounded by a random heterogeneous medium, or more simply, just
scattering by the random medium. The random coefficients are assumed to be
affine separable expansions with infinitely many independent uniformly
distributed and bounded random parameters. As quantities of interest for
the UQ, we consider the expectation of general linear functionals of the
solution, with a special case being the far-field pattern of the scattered
field. The numerical method consists of (a) dimension truncation in
parameter space; (b) application of an adapted QMC method to compute
expected values and (c) computation of samples of the PDE solution via PML
truncation and FEM approximation. Our error estimates are explicit in $s$
(the dimension truncation parameter), $N$ (the number of QMC points), $h$
(the FEM grid size) and (most importantly), $k$~(the Helmholtz wavenumber). The
method is also exponentially accurate with respect to the PML truncation
radius. Illustrative numerical experiments are given.}


  \maketitle

\else
  \title{Quasi-Monte Carlo methods for uncertainty quantification of \\
  wave propagation and scattering problems\\ modelled by the Helmholtz equation}

 \author{Ivan G.~Graham\footnote{
 Department of Mathematical Sciences, University of Bath, United Kingdom. Email: i.g.graham@bath.ac.uk}
 , 
 Frances~Y.~Kuo\footnote{
 School of Mathematics and Statistics, UNSW Sydney, Australia. Email: f.kuo@unsw.edu.au} 
 , 
 Dirk Nuyens\footnote{
 Department of Computer Science, KU Leuven, Belgium. Email: dirk.nuyens@kuleuven.be} 
 ,
 Ian H.~Sloan\footnote{
 School of Mathematics and Statistics, UNSW Sydney, Australia. Email: i.sloan@unsw.edu.au} 
 , 
 Euan A.~Spence\footnote{
 Department of Mathematical Sciences, University of Bath, United Kingdom. Email: e.a.spence@bath.ac.uk}
 }

 \date{November 2025}

\maketitle

\begin{abstract}
We analyse and implement a quasi-Monte Carlo (QMC) finite element method
(FEM) for the forward problem of uncertainty quantification (UQ) for the
Helmholtz equation with random coefficients, both in the second-order and
zero-order terms of the equation, thus modelling wave scattering in random
media. The problem is formulated on the infinite propagation domain, after
scattering by the heterogeneity, and also (possibly) a bounded
impenetrable scatterer. The spatial discretization scheme includes
truncation to a bounded domain via a perfectly matched layer (PML)
technique and then FEM approximation. A special case is the problem of an
incident plane wave being scattered by a bounded sound-soft impenetrable
obstacle surrounded by a random heterogeneous medium, or more simply, just
scattering by the random medium. The random coefficients are assumed to be
affine separable expansions with infinitely many independent uniformly
distributed and bounded random parameters. As quantities of interest for
the UQ, we consider the expectation of general linear functionals of the
solution, with a special case being the far-field pattern of the scattered
field. The numerical method consists of (a) dimension truncation in
parameter space; (b) application of an adapted QMC method to compute
expected values and (c) computation of samples of the PDE solution via PML
truncation and FEM approximation. Our error estimates are explicit in $s$
(the dimension truncation parameter), $N$ (the number of QMC points), $h$
(the FEM grid size)
and (most importantly), $k$~(the Helmholtz wavenumber). The
method is also exponentially accurate with respect to the PML truncation
radius. Illustrative numerical experiments are given.
\end{abstract}


\smallskip\noindent
\textbf{MSC2020 Subject Classification:} 60G60, 65D30, 65D32, 65N30, 35Q86


\fi

\section{Introduction}\label{sec:intro}

\subsection{Informal statement of the exterior Dirichlet problem for the Helmholtz equation
and the UQ algorithm}

In this paper we study the uncertainty quantification of problems governed
by an exterior Dirichlet problem for the Helmholtz equation, posed in the
infinite region  $D_+:= \mathbb{R}^d \setminus \overline{\Dminus}$ where
$D \subset \mathbb{R}^d$, $d \in \{2,3\}$, is a bounded domain (which may
be empty). When $D \not = \emptyset$ we impose a homogeneous Dirichlet
condition on its boundary $\Gdir$. Thus the general form of the problem is
to find $u : D_{+} \to \CC$ such that
\begin{align}
 \grad \cdot \mleft(A \grad u\mright) + k^2 n u &\,=\, -f \quad\,\tin D_+,  \label{eq:intro_tedp}
 \\
u &\,=\, 0 \qquad\ton \Gdir, \label{eq:intro_dbc}
\end{align}
where $A$ is a positive-definite matrix-valued function, $n$ is a positive
scalar-valued function, $f$ is a forcing term, and $k>0$ is the Helmholtz wavenumber. We also require
$u$ to satisfy the Sommerfeld radiation condition,
 \beq\label{eq:src}
 \pdiff{u}{r}(\bsx) - \ri k u(\bsx) = o \Big(\frac{1}{r^{(d-1)/2}}\Big)
 \eeq
as $r:= |\bsx|_{2}\tendi$ (uniformly in $\widehat{\bsx}:= \bsx/r$),
thus ensuring that the asymptotic solution is an outgoing spherical wave,
and making the solution unique. The cost of solving the Helmholtz equation numerically increases as $k$ increases; we therefore pay particular attention to how our method and error depend on~$k$.

In this paper the coefficients $A$ and $n$ are heterogeneous (with the
random case being of  most interest), but we always assume that $f$ has
bounded support  and that the problem is homogeneous in the far
field.

Such problems arise widely in applications. The propagation of acoustic
waves in a homogeneous medium is often modelled by \eqref{eq:intro_tedp} with $A=I$ and $n =
c^{-2}$, where $c$ is the sound speed, and $k$ is the angular
frequency of the excitation, also called the wavenumber, see, e.g.,
\cite[Chapter 8]{CoKr:13}.

In seismology, problem \eqref{eq:intro_tedp}--\eqref{eq:src} arises when
elastic waves are computed using the `acoustic aproximation', in which
case $A = \rho^{-1}I$ and $n = \rho^{-1} c_p^{-2}$, where $\rho$ is the
density and $c_p$ is the speed of longitudinal waves (also called
$P$-waves) in the medium, see, e.g., \cite{Ta:84}. In `full waveform
inversion' the fields $n$ and possibly $A$ are to be recovered in an
inverse problem, requiring many solves of the above problem, and, for each
solve, $f$ is given by a point source.

For any radius $R>0$, we define the open ball
\begin{align} \label{not_BR}
 B_{R}:= \{ \bsx \in \mathbb{R}^d: |\bsx|_2 < R\},
 \quad \text{and the domain} \quad  D_R : = B_R\backslash \overline{D}.
\end{align}
Since the obstacle $D$ is bounded, the region $D_R$ is non-empty for a
large enough $R$.

We assume throughout that the supports of $A-I$ and $n-1$ are bounded.
Thus we can choose a radius $R_0>0$ such that these supports, and the
closure of the obstacle $\overline{D}$, are all contained inside the ball
$B_{R_0}$ and are sufficiently far from its boundary $\partial B_{R_0}$.
Thus we have $A \equiv I$ and $n \equiv 1$ outside $B_{R_0}$.
(Assumption~\ref{ass:1} and Definition~\ref{def:EDP} below give a more
precise statement of the problem, including regularity requirements.) See
Figure~\ref{fig:intro} for an illustration.

\tikzset{
  ring shading/.code args={from #1 at #2 to #3 at #4}{
    \def\colin{#1}
    \def\radin{#2}
    \def\colout{#3}
    \def\radout{#4}
    \pgfmathsetmacro{\proportion}{\radin/\radout}
    \pgfmathsetmacro{\outer}{.8818cm}
    \pgfmathsetmacro{\inner}{.8818cm*\proportion}
    \pgfmathsetmacro{\innerlow}{\inner-0.01pt}
    \pgfdeclareradialshading{ring}{\pgfpoint{0cm}{0cm}}%
    {
      color(\inner)=(#1);
      color(\outer)=(#3)
    }
    \pgfkeysalso{/tikz/shading=ring}
  },
}
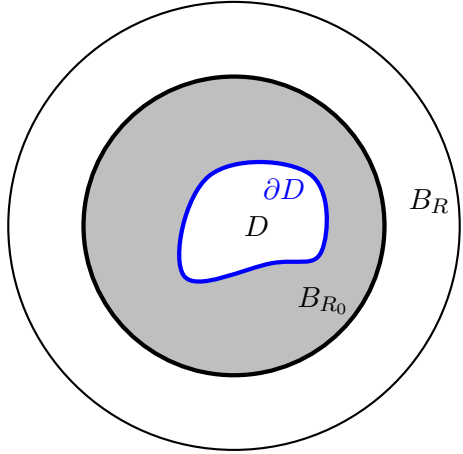
\begin{figure}[t]
\begin{center}
\begin{tikzpicture}[scale=0.33]
 \draw[thick] (0,0) circle (9); 
 \draw[ultra thick,fill=lightgray] (0,0) circle (6);   
 \draw[ultra thick, pattern={Lines[angle=45,distance=8pt]}, blue, fill=white]
   plot[smooth cycle, tension=.7]
   coordinates {(-2,-2) (-1,2) (3,2.1) (3.5,-1) (1.5,-1.5)};
 \node[right] at (0,0) {$D$};
 \node[blue] at (2,1.5) {$\partial D$};
 \node[left] at (5,-3) {$B_{R_0}$};
 \node[left] at (9,1) {$B_{R}$};
\end{tikzpicture}
\end{center}
\caption{The domain $D_{R_0}$ is the (open) shaded region inside the ball $B_{R_0}$,
excluding the closure of the obstacle~$D$.
The supports of $A-I$ and $n-1$ are in $D_{R_0}$. The support of the source term
$f$ is in $\overline{D_{R_0}}$.
The quantity of interest $G$ will be  a linear functional on $D_R$ for $R \geq R_0$.}
\label{fig:intro}
\end{figure}

\paragraph{Random coefficients and the UQ problem.}

While there is a further brief discussion in \S\ref{subsec:determin} of
the deterministic version of \eqref{eq:intro_tedp}--\eqref{eq:src}, we mainly study the random case, in which the coefficients $A$ and $n$ 
are taken to be random fields with the infinite affine expansions
\begin{align}
 A(\bsx,\bsy) &\,=\, {A}_0(\bsx) + \xi_A\,
  \sum_{j=1}^\infty y_j\,\Psi_j(\bsx), && \bsx\in D_+, \; \bsy\in U, \label{eq:Axy}\\
 n(\bsx,\bsy) &\,=\, {n}_0(\bsx) + \xi_n\,
  \sum_{j=1}^\infty y_j\,\psi_j(\bsx), && \bsx\in D_+, \; \bsy\in U, \label{eq:nxy}
\end{align}
where $U := [-\frac{1}{2},\frac{1}{2}]^\bbN$ and $\bsy = (y_1, y_2, \ldots )$ is an
infinite sequence of independent uniformly distributed random parameters
$y_j \in [-\frac{1}{2},\frac{1}{2}]$, while $\xi_A$, $\xi_n$ are nonnegative scaling
parameters. In \S\ref{subsec:random_prob}, conditions are given to ensure
that the infinite expansions in \eqref{eq:Axy} and \eqref{eq:nxy} are
absolutely convergent.

With this formulation, the fields $A$ and $n$ can be either correlated or
not. To model two fields that are not correlated, it suffices to set
$\Psi_{2j-1} \equiv 0$ and $\psi_{2j} \equiv 0$, for all $j \geq 1$. On
the other hand, if there exists a $j \geq 1$ such that $\Psi_j \not\equiv
0$ and $\psi_j \not\equiv 0$, then $A$ and $n$ will be correlated. In
any case the functions $\Psi_j$ that are nonzero are assumed to be
linearly independent and known, and likewise the nonzero functions
$\psi_j$.

The output solution $u(\bsx,\bsy)$ depends on the spatial variable
$\bsx\in D_+$, as well as  the  random parameters $\bsy \in U$, and the UQ
problem that we study is to efficiently compute (and give error estimates
for)  approximations of the expected value of $Gu(\cdot, \bsy)$, where $G$
is some linear functional acting with respect to the spatial variable
$\bsx$ in $D_R = B_R\setminus\overline{D}$ with $R\ge R_0$. The expected
value of any function $\Theta: U\rightarrow \mathbb{C}$ is expressed
as an integral,
\begin{align} \label{eq:def_exp}
  \bbE[\Theta] \ :=\  \int_{U} \Theta(\bsy) \, \rd \bsy
  \ =\  \lim_{s\to\infty}
  \int_{U_s} \Theta(y_1,\ldots,y_s, 0,0, \ldots)c \, \rd y_1 \cdots \rd y_s,
\end{align}
when this limit exists, where $U_s := [-\frac{1}{2},\frac{1}{2}]^s$. 
This limit exists under the conditions required for Legesgue's dominated convergence theorem, see e.g., \cite[Theorems~1 and~2]{KNPSW17}.

\paragraph{The algorithm.}

We start by truncating the affine expansions \eqref{eq:Axy} and
\eqref{eq:nxy} to $s$ terms, yielding approximations to $A$ and $n$ that
depend only on $\bsy_s: = (y_1, \ldots, y_s)$, and then denote the
corresponding approximate solution of
\eqref{eq:intro_tedp}--\eqref{eq:src} by
\begin{align} \label{eq:defus}
  u_s(\cdot, \bsy_s): = u(\cdot, (\bsy_s,\boldsymbol{0})).
\end{align}
Taking the expectation of this yields a first (theoretical) approximation
\begin{align} \label{eq:theoretical}
 \bbE[Gu] \approx \bbE[Gu_s] .
\end{align}

To attack the $s$-dimensional integral arising in the computation of the
expectation in \eqref{eq:theoretical}, we will consider a family of
first order quasi-Monte Carlo (QMC) quadrature rules known as ``randomly
shifted lattice rules'', as well as a family of higher order QMC rules
known as ``interlaced polynomial lattice rules''. In the case of randomly
shifted lattice rules, we construct $N$ suitable lattice points in
the $s$-dimensional unit cube $[0, 1]^s$ with specially chosen generating
vector $\bsz \in \bbN^s$ and a uniformly-distributed random shift
$\bsDelta \in [0, 1)^s$ (cf., \cite{DiKuSl:13}).
The entire pointset is then translated into the cube $[-\frac{1}{2},\frac{1}{2}]^s$.
The estimate of $\bbE[Gu_s]$ is then the (equal-weight) average of the
values of $G u_s$ at these $N$ shifted points in~$U_s$, which we denote by
$Q_{s,N,\bsDelta} (Gu_s)$. Full details of the QMC algorithm are
in~\S\ref{sec:QMC}. In general $N$ may need to increase with
increasing $s$ or increasing $k$.

Finally, to make this practical, we employ a finite element method (FEM)
to approximate each sample $u_s$ (for $\bsy_s \in U_s$) by a spatially
discrete approximation denoted $u_{s,{\rm PML}, h}$, where $h$ denotes
mesh diameter and PML indicates that we employ a perfectly matched layer
(PML) technique to approximate the PDE over the infinite domain $D_+$
by an accurate approximation on a finite domain; details are
given in \S\ref{sec:FEM}. This yields the practical approximation of the
required expected value, which is denoted $Q_{s,N,\bsDelta} (Gu_{s,{\rm
PML}, h})$.

To measure errors we need appropriate function spaces.

\begin{notation}[Sobolev spaces and dual spaces] \label{not:Sobolev}
Throughout the paper, for any bounded domain $\Omega$ and any $r \geq 0$,
$H^r(\Omega)$ denotes the usual Sobolev space of order $r\geq 0$, with
$H^0(\Omega) = L^2(\Omega)$ and $H^r(\Omega)'$ denoting the dual
space of $H^r(\Omega)$ (i.e., the space of all bounded linear functional
on $H^r(\Omega)$). Standard norms on $H^r(\Omega)$ and $H^r(\Omega)'$  are
denoted $\Vert \cdot \Vert_{H^r(\Omega)}$ and   $\Vert  \cdot
\Vert_{H^r(\Omega)'}$. To be precise, the scales $H^r(\Omega)$ and $H^r(\Omega)'$ with varying $r$ are
interpolation spaces \cite[\S 14.2.3]{BrSc:08}, \cite[\S 4.5]{BeLo:76}.

In the special case $r = 1$, the standard norm is $\|v\|_{H^1(\Omega)}^2
:= \|\nabla v\|_{L^2(\Omega)}^2 + \|v\|_{L^2(\Omega)}^2$, but we 
also make extensive use of the $k$-weighted $H^1(\Omega)$ norm, defined for
$k>0$ by
\begin{align} \label{eq:weighted}
 \Vert v \Vert_{H^1_k(\Omega)}^2
 := \Vert \nabla v \Vert_{L^2(\Omega)}^2  + k^2 \Vert v \Vert_{L^2(\Omega)}^2 .
\end{align}
We also use the space $H^1_0(\Omega) := \{ v\in H^1(\Omega) : v = 0
\mbox{ on } \partial\Omega\}$ equipped with the norm \eqref{eq:weighted}.

Although the solution $u$ of problem \eqref{eq:intro_tedp}--\eqref{eq:src}
is defined on all of $D_+$, we will typically be concerned with its
restriction to the domain $D_R$ for some $R \geq R_0$, see \eqref{not_BR}
and Figure~\ref{fig:intro}. In this context, we write $u \in
H_{0,D}^1(D_R)$ to indicate that $u$ has vanishing trace on $\partial D$
(see \eqref{eq:spaceEDP}).
\end{notation}

\subsection{The main results of the paper}\label{sec:main_results}

In this introductory subsection we restrict ourselves to an illustration
of the main results, with pointers to the full details given later. A key
assumption used throughout this paper (and in many other papers about UQ
of PDEs, see e.g., \cite{CD15}) is that the expansions \eqref{eq:Axy} and
\eqref{eq:nxy} converge quickly enough. The speed of convergence is
controlled by a ``summability exponent'' $p\in (0,1)$, see \eqref{eq:bj}
below, with a smaller value of $p$ corresponding to a faster convergence
of \eqref{eq:Axy} and \eqref{eq:nxy}.

In our error analysis, we also impose conditions on the coefficients
$A,n$, and obstacle $D$ to ensure good (i.e., nontrapping) behaviour
of the Helmholtz solution operator (see Definition~\ref{def:nontrapping}
below). Although we do not impose all of these conditions in our numerical
experiments, we do not observe any ``bad'' behaviour of the Helmholtz
solution operator. This is because such bad behaviour is very specific to
particular values of $k$ and particular data, and so it is rare that one
actually sees this behaviour in practice; see the discussion in
\S\ref{sec:trapping}.

In this illustration we restrict ourselves to the case when $G \in
L^2(D_R)'$, although bounded linear functionals on $H^r(D_R)$ for any
$r \in [0,1]$ are also treated in the detailed theorems. We restrict
our illustration here to randomly shifted lattice rules, while the theory of interlaced
polynomial lattice rules is also considered in \S\ref{sec:QMC}.

Our overall error estimate for the root-mean-square error has the
following form:
\begin{align}\label{eq:tentative}
  & \left( \bbE_{\bsDelta} \left[\,
  \left|\,\bbE[G u]-  Q_{s,N, \bsDelta} (G u_{s,{\rm PML},h}) \,\right|^2\right] \,
  \right)^{1/2}
  \ \leq \ \sqrt{2} \left( \Etrunc  + \EQMC + \EFEM \right),
\end{align}
where $\bbE_{\bsDelta}$ denotes expectation with respect to the random shift $\bsDelta$ and 
\begin{align} \label{eq:tentative1}
  \Etrunc & := \left|\, \bbE [ G(u -u_s)] \,\right|,\nonumber \\
  \EQMC & :=\left( \bbE_{\bsDelta}  \left[ \left| \bbE [Gu_s]  - Q_{s,N,\bsDelta} (G u_s) \right|^2
          \right]\right)^{1/2},\nonumber \\
  \EFEM & :=
  \sup_{\bsy_s\in U_s} \left| G (u_s - u_{s,{\rm PML},h})(\bsy_s) \right|.
\end{align}

Our illustration includes the case when the linear functional $G$ is the
far-field pattern of the scattered wave produced when a plane wave is
incident on a heterogeneous medium, with or without an impenetrable
obstacle.

In the analysis we use the notation
\begin{align}\label{eq:def_xi}
  \xi := \xi_A + \xi_n \ ,
\end{align}
where $\xi_A$, $\xi_n$ are the scaling parameters in \eqref{eq:Axy},
\eqref{eq:nxy}, and we shall see below  that $k\,R\,\xi$ is an important critical
parameter in the error estimates.

Throughout the paper, for two expressions $A$ and $B$, we write $A
\lesssim B$ to mean $A \leq  cB $ where $c$ is a constant independent of
$k,R,\xi,N$, $h$ and $s$, and we write $A \sim B$ to mean $A \lesssim B$
and $B \lesssim A$.

\paragraph{Dimension truncation.}

The detailed estimate for the dimension truncation error is given in
Theorem~\ref{thm:dim_trunc}. In our case, assuming that $G\in L^2(D_R)'$
(the dual space of $L^2(D_R)$, which can be identified with $L^2(D_R)$) for some $R \geq R_0$, the estimate with the fastest decay with respect to
$s$ is
$$
\Etrunc \ \lesssim \  R\,k^{-1}\, c(k\,R\,\xi) \, s^{-2/p +1} \, \Vert f \Vert_{L^2(D_{R_0})} \,
\Vert G \Vert_{L^2 (D_R)'},
$$
although here $c$ is a function that grows rapidly with respect to
its argument. In the case when $k$ is large the estimate
$$
 \Etrunc \ \lesssim \ R^2 \,\xi \, s^{-1/p +1} \, \Vert f \Vert_{L^2(D_{R_0})} \, \Vert G \Vert_{L^2(D_R)'}
$$
is available, providing a slower decay rate with respect to $s$ but
independent of $k \rightarrow \infty$.

\paragraph{QMC error.}

The QMC error is analyzed in Theorem~\ref{thm:qmc}. Given any linear
functional $G\in L^2(D_R)'$ for some $R\ge R_0$, we show that randomly
shifted lattice rules can be constructed to achieve
$$
 E^{\rm QMC}
 \ \le \ R\,k^{-1}\, C(s,k,R,\xi,\delta)\,N^{-1+\delta}\,
 \Vert f \Vert_{L^2(D_{R_0})} \, \Vert G \Vert_{L^2(D_R)'}
 , \quad \text{if $p \le 2/3$},
$$
where the constant $C(s,k,R,\xi,\delta)$ grows unboundedly as $\delta \to 0$. Moreover: (i) if
$k$ is bounded, then $C(s,k,R,\xi,\delta)$ is bounded independently of
$s$; (ii) if $s$ is fixed and $k$ grows, then, with $\Cstab$ denoting the
stability constant for the boundary value problem
\eqref{eq:intro_tedp}--\eqref{eq:src} written in variational form (see
Corollary~\ref{cor:random}), we have $C(s,k,R,\xi,\delta) \sim (\Cstab\,
k\,R\,\xi)^s$, which grows quickly with increasing dimension unless
\begin{align}
 \xi \leq \frac{1}{\Cstab\,k\,R}, \label{weakly_stoch}
\end{align}
in which case the magnitude of the random perturbation dies away as
$k\to\infty$. A higher order convergence rate $N^{-1/p}$ can be
achieved using interlaced polynomial lattice rules.

\paragraph{FEM error.}

The PML truncation and FEM approximation of the problem are described in
\S\ref{sec:FEM}. Given any linear  functional $G \in L^2(D_R)'$  for some
$R \geq R_0$, the PML truncation should be performed in a spherical
annulus  $B_{R_2}\backslash B_{R_1}$ with $R_2 > R_1 > R$. We approximate
the truncated problem using $H^1$-conforming Lagrange (in the sense of \cite[p.95]{Ci:91}) 
finite elements of
degree $m\geq 1$. Then, assuming that all the data $A$, $n$ and $\partial
D$ are sufficiently smooth, and provided the mesh-resolution condition
\begin{align} \label{eq:resolution}
  (hk)^{2m} k \quad \text{sufficiently small}
  \end{align}
is satisfied, the dimension-truncated and PML-truncated problem has a
unique FEM solution $u_{s,{\rm PML}, h}$. Our estimate of $\EFEM$ then
takes the form (see Corollary \ref{cor:functional})
\begin{align} \EFEM \ \lesssim\ k^{-1}\,
 \Big(\exp(-C_{\rm PML}\,k) +  ((hk)^{m+1} + (hk)^{2m}\, k)\Big)\,
 \Vert f \Vert_{L^2(D_{R_0})}\, \Vert G \Vert_{L^2(D_R)'} , \label{eq:introFEM}
\end{align}
where both the hidden positive constant and the explicit positive constant
$C_{\rm PML}$ depend on the PML parameters but not on $h,s$ or $k$.

For the particular case when $G$ is the far-field pattern of the scattered
field in a plane-wave sound-soft scattering problem (see
\S\ref{subsec:SSSP} and \S\ref{sec:SSSP}), then the estimate of $\EFEM$
takes the form (Theorem \ref{thm:FEM-ffp})
\begin{align}
 \EFEM \ \lesssim\ c(d,k)\,k\,
 \Big(\exp(-C_{\rm PML}\,k) + ((hk)^{m+1} + (hk)^{2m}\, k)\Big),  \label{eq:introFEMFFP}
\end{align}
where $c(d,k)$ is independent of $k$ when $d = 3$ and decays with
$\mathcal{O}(k^{-1/2})$ when $d = 2$.

\subsection{Discussion of the main results in the context of the literature}

While there have been many papers analysing the application of QMC to PDE
problems, much of this has dealt with coercive PDEs, thus not applicable to high-frequency Helmholtz. 
However we want to mention here three recent papers on the $k$-explicit UQ of the Helmholtz equation:~\cite{GaKuSl}, \cite{HSS}, \cite{GrMaZh:20}. 
\bit
\item 
Reference \cite{GaKuSl} considers the Helmholtz equation \eqref{eq:intro_tedp} with $A=I$ and variable coefficient $n$, and an impedance boundary condition. The novelties of the present paper relative to \cite{GaKuSl} are that 
(i) we do not study impedance boundary conditions (which give $k$-independent errors when approximating the radiation condition \eqref{eq:src} for high-frequency waves \cite{GLS2}) but instead study PML truncation (which gives exponentially-small errors in approximating \eqref{eq:src}); 
(ii) we consider the physically-relevant quantity of the far-field pattern and give detailed numerical experiments; and 
(iii) the finite element error estimates in \cite{GaKuSl} are not $k$-explicit  (see \cite{GaKuSl:cor}),  while here we give fully $k-$explicit error estimates for all approximations. 

\item Reference \cite{HSS} gives $k$-explicit theory for shape UQ (different from the setting in this paper) for the Helmholtz transmission problem. The novelties of the present paper relative to \cite{HSS} are that 
(i) we have randomness in the highest-order term of the PDE (\cite{HSS} only has randomness in lowest-order term); and 
(ii) we give numerical experiments supporting the theory (while \cite{HSS} is theory only).

\item While reference \cite{GrMaZh:20} presents UQ methods with error independent of $k$ (or even improving as $k \rightarrow \infty$) for a high random-dimensional Helmholtz problem (using its connection to  oscillatory multivariate integration), the results in \cite{GrMaZh:20} are restricted to the case of 1 space dimension,
whereas the analysis in the present paper holds in general space dimension.
\eit

Thus the present paper is the first to give an analysis of all sources of error
explicitly in the (possibly large) parameter $k$ for random Helmholtz
problems arising in wave scattering with both $A$ and $n$ random, and the
first to apply this to practical scattering problems posed on an infinite domain.

Finally we mention that in \cite{GaKuSl}
the stochastic perturbation of $n$ was constrained to decay with
$\mathcal{O}(1/k)$ as $k \rightarrow \infty$ (and this condition can  play a role in our analysis
-- see \eqref{weakly_stoch} above). 
A condition similar to \eqref{weakly_stoch} (there
called ``weakly stochastic'') played a prominent r\^{o}le in the
analysis of the `multimodes' algorithm discussed in \cite{FeLiLo:15}, see
also~\cite{FeLiLo:19}.
There is some numerical evidence that \eqref{weakly_stoch} is actually
necessary to keep the QMC error bounded when $\xi$ and $s$ are fixed but
$k$ grows: a numerical experiment in \cite[\S4.6.3]{Pe:19} found that the number of QMC points $N$ needed to keep the error bounded below a constant increases quickly (apparently
exponentially) with $k$, although only a moderate range of $k$ was tested.

\subsection{The difficulty caused by trapping/nontrapping} \label{sec:trapping}

In this paper we have to deal with the fact that the $k$-dependence of the
Helmholtz solution operator depends in a complicated, nonlinear way on the
coefficients $A,n$, and the impenetrable obstacle $D$. Indeed, if $A,n$
and  $D$ support  trapped rays, then for certain data $f$ (with
$\|f\|_{L^2(D_{R_0})}=1$) and a certain unbounded sequence of wavenumbers,
the norm of the solution to \eqref{eq:intro_tedp}--\eqref{eq:src} grows
exponentially. However, if $A,n,$ and $D$ are such that no rays are
trapped (i.e., the problem is nontrapping), then the solution to
\eqref{eq:intro_tedp}--\eqref{eq:src} is bounded uniformly in $k$ (in a
suitable $k$-weighted norm); see Lemma~\ref{lem:boundsDtN}(i) below.

The theory in this paper is restricted to a certain class of $A$ and $n$
and star-shaped $D$, for which bounds on the Helmholtz solution operator
that are explicit in both $k$ and the coefficients $A$ and $n$ are
available from \cite{GrPeSp:19}. The coefficients $A$ and $n$ in this
class satisfy pointwise bounds on the fields and their derivatives; see
Definition \ref{def:nontrapping} below. It is important to emphasise that,
because the properties of trapping/nontrapping depend on the \emph{global}
properties of $A,n,$ and $D$ \emph{in combination}, any pointwise
conditions imposed on $A$ and $n$ that ensure nontrapping (such as in
Definition \ref{def:nontrapping}) are, \emph{by construction}, sufficient
but not necessary. As mentioned above, the bad behaviour due to trapping
is very specific to both the value of $k$ (see \cite{LaSpWu:21}) and the
particular data; however, establishing rigorous convergence results about
QMC methods applied to the Helmholtz equation with ``generic''
coefficients, data, and wavenumbers seems currently out of reach.

\section{The Exterior Dirichlet Problem (EDP)}
\label{subsec:EDP}

In this section we define our model wave scattering problem and give the
basic theory for its solution. The model problem involves scattering both
by the assumed heterogeneity of the medium (characterised by variable
coefficients $A$ and $n$) and also (possibly) by an impenetrable obstacle
$D$. Later we focus on the plane-wave scattering problem, where the
scattering is initiated by an incident plane wave $u^I$.

The following assumption holds throughout the paper.

\begin{assumption}[The domain and the source term] \label{ass:1}
For $d\in \{2,3\}$, let $D\subset\bbR^d$ be a bounded (possibly empty)
Lipschitz domain such that its open complement $D_+:= \mathbb{R}^d
\setminus \overline{D}$ is connected (and so the scattering problem in
$D_+$ makes sense), and denote its boundary as usual by $\partial D$.
Moreover, let $B_{R_0}$ be an open ball of radius $R_0>0$ that contains
$\overline{D}$ and also contains the supports of $A-I$ and $n-1$, and
let $D_{R_0} := B_{R_0}\setminus\overline{D}$. Furthermore, let the
support of $f$ be in $\overline{D_{R_0}}$ and let $f\in
L^2(D_{R_0})$. (Throughout, the support of a function is defined as the
closure of the set of points at which the function has non-zero value.)
\end{assumption}

\begin{notation} \label{not:oneinf}
Let $\vert \cdot \vert_2$ denote both the Euclidean norm on $\mathbb{C}^d$
and the induced matrix norm (i.e., the spectral norm) on $d \times d$
matrices. With $\Omega$ an arbitrary bounded Lipschitz domain, let $\Vert
\cdot \Vert_{L^\infty(\Omega)}$ and $\Vert v \Vert_{W^{1,\infty}(\Omega)}
:= \max(\|v\|_{L^\infty(\Omega)}, \||\nabla
v(\cdot)|_2\|_{L^\infty(\Omega)})$ be the standard $L^\infty$ and
$W^{1,\infty}$ norms on scalar-valued functions on $\Omega$. We extend
these definitions to real symmetric matrix-valued functions $A$ on
$\Omega$ by
$$
 \Vert A \Vert_{L^\infty(\Omega)} = \Vert \sigma_{\max}(A)  \Vert_{L^\infty(\Omega)}
 \quad \text{and} \quad
 \Vert A \Vert_{W^{1,\infty}(\Omega)} = \Vert \sigma_{\max}(A)   \Vert_{W^{1,\infty}(\Omega)},
$$
where $A$ is a real symmetric $d\times d$ matrix, and
$\sigma_{\max}(A)$ is its maximum singular value.
\end{notation}

\subsection{The deterministic problem}

\label{subsec:determin}

We first describe the deterministic setting, with the random case
described in \S\ref{subsec:random_prob}.

\begin{assumption}[The coefficients of the PDE] \label{ass:2}
The function $n \in W^{1,\infty}(D_+)$ is real-valued, and for almost
all $\bsx \in D_+$ there exist $n_{\min}, n_{\max}>0$ such that
\begin{align} \label{bounds_n}
 n_{\min} \leq n(\bsx) \leq n_{\max}.
\end{align}
Moreover, $A \in W^{1,\infty}(D_+)$ and for almost all $\bsx \in D_+$,
$A(\bsx)$ is real symmetric and there exist $A_{\min}, A_{\max}>0$ such
that
\begin{align}\label{bounds_A}
 A_{\min} \vert \bszeta \vert_2^2
 \leq  (A(\bsx) \bszeta) \cdot  \bszeta
 \leq A_{\max} \vert \bszeta \vert_2^2  \quad \tfa 
 \bszeta\in \mathbb{R}^d.
\end{align}
\end{assumption}

We now define the exterior Dirichlet problem from the introduction in a
more precise way.

\begin{definition}[The exterior Dirichlet Problem (EDP)]\label{def:EDP}
Given $\Dminus$, $A$, $n$, and $R_0$ satisfying Assumptions~\ref{ass:1}
and~\ref{ass:2}, wavenumber $k > 0$, and data $f \in L^2(D_{R_0})$, we say
that a function $u : D_{+} \to \CC$ is a solution of the \emph{exterior
Dirichlet problem} if, for every $R \geq R_0$, the restriction of $u$ to
the domain $D_R = B_R\setminus\overline{D}$ satisfies $u|_{D_R} \in
H^1(D_R)$, and also $u$ satisfies the PDE \eqref{eq:intro_tedp} (in the
standard weak sense, i.e., under testing with $C^\infty$ functions with
compact support; see, e.g., \cite[p.2874]{GrPeSp:19}), as well as the
boundary condition \eqref{eq:intro_dbc}, and the radiation condition
\eqref{eq:src}.
\end{definition}

Requiring $u|_{D_R}\in H^1(D_R)$ to satisfy \eqref{eq:intro_dbc} is
legitimate here, since the trace of $u$ on $\partial D$ is a well-defined
function. Moreover, requiring \eqref{eq:src} is also legitimate, since,
although $u|_{D_R}$ is only \emph{a priori} in $H^1(D_{R})$ for every
$R\geq R_0$, interior regularity results for elliptic PDEs imply that $u$
is $C^{\infty}$ outside~$B_{R_0}$.

To obtain the variational formulation of the EDP we first introduce the
following Exterior Dirichlet to Neumann map. Using this, we formulate the
EDP on a bounded domain $D_R$ for any $R \geq R_0$, in the variational
form.

\begin{definition}[Exterior Dirichlet-to-Neumann map $\DtN$ on $\partial{B}_R$]
Let $\Dminus$, $A$, $n$ and $R_0$ satisfy Assumptions~\ref{ass:1}
and~\ref{ass:2}. For any $R\ge R_0$ and $g\in H^{1/2}(\partial B_R)$, let
$v$ be the solution of $(\Delta +k^2)v=0$ in the infinite domain exterior
to $B_R$, with $v$ also satisfying the Dirichlet condition $v = g$ on
$\partial B_R$, as well as the Sommerfeld radiation condition
\eqref{eq:src}. Then $\DtN g \in H^{-1/2}(\partial B_R)$ is defined to be
the Neumann trace (with outward-pointing normal) of $v$ (recall that
the Neumann trace equals the outward-facing normal derivative, provided
that the function has sufficient regularity).
\end{definition}

Note that although $\DtN$ depends on both $k$ and $R$, our notation only
reflects its $k$-dependence.

\paragraph{Variational form of \eqref{eq:intro_tedp}--\eqref{eq:src}.}

Green's identity (see, e.g., \cite[Lemma 4.3]{Mc:00}) implies that the
restriction $u_R := u|_{D_R}$ of the solution of the EDP to $D_{R}$ for
any $R\geq R_0$ is also the solution of the variational problem:
\beq\label{eq:vp}
 \text{ find } u_R \in H_{0,D}^1(D_R )
 \text{ such that } \,a_k(u_R,v) = F(v) \,\tfa v \in H_{0,D}^1(D_R )
 \eeq
where
 \beq\label{eq:spaceEDP}
 H_{0,D}^1(D_R ):= \big\{ v\in H^1(D_R ) : v=0 \ton \Gdir\big\},
 \eeq
 \beq\label{eq:agen}
 a_k(w,v) \de \int_{D_R } \mleft(\mleft(A \grad w\mright)\cdot\grad \vbar - k^2 n\minispace w \vbar\mright)
 - \int_{\partial B_R} \DtN( w \vert_{\partial B_R})\, \overline{v}
 \quad\tand \quad
 F(v) = \int_{D_R } f \,\vbar .
 \eeq

Note that the $\DtN$ operator plays a theoretical r\^{o}le, providing us with the
variational form \eqref{eq:vp}--\eqref{eq:agen}. It is not computed explicitly, but later we compute
approximations of \eqref{eq:vp} by approximating $\DtN$ using a
perfectly-matched layer; see \S\ref{sec:FEM} below.

Under the conditions of Definition~\ref{def:EDP}, the EDP has a unique
solution $u$ and $\Vert u \Vert_{H^1_k(D_R)}$ is finite for all $R \geq
R_0$ (see \cite[Theorem 2.5]{GrPeSp:19}). Moreover the variational form
\eqref{eq:vp}--\eqref{eq:agen} also has a unique solution $u_R\in
H_{0,D}^1(D_R)$, with $\Vert u \Vert_{H^1_k(D_R)}$ also finite, and so
$u_R$ coincides with the restriction of $u$ to $D_R$. Hence bounds on the
solution of the variational form imply bounds on the solution to the EDP.

However the variational form permits the inclusion of more general
functionals $F$ in \eqref{eq:vp} than the special $F$ appearing in
\eqref{eq:agen}. Indeed, any $F\in H_{0,D}^1(D_R)'$ (the dual of
$H_{0,D}^1(D_R)$) is allowed in \eqref{eq:vp}. Such more-general
right-hand sides are needed in some of our arguments see, for example
\eqref{varderiv} and~\eqref{Ftilde}.

\paragraph{Bounds on the solution of the EDP explicit in $k$, $A$, and $n$.}

The $k$-dependence of the bounds on the solution of the EDP in terms of
the data depends crucially on whether or not $A$, $n$, and $\Dminus$ allow
trapped rays. In the nontrapping case, the solution operator is bounded
uniformly in $k$ (in suitable norms), whereas
in the trapping case the solution operator can grow exponentially in $k$;
see the reviews in \cite[Section 6]{MoSp:19},
\cite[Section~1]{GrPeSp:19}, \cite[\S1]{LaSpWu:21} and the references
therein.

We now discuss conditions on $A$ and $n$ such that \emph{both}
the problem is nontrapping (recall the discussion in
\S\ref{sec:main_results}) \emph{and} that we know explicitly how the
constants in the bounds on the solution of the EDP in terms of the data
depend on $A$ and $n$; this latter condition is important when we consider
random coefficients in the next subsection.

\bde[A particular class of nontrapping $D$, $A$, and $n$]
\label{def:nontrapping} Let $D$, $A$, and $n$ satisfy
Assumptions~\ref{ass:1} and~\ref{ass:2}. Furthermore let $D$ be
star-shaped with respect to the origin. Then, given $\mu_A,\mu_n>0$ we say
that $A \in \NTA{\mu_A}$ and $n \in \NTn{\mu_n}$ if
\begin{align}\label{eq:hh-Acond}
 \bszeta^\top (A(\bsx) - \mleft(\bsx\cdot\grad\mright)A(\bsx))\, \bszeta \
 &\geq \ \mu_{A} ,
 \quad \text{for all} \quad \bszeta \in \mathbb{R}^d
 \quad \text{with} \quad  \vert \bszeta \vert_{2} = 1, \quad \tand \nonumber\\
 n(\bsx) + \bsx\cdot\grad n(\bsx)\ &\geq\  \mu_n
\end{align}
for almost every $\bsx \in D_+$. \ede

The relationship between Definition~\ref{def:nontrapping} and the
definition of nontrapping in terms of rays is explained in
\cite[\S7]{GrPeSp:19}. Under the condition~\eqref{eq:hh-Acond},
\cite{GrPeSp:19} gives results that show explicitly how the constant in
the bound of the solution of the EDP in terms of the data $f$ depends on
$\mu_A$, $\mu_n$. The more recent results in \cite{GaSpWu:20} generalise
those in \cite{GrPeSp:19}, but require more smoothness of $A,n,$ and
$\Gdir$ than assumed here. We now state these results applied to the
particular situation of $D,\Gdir, A,$ and $n$ satisfying
Definition~\ref{def:nontrapping}.

\ble[Bounds on the solution of the EDP]\label{lem:boundsDtN} Let $D, A,$
$n$, and $R_0$ be as in Definition~\ref{def:nontrapping}. Given $k_0>0$
let
\[
  C_1
  :=2 \sqrt{ \frac{1}{\mu_A} + \frac{1}{\mu_n}
  \bigg( 1 + \frac{d-1}{2\,k_0\, R_0}\bigg)^2}.
\]
Then, for all $k \geq k_0$ and $R\geq R_0$, the following estimates hold.
\begin{enumerate}
\item [\textnormal{(i)}]%
Given $f\in L^2(D_{R_0})$, the solution of the EDP of Definition
\ref{def:EDP} satisfies \beq\label{eq:EDPbound1} \N{u}_{H^1_k(D_R)}
\leq \frac{C_1}{\sqrt{\min\{\mu_A,\mu_n\}}}\, R\,
 \N{f}_{L^2(D_{R_0})},
  \eeq
\item [\textnormal{(ii)}]%
Given $F \in H_{0,D}^1(D_R )'$, the solution of the variational
problem \eqref{eq:vp} satisfies
\begin{align}\label{eq:dualbound}
 \N{u}_{H^1_k(D_R)}
 \ \leq \ \frac{1}{\min\{ A_{\min}, n_{\min}\}}
 \bigg( \frac{1}{k_0\,R_0} + 2
 \frac{C_1}{\sqrt{\min\{\mu_A,\mu_n\}}}\, n_{\max} \bigg)\,k\,R\,
 \N{F}_{H^1_k(D_R )'},
\end{align}
where $\N{F}_{H^1_k(D_R)'}$ denotes the norm of the functional $F$
with respect to the $H^1_k$ norm on $D_R$.
\end{enumerate}
\ele

\bpf Part (i) is \cite[Theorem~2.5]{GrPeSp:19} and Part (ii) is
\cite[Corollary~2.16]{GrPeSp:19}.
\epf

\subsection{The plane-wave sound-soft scattering problem}
\label{subsec:SSSP}

An important scattering problem in applications is the plane-wave
sound-soft scattering problem.

\begin{definition}[Plane-wave sound-soft scattering problem]\label{def:SSSP}
Given $D$, $A$, $n$, and $R_0$ satisfying Assumptions~\ref{ass:1}
and~\ref{ass:2}, wavenumber $k
> 0$, and unit vector $\bsalpha\in\bbR^d$, $\vert \bsalpha \vert_{2} = 1$, let
\begin{align} \label{eq:incident}
 u^I(\bsx) := \exp(\ri k \bsx\cdot\bsalpha), \quad \bsx\in\bbR^d,
\end{align}
be an incident plane wave in the direction $\bsalpha$. The total field
$u^T : D_{+} \to \CC$ is the solution of the plane-wave sound-soft
scattering problem if, for every $R\ge R_0$, the restriction of $u^T$ on
$D_R$ satisfies $u^T|_{D_R} \in H^1(D_R)$, and $u^T$ satisfies
\begin{align}
 \grad \cdot \mleft(A \grad u^T\mright) + k^2 n u^T = 0
 \quad\tin D_+, \qquad u^T = 0 \quad\ton \Gdir.
 \label{eq:PWSS}
\end{align}
In addition, the scattered field $u^S:=u^T-u^I|_{D_+}$ satisfies
$u^S|_{D_R} \in H^1(D_R)$ for every $R\ge R_0$ and it satisfies the
Sommerfeld radiation condition \eqref{eq:src} (with $u$ replaced by
$u^S$).
\end{definition}

A natural object of interest  for the plane-wave sound-soft scattering
problem is the far-field pattern of the scattered field $u^S$ defined by
 \beq
 u^{S}_\infty(\widehat{\bsx}) :=
 \lim_{|\bsx|_2\rightarrow \infty} \Big( u^S(\widehat{\bsx}\,|\bsx|_2)
 \,\exp(-\ri k |\bsx|_2)\, |\bsx|_2^{(d-1)/2} \Big).
 \label{eq:FFP} \eeq
We will return to this in \S\ref{sec:SSSP}.

\subsection{The random problem}
\label{subsec:random_prob}

We now study uncertainty quantification for the EDP with random
coefficients in the stochastic domain $U :=
[-\frac{1}{2},\frac{1}{2}]^\bbN = \{ \bsy = (y_1, y_2, y_3, \ldots ): y_j
\in [-\frac{1}{2},\frac{1}{2}]$\}. The coefficients are required to satisfy the following
assumptions.

\begin{assumption}[Random coefficients] \label{ass:3}
The coefficients $A(\bsx, \bsy), n(\bsx,\bsy)$ are given by the affine series
\eqref{eq:Axy} and \eqref{eq:nxy} under the following conditions.
\begin{itemize}
\item The mean fields  $A_0, n_0$ satisfy Assumption \ref{ass:2} with
    lower and upper bounds (analogous to  \eqref{bounds_n} and
    \eqref{bounds_A}), denoted $A_{0,\min}, A_{0,\max}, n_{0,\min},
    n_{0,\max}$.
\item The obstacle $D$ and the mean fields $A_0, n_0$ satisfy the
    nontrapping condition in Definition~\ref{def:nontrapping}, i.e.,
    $D$ is star-shaped with respect to the origin, and there exist
    $\mu_A$, $\mu_n$ such that $A_0 \in\NTA{\mu_A}$ and $n_0 \in
    \NTn{\mu_n}$.
   \item The supports of $A_0 - I$, $n_0 - 1$, $\Psi_j$ and $\psi_j$
       for all $j\geq 1$ are all compactly contained in $B_{R_0}$.

\item The scaling parameters are nonnegative: $\xi_A\ge 0$ and
    $\xi_n \ge 0$.

\item The fluctuations $\Psi_j \in W^{1,\infty}(D_{R_0})$ are real
    symmetric matrices with (recalling Notation~\ref{not:oneinf}),
    $\|\Psi_1\|_{L^\infty(D_{R_0})}=1$ and \beq\label{eq:Apsimeaspos}
    \xi_A \sum_{j=1}^\infty \N{\Psi_j}_{L^{\infty}(D_{R_0})}
    \le 
    A_{0,\min} \quad\tand\quad\sum_{j=1}^\infty
    \N{\Psi_j}_{W^{1,\infty}(D_{R_0})} < \infty.
     \eeq

\item The fluctuations $\psi_j \in W^{1,\infty}(D_{R_0})$ satisfy
    $\|\psi_1\|_{L^\infty(D_{R_0})}=1$, \beq\label{eq:npsimeaspos}
    \xi_n \sum_{j=1}^\infty \N{\psi_j}_{L^\infty(D_{R_0})}
    \le
    n_{0,\min} \quad\tand\quad\sum_{j=1}^\infty
    \N{\psi_j}_{W^{1,\infty}(D_{R_0})} < \infty.
     \eeq

\item There exists a summability exponent $p \in (0,1)$ for the
    sequence $\bsb = (b_j)_{j\ge 1}$ such that
\begin{align} \label{eq:bj}
 \sum_{j=1}^\infty b_j^p \ <\  \infty ,
 \quad \text{where} \quad
 b_j \,:=\, \max\big(\|\psi_j\|_{L^\infty(D_{R_0})},\|\Psi_j \|_{L^{\infty}(D_{R_0})}\big).
\end{align}

\end{itemize}
\end{assumption}

The following lemma is an immediate consequence of Assumption \ref{ass:3} and the triangle inequality.

\begin{lemma}\label{lem:easy}
Under Assumption~\ref{ass:3}, $A(\bsx,\bsy)$ and $n(\bsx,\bsy)$ defined by
\eqref{eq:Axy} and \eqref{eq:nxy} satisfy Assumption~\ref{ass:2} uniformly
in $\bsy\in U$, with
\begin{align*}
 A_{\min} & := \frac{A_{0,\min}}{2} \quad  \text{and} \quad A_{\max} :=  A_{0,\max}+\frac{A_{0,\min}}{2}, \\
 n_{\min} & := \frac{n_{0,\min}}{2} \quad \text{and} \quad n_{\max} :=  n_{0,\max} + \frac{n_{0,\min}}{2}.
\end{align*}
\end{lemma}

\begin{remark} \label{exis_unique}
Only the first estimates in \eqref{eq:Apsimeaspos} and
\eqref{eq:npsimeaspos} are needed to prove Lemma~\ref{lem:easy}. The
second estimates were used in \cite[Appendix C]{PeSp:20a} to prove that
$A$ and $n$ are measurable. Under Assumption \ref{ass:3}, existence and
uniqueness of the solution of the EDP of Definition \ref{def:EDP} is
proved in \cite[Theorem 1.4, \S1.2, Remark 1.12]{PeSp:20}.
\end{remark}

We now impose conditions on $A-A_0$ and $n-n_0$ so that the coefficients
$A$ and $n$ belong to the class of coefficients in Definition
\ref{def:nontrapping}, and are thus nontrapping for all realisations of
$\bsy\in U$. As in \S\ref{sec:main_results}, we emphasise that, because
the properties of trapping/nontrapping depend on the global properties of
$A,n,$ and $D$ in combination, any pointwise conditions imposed on $A$ and
$n$ that ensure nontrapping (such as in Definition \ref{def:nontrapping})
will be sufficient but not necessary.

\begin{assumption}[Conditions on the fluctuations ensuring nontrapping]\label{ass:4}
Suppose that Assumption~\ref{ass:3} holds. Then assume further that
\begin{align} \label{eq:nontrap}
 \xi_A
 \sum_{j=1}^\infty 
 \esssup_{\bsx \in D_{R_0}}
  | \Psi_j(\bsx) -
  \mleft(\bsx\cdot\grad\mright)\Psi_j(\bsx)|_2 & \leq \mu_A, 
  \quad\mbox{and} \nonumber\\
 \xi_n
 \sum_{j=1}^\infty
 \esssup_{\bsx \in D_{R_0}} |\psi_j(\bsx) + \bsx\cdot\grad\psi_j(\bsx) |
 &\leq \mu_n, 
\end{align}
where $\mu_A$, $\mu_n$ are as given in Assumption~\ref{ass:3}.
\end{assumption}

The fact that the corresponding EDP is nontrapping is guaranteed by the
following corollary, whose proof follows from Lemma~\ref{lem:boundsDtN},
Definition~\ref{def:nontrapping} and the triangle inequality. To
simplify the statement, we first replace the constant in
\eqref{eq:EDPbound1} by the larger constant in \eqref{eq:dualbound}.

\begin{corollary} \label{cor:random}
Let $D$, $A$, $n$, and $R_0$ satisfy Assumptions~\ref{ass:1},
\ref{ass:3}~and~\ref{ass:4}, and let $k_0>0$. Then we have $A \in
\NTA{\mu_A/2}$ and $n \in \NTn{\mu_n/2}$ for all $\bsy \in U$. Consequently, the bounds \eqref{eq:EDPbound1} and \eqref{eq:dualbound} (with $\mu_A$, $\mu_n$ replaced by $\mu_A/2$,
$\mu_n/2$)
hold for the solution $u = u(\cdot, \bsy)$ of the EDP with random
coefficients $A(\cdot, \bsy)$, $n(\cdot, \bsy)$.

In particular, introducing the stability constant
\begin{align} \label{eq:Cstab}
\Cstab:=
 \frac{1}{\min\{ A_{\min},
    n_{\min}\}} \bigg( \frac{1}{k_0\, R_0} +
    4 \sqrt{ \frac{1}{\mu_A/2} + \frac{1}{\mu_n/2}
    \bigg( 1 + \frac{d-1}{2\,k_0\, R_0}\bigg)^2}
    \frac{n_{\max}}{\sqrt{\min\{\mu_A/2, \mu_n/2\}}}
    \bigg),
\end{align}
then we have, for all $k\ge k_0$, $R\ge R_0$, and $\bsy \in U$:
\begin{itemize}
\item[\textnormal{(i)}]%
Given $f\in L^2(D_{R_0})$, the solution of the EDP of
Definition \ref{def:EDP} satisfies
\begin{align} \label{eq:EDP1}
 \N{u(\cdot,\bsy)}_{H^1_k(D_R)} \ \leq \
\Cstab \,   R\,
 \N{f}_{L^2(D_{R_0})}.
\end{align}

\item[\textnormal{(ii)}]%
Given $F\in H_{0,D}^1(D_R)'$, the solution of the variational problem
\eqref{eq:vp} satisfies
\begin{align} \label{eq:EDP2}
 \Vert u(\cdot, \bsy) \Vert_{H^1_k(D_R)} \ \ \leq \ \Cstab\,k\,R\,   \Vert F \Vert_{H_k^1(D_R)'}. 
\end{align}
       \end{itemize}
\end{corollary}

\section{Parametric regularity} \label{sec:regu}

In this subsection we estimate the partial derivatives of $u(\cdot, \bsy)$
with respect to $\bsy$. These estimates are a crucial ingredient in the analysis of the dimension truncation error in \S\ref{sec:dimtrunc} and the QMC error \S\ref{sec:QMC}. In particular,  Corollary~\ref{cor:regu} is used directly in the derivation of \eqref{reg_est}, which in turn is used in \eqref{eq:taylor3}, \eqref{eq:taylor5}, \eqref{eq:Gu-norm} and \eqref{eq:error2}.

Let $\bbN_0 = \bbN \cup \{0\}$ denote the set of
nonnegative integers. For any multi-index $\bsnu\in\bbN_0^\infty$ with
finite order $|\bsnu| := \sum_{j\ge 1} \nu_j < \infty$, let
\[
  \partial^\bsnu \,=\, \partial^\bsnu_\bsy \,=\, \prod_{j\ge 1} \bigg(\frac{\partial}{\partial y_j}\bigg)^{\nu_j}
\]
denote the mixed partial derivative with respect to $\bsy$.

\begin{theorem} \label{thm:regu}
Suppose that $u(\cdot,\bsy)$ is the solution of the EDP with $D$, $A$,
$n$, and $R_0$ satisfying Assumptions~\ref{ass:1},
\ref{ass:3}~and~\ref{ass:4}, and let $k_0>0$. Then, for any
multi-index $\bsnu\in\bbN_0^\infty$ with $|\bsnu|<\infty$, all $k\ge
k_0$, $R \geq R_0$ and all $\bsy\in U$,
\begin{align} \label{eq:regu}
  \|\partial^\bsnu u(\cdot,\bsy)\|_{H^1_k(D_R)} \,\le\,
  (\Cstab \,k R\,\xi)^{|\bsnu|}\, |\bsnu|!\, \bsb^\bsnu\,
  \| u(\cdot, \bsy )\|_{H^1_k(D_R)},
\end{align}
with $\xi = \xi_n + \xi_A$, and $\bsb^\bsnu := \prod_{j\ge 1}
b_j^{\nu_j}$, where $b_j$ is defined in \eqref{eq:bj} and $\Cstab$ is
defined in \eqref{eq:Cstab}.
\end{theorem}

\begin{proof}
We prove this result by induction on $|\bsnu|$. The result holds trivially
for $\bsnu=\bszero$. For $\bsnu\ne \bszero$, we first note that, with
$a_k$ and $F$ as in \eqref{eq:agen}, the solution $u(\cdot, \bsy) \in
\Honezero$ satisfies the variational formulation
\begin{align} \label{var}
  a_k (u,v) = F(v) \quad \text{for all} \quad v \in \Honezero .
\end{align}
To obtain the result, we use the following consequence of the Leibnitz
rule.
\begin{align} \label{Leibnitz}
 \partial^{\bsnu} (A \nabla u)
 = \sum_{\boldsymbol{0}\le \bsm\le\bsnu} \binom{\bsnu}{\bsm}
 (\partial^\bsm A) \, \nabla  (\partial^{\bsnu-\bsm} u),
 \quad
 \partial^{\bsnu} (n u)
 = \sum_{\boldsymbol{0} \le \bsm\le\bsnu} \binom{\bsnu}{\bsm}
 (\partial^\bsm n) \, (\partial^{\bsnu-\bsm} u);
\end{align}
here $\boldsymbol{0} \le \bsm\le\bsnu$ means that $0 \le m_j\le \nu_j$ for
all $j\ge 1$, and $\binom{\bsnu}{\bsm} = \prod_{j\ge 1}
\binom{\nu_j}{m_j}$.

Noting the linearity of $n(\bsx,\bsy)$ with respect to $\bsy$ in
\eqref{eq:nxy}, it is easy to see that if we differentiate $n(\bsx,\bsy)$
with respect to $y_j$, we obtain $\xi_n\,\psi_j(\bsx)$, and if we
differentiate a second time with respect to any variable we get $0$. The
derivatives of $A$ behave analogously. Thus we have
\begin{align} \label{diff}
  (\partial^{\bsm}n)(\bsx,\bsy) \,=\,
 \begin{cases}
 n(\bsx,\bsy) &\mbox{if } \bsm=\bszero, \\
 \xi_n\, \psi_j(\bsx) &\mbox{if } \bsm=\bse_j, \\
 0 &\mbox{otherwise},
 \end{cases}
 \quad\mbox{and}\quad
 (\partial^{\bsm}A)(\bsx,\bsy) \,=\,
 \begin{cases}
 A(\bsx,\bsy) &\mbox{if } \bsm=\bszero, \\
 \xi_A\, \Psi_j(\bsx) &\mbox{if } \bsm=\bse_j, \\
 0 &\mbox{otherwise},
 \end{cases}
\end{align}
where $\bse_j$ denotes the multi-index whose $j$th component is $1$ and
all other components are $0$.

We now assume the inductive hypothesis, that the estimate \eqref{eq:regu}
holds for all $\bsnu$ of order $\vert \bsnu \vert \le M$, and then choose
any $\bsnu$ with $\vert \bsnu\vert = M+1$. Then, applying
$\partial^{\bsnu}$ to \eqref{var}, the right-hand side vanishes (since it
is independent of $\bsy$). Applying $\partial^\bsnu$ to   the left-hand
side, inserting the formulae \eqref{Leibnitz} and using \eqref{diff}, and
then finally moving the terms corresponding to $\bsm\ne\bszero$ to the
right-hand side, we obtain
\begin{align} \label{varderiv}
 a_k(\partial^{\bsnu} u,v) = F_\bsnu(v) \quad \text{for all } \quad v \in \Honezero,
\end{align}
with
\begin{align} \label{Ftilde}
  F_\bsnu(v) := - \xi_A
  \sum_{j\ge 1} \nu_j \int_{D_R}
  (\Psi_j \nabla \partial^{\bsnu - \bse_j} u ) \cdot \nabla \overline{v}
  \;+\; k^2\, \xi_n \sum_{j\ge 1} \nu_j
  \int_{D_R} \psi_j \, \partial^{\bsnu - \bse_j} u \, \overline{v}.
\end{align}
Due to linearity of the $\DtN$ map, the boundary integral in
\eqref{eq:agen} becomes $\int_{\partial B_R} \DtN(\partial^{\bsnu}
u|_{\partial B_R})\,\overline{v}$ which is part of the left-hand side
$a_k(\partial^{\bsnu} u,v)$, while no boundary term appears in the
right-hand side $F_\bsnu(v)$. From Corollary~\ref{cor:random}(ii) we then
conclude that
\begin{align} \label{eq:almost}
  \|\partial^{\bsnu} u(\cdot, \bsy)\|_{H^1_k(D_R)}
  \le \Cstab\, k\, R\, \|F_\bsnu\|_{H^1_k(D_R)'}.
\end{align}

We now  estimate the norm on the right-hand side of \eqref{eq:almost}.   
For the first sum on the right-hand side of \eqref{Ftilde}, we apply Cauchy--Schwarz
inequality and then the inductive hypothesis to obtain (recalling Notation
\ref{not:oneinf}):
\begin{align*}
 &\bigg| \xi_A\, \sum_{j\ge 1} \nu_j\,
 \int_{D_R} (\Psi_j \nabla \partial^{\bsnu - \bse_j} u)
 \cdot \nabla \overline{v} \bigg|
 \,\le\, \xi_A\, \sum_{j\ge 1} \nu_j\,  \|\Psi_j\|_{L^{\infty}(D_R)}\,
 \|\partial^{\bsnu - \bse_j} u\|_{H^1_k(D_R)}\, \|v\|_{H^1_k(D_R)} \nonumber\\
 &\,\le\, \xi_A\, \sum_{j\ge 1} \nu_j\, \|\Psi_j\|_{L^{\infty}(D_R)}\,
    (\Cstab\, k\, R\, \xi)^{|\bsnu - \bse_j|} \, |\bsnu-\bse_j|! \, \bsb^{\bsnu-\bse_j} \,
   \| u(\cdot, \bsy )\|_{H^1_k(D_R)}\,
   \|v\|_{H^1_k(D_R)} \nonumber\\
 &\,\le\, \xi_A \,  (\Cstab\, k\,R\, \xi)^{|\bsnu| -1}\,|\bsnu|! \,
  \bsb^{\bsnu}
  \| u(\cdot, \bsy )\|_{H^1_k(D_R)}\,
  \|v \|_{H^1_k(D_R)},
\end{align*}
where we used $\|\Psi_j\|_{L^{\infty}(D_R)} \leq b_j$, $|\bsnu-\bse_j| =
|\bsnu|-1$, and $\sum_{j\ge1} \nu_j = |\bsnu|$.

We can apply an analogous process to estimate the second term on the
right-hand side of \eqref{Ftilde}, but in this case we apply
Cauchy--Schwarz inequality in $L^2(D)$. Then the $k^2$ factor appearing
there is cancelled out by the fact that (recalling \eqref{eq:weighted}),
$\|w\|_{L^2(D_R)} \le k^{-1} \|w\|_{H^1_k(D_R)}$ for any $w \in H^1(D_R)$.
Putting these estimates together and using $\xi = \xi_n + \xi_A$, we
obtain
$$
 \|F_\bsnu \|_{H^1_k(D_R)'}
 \leq \xi \,  (\Cstab\, k\, R\, \xi)^{|\bsnu|-1}\,
 \vert \bsnu \vert ! \, \bsb^{\bsnu}\,
 \| u(\cdot, \bsy )\|_{H^1_k(D_R)}.
$$
Substituting this into \eqref{eq:almost} yields the desired result.
\end{proof}

\begin{corollary}\label{cor:regu}
Suppose that $u(\cdot,\bsy)$ is the solution of the EDP with $f\in
L^2(D_{R_0})$ and $D$, $A$, $n$, and $R_0$ satisfying
Assumptions~\ref{ass:1}, \ref{ass:3}~and~\ref{ass:4}, and let
$k_0>0$. Then, for all $k\geq k_0, R\geq R_0$, $r \in [0,1]$, and
$\bsy\in U$,
\begin{align} \label{eq:regu_cor2}
  \|\partial^\bsnu u(\cdot,\bsy)\|_{H^r(D_R)} \,\le\,
  (\Cstab \,k\, R\,\xi)^{|\bsnu|}\, |\bsnu|!\, \bsb^\bsnu\,
  \Cstab\, R\, k^{r-1}\,  \N{f}_{L^2(D_{R_0})}.
\end{align}
where $b_j$ is defined in \eqref{eq:bj} and $\Cstab$ is defined in
\eqref{eq:Cstab}.
\end{corollary}

\begin{proof}
Combining Theorem~\ref{thm:regu} with the bound given by
Corollary~\ref{cor:random}(i), we obtain
\begin{align}\label{eq:const_frac}
 \|\partial^\bsnu u(\cdot,\bsy)\|_{H^1_k(D_R)} \,\le\,
 (\Cstab \,k\, R\,\xi)^{|\bsnu|}\, |\bsnu|!\, \bsb^\bsnu\,
 \Cstab\, R\, \N{f}_{L^2(D_{R_0})} .
\end{align}
The estimates \eqref{eq:regu_cor2} for $r =0$ and $r = 1$ then follow
on recalling the definition of the $H^1_k$-norm in \eqref{eq:weighted}.
The case of general $r \in (0,1)$ is obtained by applying operator
interpolation (see, e.g., \cite[\S 14]{BrSc:08}) to the linear operator
  $T : L^2(D_{R_0}) \rightarrow H^r(D_R)$ defined by
  $T f := \partial^{\bsnu} u$, where $u$ is the solution of the EDP.
\end{proof}

\section{Dimension truncation}
\label{sec:dimtrunc}

In this section we estimate the quantity $\Etrunc$ defined in \eqref{eq:tentative1}. For $G \in H^r(D_R)'$ and $r \in [0,1]$ we define
$$
 \calG(\bsy) = G(u(\cdot, \bsy)|_{D_R}), \quad \bsy \in U ,
$$
where $u$ is the solution to the EDP from Definition \ref{def:EDP}. Also
for any finite $s\in \mathbb{N}$, we define the truncated functional
$$
 \cG_s(\bsy_s) = G(u(\cdot, (\bsy_s; \bszero))|_{D_R}), \quad
 \text{where} \quad \bsy_s = (y_1, \ldots y_s).
$$
Then
$$
 \Etrunc = \left\vert \, \int_U(\calG(\bsy) - \calG_s(\bsy_s) \, \rd \bsy \, \right\vert.
$$

\begin{theorem} \label{thm:dim_trunc}
Let $D$, $A$, $n$, and $R_0$ satisfy Assumptions~\ref{ass:1},
\ref{ass:3}~and~\ref{ass:4}, and let $k_0>0$. Let $\bsb = (b_j)_{j\ge
1}$ be the sequence defined in \eqref{eq:bj}. Let $f\in L^2(D_{R_0})$ and
let $G \in H^r(D_R)'$ for some $r \in [0,1]$. Then for any finite
$s\in\bbN$, for all integers $M \geq 0$, and all $k\ge k_0$, $R \geq
R_0$, there is a constant $C(M, \bsb)$ (depending on $M$ and $\bsb$, but
independent of other parameters) such that
\begin{align}  \label{trunc_err1}
 &\Etrunc \leq   C(M, \bsb) \, \Cstab\,R\,k^{r-1}\,K_r(f,G) \nonumber \\
 &\quad\times \Big[ \Lambda_M\,
 (\Cstab\,k\,R\,\xi)^{2}\max\big\{ 1, (\Cstab\,k\,R\,\xi)^{M-2} \big\} \,s^{-2/p+1}
 + (\Cstab\,k\,R\,\xi)^{M+1}  s^{(-1/p+1)(M+1)} \Big],
\end{align}
where $\Cstab$ is defined in \eqref{eq:Cstab},
\begin{align}\label{eq:def-K}
 K_r(f,G) := \Vert f \Vert_{L^2(D_{R_0})}  \,  \Vert G \Vert_{H^r(D_R)'}
\end{align}
and
\[
 \Lambda_M = \begin{cases}
 0 & \text{when } M= 0,1, \\
 1 & \text{otherwise}.
 \end{cases}
\]
\end{theorem}

\begin{remark}
The choices $M = 0,1$ respectively yield the bounds of order:
$$
 \Cstab\, R\, k^{r-1} (\Cstab\,k\,R\,\xi) s^{-1/p+1} \quad \text{and} \quad
 \Cstab\, R\, k^{r-1} (\Cstab\,k\,R\,\xi)^2 s^{-2/p+2}
$$
(omitting constants independent of $\Cstab, k, R, \xi, s$).
For $M \geq 2$ we obtain estimates with two terms which can be balanced in
various ways, either to ensure the fastest decay with $s$ or to ensure the
two terms are equal in magnitude. For example, the choice $M = \lceil
\frac{1}{1-p}\rceil$ gives the optimal rate as $s \rightarrow \infty$ of
order
$$
 \Cstab\, R\, k^{r-1} (\Cstab\,k\,R\,\xi)^{\lceil \frac{1}{1-p}\rceil+1}  s^{-2/p+1}.
$$
\end{remark}

\begin{proof}
We give the proof for $r \in \{0,1\}$; and then extend to  $r \in
(0,1)$ using  operator interpolation. The proof uses an extension of
arguments found in \cite[Theorem~4.1]{GiGrKuScSl:19}
to the Helmholtz case, with the emphasis here on obtaining estimates
explicit in the wavenumber $k$ as well as the dimension truncation
parameter $s$. (See also \cite{GuKa:22} for further extensions and
generalisations of these arguments.)

Using the linearity of $G$ together with Corollary~\ref{cor:regu} and
\eqref{eq:def-K}, we have for any multi-index $\bsnu$ with $|\bsnu| <
\infty$ and $\bsy\in U$,
\begin{align} \label{reg_est}
 \vert (\partial^{\bsnu} \cG)(\bsy) \vert
 = \vert G(\partial^{\bsnu} u(\cdot, \bsy)|_{D_R})\vert
 & \le \|G\|_{H^r(D_R)'}\,
 \|\partial^{\bsnu} u(\cdot, \bsy)\|_{H^r(D_R)}
 \nonumber\\
 & \le (\Cstab \,k\, R\,\xi)^{|\bsnu|}\, |\bsnu|!\, \bsb^\bsnu\,
   \Cstab\,R\,k^{r-1}\,K_r(f,G).
\end{align}
With this preliminary, the next step in the proof is to write $\cG(\bsy)$
in terms of its Taylor polynomial of order $M$ about the point
$(\bsy_s;\bszero)$,  plus remainder,  yielding
\begin{align*} 
 &\int_U (\cG(\bsy)-\cG_s(\bsy_{ s}))\,{\rm d}\bsy \nonumber \\
 &=
 \underbrace{\sum_{\ell=1}^M \sum_{\substack{|\bsnu|=\ell\\ \nu_j=0~\forall j\leq s}}\frac{1}{\bsnu!}
 \int_U \bsy^{\bsnu}\,\partial^{\bsnu}\cG(\bsy_{ s};\mathbf 0)\,{\rm d}\bsy}_{=:\,T_1}
 +
 \underbrace{\sum_{\substack{|\bsnu|=M+1\\ \nu_j=0~\forall j\leq s}}\frac{M+1}{\bsnu!}
 \int_U\int_0^1(1-t)^M\,\bsy^{\bsnu}\,\partial^{\bsnu}\cG(\bsy_{s};t\,\bsy_{s+})\,{\rm d}t\,{\rm d}\bsy}_{=:\, T_2},
\end{align*}
where $\bsy_{s+}: = (y_{s+1}, y_{s+2}, \ldots)$, and so
\begin{align} \label{eq:taylor2}
   \Etrunc \leq \vert T_1\vert \ + \ \vert T_2 \vert.
\end{align}
We estimate each of the terms $T_1, T_2$ separately.

The $\ell=1$ term in $T_1$ vanishes, since when $\vert \bsnu \vert
=1$ and $\nu_j = 0$ for all $j \leq s$ we have
\begin{align} \label{gantner1}
 \int_U \bsy^{\bsnu}\,\partial^{\bsnu}\cG(\bsy_{s};\mathbf 0)\,{\rm d}\bsy
 = \bigg(\int_{U_{s}} \partial^{\bsnu}\cG(\bsy_{s};\mathbf 0)\,\,{\rm d}\bsy_{s}\bigg)
 \bigg(\prod_{j>s}\int_{-1/2}^{1/2}y_j^{\nu_j}\,{\rm d}y_j\bigg),
\end{align}
where $U_s := [-\frac{1}{2},\frac{1}{2}]^s$ as before and the product on the right-hand side vanishes because 
$\nu_j=1$ for one~$j$. This fact was first observed by \cite{Gan18}. Thus we can
start the summation in $T_1$ from $\ell=2$, and trivially $T_1 = 0$ if
$M=0$ or $M=1$. For $M\ge 2$, we use the facts that (inside the double
sum), $\vert y_j \vert \leq 1/2$, $\bsnu ! \geq 1 $ and $\vert \bsnu \vert
\geq 2$, thus obtaining
\begin{align}
 \vert T_1 \vert \
 &\leq  \sum_{\ell = 2}^M \sum_{\substack{|\bsnu|=\ell\\ \nu_j=0~\forall j\leq s}}
 \frac{1}{2^{\vert \bsnu\vert }\bsnu!} \, \sup_{\bsy \in U}
 \vert (\partial^{\bsnu} \cG)(\bsy) \vert
 \leq   \frac{1}{4} \, \sum_{\ell = 2}^M
 \sum_{\substack{|\bsnu|=\ell\\ \nu_j=0~\forall j\leq s}} \, \sup_{\bsy \in U}
 \vert (\partial^{\bsnu} \cG)(\bsy) \vert \nonumber \\
 & \leq \frac{1}{4} M!\,
 \max \big\{ (\Cstab\,k\,R\,\xi)^2, (\Cstab\,k\,R\,\xi)^M\big\}\, \bigg( \sum_{\ell = 2}^M
 \sum_{\substack{|\bsnu|=\ell\\ \nu_j=0~\forall j\leq s}} \bsb^{\bsnu} \bigg)\,
 \Cstab\,R\,k^{r-1}\, K_r(f,G),
 \label{eq:taylor3}
\end{align}
where in the second step we used \eqref{reg_est}. The double sum on the
right-hand side of \eqref{eq:taylor3} is estimated in \cite[eq
(4.10)]{GiGrKuScSl:19} by
\begin{align} \label{eq:taylor4}
  \sum_{\ell = 2}^M \sum_{\substack{|\bsnu|=\ell\\ \nu_j=0~\forall j\leq s}}
 \bsb^{\bsnu}  \leq C s^{-2/p+1},
\end{align}
where $C$ depends on $\bsb$ but not $s$. Combining \eqref{eq:taylor3} and
\eqref{eq:taylor4}, we obtain a bound for $|T_1|$ in~\eqref{eq:taylor2},
giving the first term in the bound \eqref{trunc_err1}.

Also, using similar arguments, $T_2$ can be estimated by
\begin{align}
 \vert T_2 \vert \
 &\leq \sum_{\substack{|\bsnu|=M+1 \\ \nu_j=0~\forall j\leq s}}
 \frac{M+1}{2^{\vert \bsnu\vert }\bsnu!} \,  \left(\int_0^1 (1-t)^M \rd t \right) \,
 \sup_{\bsy \in U} \vert (\partial^{\bsnu} \cG)(\bsy) \vert \nonumber\\
 & \leq (M+1)! \left(\frac{\Cstab\,k\,R\,\xi}{2} \right)^{M+1}
 \bigg(\sum_{\vert \bsnu \vert = M+1} \bsb^{\bsnu} \bigg) \, \Cstab\,R\,k^{r-1}\,K_r(f,G),
 \label{eq:taylor5}
\end{align}
and again we used \eqref{reg_est} in the last step.
The sum appearing in \eqref{eq:taylor5} is also estimated in \cite[eq
(4.11)]{GiGrKuScSl:19}, with the result
\begin{align} \label{eq:taylor6}
 \sum_{\vert \bsnu \vert = M+1} \bsb^{\bsnu}
 \leq C s^{(-1/p+1)(M+1)},
\end{align}
with a constant depending on $\bsb$. Combining \eqref{eq:taylor5} and
\eqref{eq:taylor6} yields a bound for $|T_2|$ in \eqref{eq:taylor2},
giving the second term in \eqref{trunc_err1}. This completes the proof for
$r \in \{0,1\}$.

To extend the result to $r \in (0,1)$, we apply operator interpolation to
the linear functional $T : H^r(D_R)' \rightarrow \mathbb{C}$ defined by
\[
 TG = \int_U G(u(\cdot, \bsy) - u(\cdot, \bsy_s))\,\rd \bsy ,
  \quad \text{for any } \quad G \in H^r(D_R)'.
\]
The above estimates show that there is a scalar $C_*$ independent of $r$,
such that $\Vert T \Vert_{H^r(D_R)' \rightarrow \mathbb{C}}  \leq C_*
k^{r-1}$ for $r = 0,1$. Thus the result for $r\in (0,1)$ follows on
recalling that the dual spaces $H^r(D_R)'$ are interpolation spaces.
\end{proof}

\section{QMC error bound}
\label{sec:QMC}

In this section we use QMC methods to approximate the expected value of a
linear functional of our PDE solution, $G(u_s(\cdot,\bsy))$, with respect
to $\bsy\in U_s := [-\frac{1}{2},\frac{1}{2}]^s$, i.e.,
\[
  \bbE[G u_s] \,=\, \int_{U_s} G(u_s(\cdot,\bsy))\,\rd\bsy.
\]
This general problem has been analyzed in many previous papers, for
example, in \cite{KSS12} with a family of first order QMC methods known as
``randomly shifted lattice rules'', and in \cite{DKLNS14} with a family of
higher order QMC methods known as ``interlaced polynomial lattice rules''.
The surveys \cite{KN16,KN18a} cover both methods. The paper \cite{GaKuSl}
also includes both methods for a Helmholtz PDE.

An $N$-point randomly shifted lattice rule \cite{DiKuSl:13} for the
integral of a function $\Theta$ over $U_s$ is given by
\begin{align} \label{eq:lattice}
  Q_{s,N,\bsDelta}(\Theta)
  \,=\, \frac{1}{N} \sum_{i=1}^N \Theta(\{\bst_i+\bsDelta\} - \bshalf),
\end{align}
where $\bst_i = \{i\bsz/N\} \in [0,1]^s$, $\bsz\in\bbN^s$ is the lattice
generating vector, the braces indicate that we take the fractional part of
each component in a vector, and $\bsDelta\in [0,1]^s$ is a random vector
where all components are i.i.d.\ uniform random numbers in $[0,1]$. The
subtraction by $\bshalf$ translates the standard unit cube $[0,1]^s$ to
$U_s$. An $N$-point interlaced polynomial lattice rule takes the form
\[
  Q_{s,N}(\Theta)
  \,=\, \frac{1}{N} \sum_{i=1}^N \Theta(\bst_i - \bshalf),
\]
where the points $\bst_i\in [0,1]^s$ are obtained by ``interlacing'' the
points of a ``polynomial lattice rule'', which are specified by a
generating vector of ``polynomials'' rather than of integers. For the
precise details as well as implementation, see, e.g., \cite{DKLNS14,KN16}
and the references there.

The general QMC theory (e.g., \cite{DiKuSl:13}) says that the integration
error depends on the function space setting assumed for the integrand, as
well as properties of the specific family of QMC methods. In a nutshell,
the absolute integration error is bounded by a product of two factors: (i)
the ``worse case error'' measuring the quality of the QMC method in the
function space, and (ii) the norm of the integrand in this function space.
Modern QMC analysis works with ``weighted'' function spaces, allowing
``weights'' to be chosen to model the dimension structure of a given
integrand. These chosen weights then enter the expression for the norm
and the worst case error, with the latter in turn used to construct
tailored QMC methods, leading to QMC error bounds with the best
convergence rate and with implied constant independent of dimension $s$.
To apply QMC theory to PDE problems, we need to bound the norm of the
integrand using the regularity estimates for the PDE solutions. Bounds on
the worst case errors are known from QMC theory.

In the following we make use of known QMC results without too much
further detail. However, we focus on one important new point to
address, namely, the dependence of the error bound on~$k$. For simplicity we state the results for $N$ being a power of a prime.

\begin{theorem} \label{thm:qmc}
Let $D$, $A$, $n$, and $R_0$ satisfy Assumptions~\ref{ass:1},
\ref{ass:3}~and~\ref{ass:4}, and let $k_0>0$.
Let $f\in L^2(D_{R_0})$ and let $G \in H^r(D_R)'$ for some $r \in [0,1]$.
For $s\ge 1$ and $\bsy\in U_s$, let $u_s(\cdot,\bsy)$ denote the solution
of the EDP with random coefficients where the series in
\eqref{eq:Axy} and \eqref{eq:nxy} are truncated to $s$ terms. Then for all
$k\ge k_0$ and $R \geq R_0$ we have the following QMC results.
\begin{itemize}
\item [\textnormal{(i)}] Randomly shifted lattice rules can be
    constructed such that
\begin{align*}
  \sqrt{\bbE_\bsDelta \big[ |\bbE[G u_s] - Q_{s,N,\bsDelta}(G u_s) |^2 \big]}
  \le C \, R\,k^{r-1}\,K_r(f,G)\,
  N^{-\min(1-\delta,\; 1/p-1/2)},
\end{align*}
where the constant $C$ grows unboundedly as $\delta\to 0$. If
$k$ satisfies $\Cstab\,k\,R\,\xi \le \Upsilon <\infty$ then $C$ is
bounded independently of $s$ but depends on $\Upsilon$ $($see
\eqref{eq:ugly} below$)$. If $s$ is bounded then the rate of
convergence is $N^{-1+\delta}$ and the constant $C$ is
proportional to $(\Cstab\,k\,R\,\xi)^s$.

\item [\textnormal{(ii)}] Interlaced polynomial lattice rules with
    interlacing factor $\alpha = \lfloor 1/p\rfloor + 1$ can be
    constructed such that
\begin{align*}
  |\bbE[G u_s] - Q_{s,N}(G u_s) |
  \le C\, R\,k^{r-1}\,K_r(f,G)\, N^{-1/p}.
\end{align*}
If $k$ satisfies $\Cstab\,k\,R\,\xi \le \Upsilon <\infty$ then $C$ is
bounded independently of $s$ but depends on $\alpha$
and~$\Upsilon$. If $s$ is bounded then the rate of convergence is
$N^{-\alpha}$ and the constant $C$ is proportional to $(\Cstab kR
\xi)^{\alpha\,s}$.
\end{itemize}
\end{theorem}

\begin{proof}
(i) From the standard QMC theory for randomly shifted lattice rules (see
e.g., \cite[formulas (11) and (12) with translated unit cube]{KN18a} or
\cite[Theorem 2.1]{KSS12}), we have
\begin{align*}
  E^{\rm QMC}
  := \sqrt{\bbE_\bsDelta \big[ |\bbE[G u_s] - Q_{s,N,\bsDelta}(G u_s) |^2 \big]}
  \le \bigg(
  \frac{2}{N} \sum_{\emptyset\ne\setu\{1:s\}} \gamma_\setu^\lambda\,[\vartheta(\lambda)]^{|\setu|}\bigg)^{\frac{1}{2\lambda}}\,
  \|Gu_s\|_{\calW_{s,\bsgamma}},
  \quad\forall\;\lambda\in (\tfrac{1}{2},1],
\end{align*}
with $\vartheta(\lambda) := \frac{2\zeta(2\lambda)}{(2\pi^2)^\lambda}$, 
where $\zeta$ denotes the Riemann-zeta function, and with the weighted Sobolev space norm
\[
 \|Gu_s\|_{\calW_{s,\bsgamma}}
 \,:=\, \bigg(
 \sum_{\setu\subseteq\{1:s\}} \frac{1}{\gamma_\setu}
 \int_{[-\frac{1}{2},\frac{1}{2}]^{|\setu|}}
 \bigg|\int_{[-\frac{1}{2},\frac{1}{2}]^{s-|\setu|}}
 \frac{\partial^{|\setu|}}{\partial\bsy_\setu} G(u_s(\cdot,\bsy))\,\rd\bsy_{\{1:s\}\setminus\setu}\bigg|^2
 \,\rd\bsy_\setu
 \bigg)^{\frac{1}{2}}.
\]
In this norm, for each subset $\setu\subseteq\{1:s\} :=
\{1,2,\ldots,s\}$, we take the mixed first derivative of the function with
respect to the variables $\bsy_\setu := (y_j)_{j\in\setu}$, integrate out
the remaining variables $\bsy_{\{1:s\}\setminus\setu} :=
(y_j)_{j\in\{1:s\}\setminus\setu}$, and then take the squared $L^2$ norm.
The weight parameter $\gamma_\setu>0$ moderates the relative importance
between the different subsets of variables. A small $\gamma_\setu$ means
that the function depends weakly on $\bsy_\setu$.

With $G\in H^r(D_R)'$, we obtain from \eqref{reg_est} for mixed first
derivatives (taking $\nu_j=1$ for $j\in\setu$ and $\nu_j=0$ for
$j\not\in\setu$) that
\begin{align} \label{eq:Gu-norm}
  \|G u_s\|_{\calW_{s,\bsgamma}}
  \,\le\,
  \Cstab\,R\,k^{r-1}\,K_r(f,G)\,
  \bigg(\sum_{\setu\subseteq\{1:s\}} \frac{(|\setu|!)^2 \prod_{j\in\setu} \beta_j^2}{\gamma_\setu}
  \bigg)^{\frac{1}{2}},
  \quad \beta_j := \Cstab \,k\, R\,\xi\,b_j.
\end{align}
Thus
\begin{align} \label{eq:part4}
  E^{\rm QMC}
  \,\le \Cstab\,R\,k^{r-1}\,K_r(f,G)\, \bigg(
  \frac{2}{N} \sum_{\emptyset\ne\setu\{1:s\}} \gamma_\setu^\lambda\,[\vartheta(\lambda)]^{|\setu|}\bigg)^{\frac{1}{2\lambda}}\,
  \bigg(\sum_{\setu\subseteq\{1:s\}} \frac{(|\setu|!)^2 \prod_{j\in\setu} \beta_j^2}{\gamma_\setu}
  \bigg)^{\frac{1}{2}}.
\end{align}
This is exactly of the form in \cite[Page 64 Step 12]{KN18a}. By choosing
the weights to minimize this bound (or equivalently and more simply, just
equating the terms inside the two sums), we obtain the weights
\begin{align} \label{eq:pod}
  \gamma_\setu \,=\, \bigg( |\setu|!\prod_{j\in\setu}
  \frac{\beta_j}{\sqrt{\vartheta(\lambda)}} \bigg)^{\frac{2}{1+\lambda}},
\end{align}
where $\beta_j$ is as defined in \eqref{eq:Gu-norm} and $\vartheta(\lambda) := \frac{2\zeta(2\lambda)}{(2\pi^2)^\lambda}$ is as stated earlier in the proof.
Substituting \eqref{eq:pod} into \eqref{eq:part4} yields finally
\begin{align} \label{eq:error}
  E^{\rm QMC}
  &\,\le\, \Cstab\,R\,k^{r-1}\,K_r(f,G)\,  \bigg(
  \frac{2}{N}\bigg)^{\frac{1}{2\lambda}}\,
  \Bigg(\sum_{\setu\subseteq\{1:s\}}
  \bigg(|\setu|!\,\prod_{j\in\setu} \big(\beta_j\,[\vartheta(\lambda)]^{\frac{1}{2\lambda}}\big)
  \bigg)^{\frac{2\lambda}{1+\lambda}} \Bigg)^{\frac{1+\lambda}{2\lambda}},
\end{align}
which is of the form in \cite[Page 64 Step 14]{KN18a}. Note that
the above result holds for all $s$ and for all $k$.

If $k$ is bounded and satisfies
\begin{align} \label{eq:restrict-k}
  \Cstab\,k\,R\,\xi \,\le\, \Upsilon < \infty
  \qquad\mbox{for some}\quad \Upsilon > 0,
\end{align}
then we may follow the argument in \cite[Page 64 Step 15 and Page 65
Step 16]{KN18a} to choose
\begin{align} \label{eq:lambda}
  \lambda = \begin{cases}
  \frac{1}{2-2\delta} \quad\mbox{for some } \delta\in (0,\frac{1}{2}) \quad & \mbox{if } p\in (0, \tfrac{2}{3}], \\
  \frac{p}{2-p}  & \mbox{if } p\in (\tfrac{2}{3},1).
  \end{cases}
\end{align}
This leads to the convergence rate $E^{\rm QMC} =
  \calO( N^{- \min(1-\delta,\; 1/p-1/2)})$,
where the implied constant is independent of $s$, but depends on
$\Upsilon$ through the factor
\begin{align} \label{eq:ugly}
  \bigg(\sum_{\ell=0}^\infty (\ell!)^{\eta-1}
  \bigg(\Upsilon^\eta\, [\vartheta(\lambda)]^{\eta/(2\lambda)} \sum_{j=1}^\infty
  b_j^\eta\bigg)^\ell \bigg)^{1/\eta},
\end{align}
with $\eta := 2\lambda/(1+\lambda) \in [p,1)$ and $\sum_{j=1}^\infty
b_j^\eta <\infty$. The outer sum in \eqref{eq:ugly} is finite by the
ratio test.

On the other hand, if $s$ is bounded then we may pull out a factor of
$(\Cstab\,k\,R\,\xi)^s$ from the last sum in \eqref{eq:part4}, leading to
the choice of weights in \eqref{eq:pod} with $\beta_j$ replaced by $b_j$. 
This is what we will do for our numerical experiment.

Weights of the form \eqref{eq:pod} are known as POD weights (``product and
order dependent weights''). The component-by-component construction
of lattice generating vector can be done for POD weights in $\calO(s\,
N\log N + s^2N)$ operations, see \cite{KSS12}.

(ii) The error for an interlaced polynomial lattice rule satisfies
(see e.g., \cite[Theorem~5.4]{KN16} or \cite[Theorem~3.1]{DKLNS14})
\begin{align*} 
  E^{\rm QMC} := \big|\bbE [G u_s] - Q_{s,N}(G u_s)\big|
  \,\le\, \bigg(\frac{2}{N} \sum_{\emptyset\ne\setu\subseteq\{1:s\}} \gamma_\setu^\lambda\,
  [\varrho_\alpha(\lambda)]^{|\setu|} \bigg)^{1/(2\lambda)}\,\|G u_s\|_{s,\alpha,\bsgamma}
  \quad\forall\;\lambda\in (\tfrac{1}{\alpha},1],
\end{align*}
where $\alpha\ge 2$ is an integer smoothness parameter (also known as the
``interlacing factor''), $N$ is a power of~$2$, $\varrho_\alpha(\lambda)
  = 2^{\alpha\lambda(\alpha-1)/2} [(1 +
  \frac{1}{2^{\alpha\lambda}-2})^\alpha-1]$, and the norm is
\begin{align*} 
  \|G u_s\|_{s,\alpha,\bsgamma}
 :=
 \sup_{\setu\subseteq\{1:s\}}
 \sup_{\bsy_\setu \in [0,1]^{|\setu|}}
 \frac{1}{\gamma_\setu}
 \sum_{\setv\subseteq\setu} \,
 \sum_{\bstau_{\setu\setminus\setv} \in \{1:\alpha\}^{|\setu\setminus\setv|}}
 \bigg|\int_{[-\frac{1}{2},\frac{1}{2}]^{s-|\setv|}}
 (\partial^{(\bsalpha_\setv,\bstau_{\setu\setminus\setv},\bszero)} G(u_s(\cdot,\bsy)) \,
 \rd \bsy_{\{1:s\} \setminus\setv}
 \bigg|\,.
\end{align*}
Here
$\partial^{(\bsalpha_\setv,\bstau_{\setu\setminus\setv},\bszero)}$
indicates that we take the partial derivative $\alpha$ times with respect
to $y_j$ for $j\in\setv$ and $\tau_j$ times with respect to $y_j$ for
$j\in\setu\setminus\setv$. Using again \eqref{reg_est} (this time with
general $\bsnu$), we obtain instead of \eqref{eq:Gu-norm},
\begin{align} \label{eq:error2}
 \|Gu_s\|_{s,\alpha,\bsgamma}
 & \,\le\,
 \Cstab\,R\,k^{r-1}\,K_r(f,G)\,
 \sup_{\setu\subseteq\{1:s\}}
 \frac{1}{\gamma_\setu}
 \sum_{\bsnu_\setu \in \{1:\alpha\}^{|\setu|}}
 |\bsnu_\setu|!\,\prod_{j\in\setu} \big(2^{\delta(\nu_j,\alpha)} \beta_j^{\nu_j}\big)\,,
\end{align}
where $\delta(\nu_j,\alpha)$ is $1$ if $\nu_j=\alpha$ and is $0$
otherwise.

For bounded $k$ satisfying \eqref{eq:restrict-k}, we now choose
$\gamma_\setu$ so that the supremum in \eqref{eq:error2} is $1$,
i.e.,
\begin{equation} \label{eq:spod}
 \gamma_\setu
 \,=\,
 \sum_{\bsnu_\setu \in \{1:\alpha\}^{|\setu|}}
 |\bsnu_\setu|!\,\prod_{j\in\setu} \big(2^{\delta(\nu_j,\alpha)} \beta_j^{\nu_j}\big)\,.
\end{equation}
Using the above weights and following the arguments in
\cite[Pages~2694--2695]{DKLNS14}, by taking $\lambda = p$ and the
interlacing factor $\alpha = \lfloor 1/p\rfloor + 1$, we eventually arrive
at the convergence rate $E^{\rm QMC} = \calO(N^{-1/p})$, with the
implied constant independent of $s$, but depends on $\alpha$ and
$\Upsilon$ through a factor similar to \eqref{eq:ugly}.

For bounded $s$, we may pull out a factor of
$(\Cstab\,k\,R\,\xi)^{\alpha\,s}$ from the supremum in \eqref{eq:error2},
leading to the choice of weights in \eqref{eq:spod} with $\beta_j$
replaced by $b_j$.

Weights of the form \eqref{eq:spod} are called SPOD weights
(``smoothness-driven product and order dependent weights''). The
generating vector (of polynomials) can be obtained by a
component-by-component construction in $\calO(\alpha\,s\, N\log N +
\alpha^2\,s^2 N)$ operations, see \cite{DKLNS14}.
\end{proof}

\section{PML truncation and FEM approximation}
\label{sec:FEM}

\subsection{Definition of radial PML truncation}

We now approximate the EDP \eqref{eq:vp}--\eqref{eq:agen},  
using a perfectly-matched layer (PML) truncation. 
The PML is chosen to minimise spurious
reflections, and will act in the annulus $B_{R_2} \backslash
\overline{B_{R_1}}$ where 
\begin{align} \label{eq:annulus} 
 R_0< R_1 < R_{2}<\infty, 
\end{align}
with  $R_0$  as given in Assumption~\ref{ass:1}. The truncation can be
thought of as a practical approximation of the Dirichlet to Neumann map
$\DtN$ appearing in \eqref{eq:agen} (with $R = R_1$). But,
while $\DtN$ is a global operator on the surface $\partial B_R$, the PML
only requires computations with  local operators. While PMLs have been
well-known since the work of Berenger \cite{Be:94} (see also
Remark~\ref{rem:relation}), we make particular use here of very recent
theoretical results from \cite{GLS2} proving that fixed width PML
truncation is exponentially accurate as $k \rightarrow \infty$ (see
Theorem \ref{thm:GLS} below).
 
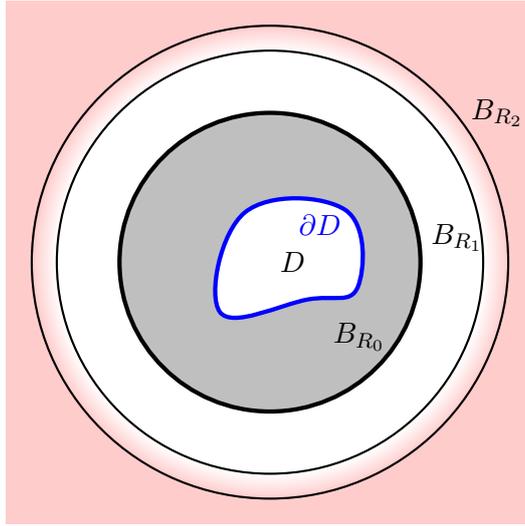
\begin{figure}[t]
\begin{center}
\begin{tikzpicture}[scale=0.33]
 \fill[red!20!white] (-10.5,-10.5) rectangle (10.5,10.5);
 \shade[even odd rule,ring shading={from white at 8.5 to red!20!white at 9.5}] (0,0) circle (8.5) circle (9.5);
 \draw[thick] (0,0) circle (9.5); 
 \draw[thick,fill=white] (0,0) circle (8.5); 
 \draw[ultra thick,fill=lightgray] (0,0) circle (6);   
 \draw[ultra thick, pattern={Lines[angle=45,distance=8pt]}, blue, fill=white]
   plot[smooth cycle, tension=.7]
   coordinates {(-2,-2) (-1,2) (3,2.1) (3.5,-1) (1.5,-1.5)};
 \node[right] at (0,0) {$D$};
 \node[blue] at (2,1.5) {$\partial D$};
 \node[left] at (5,-3) {$B_{R_0}$};
 \node[left] at (8.9,1) {$B_{R_1}$};
 \node[left] at (10.5,6) {$B_{R_2}$};
\end{tikzpicture}
\end{center}
\caption{The PML truncation problem is posed on the domain $D_{R_2}$
which is the region inside the ball $B_{R_2}$, excluding $\overline{D}$.
The red gradient shading illustrates the radial cutoff function $\varphi_{\rm PML}\in C^3(\bbR_+)$
satisfying \eqref{eq:sigma_prop}; it is increasing in the annulus
between radii $R_1$ and $R_2$.} \label{fig:PML}
\end{figure}

The PML truncated problem (which is subsequently solved computationally)
is then posed on  
\begin{align} \label{eq:defOmega}
  D_{R_2} = B_{R_2} \backslash  \overline{D}.
\end{align}
Since $R_1 > R_0$, we have $A = I$, $n = 1$ and $f = 0$ in a neighbourhood
of the region of truncation.
We illustrate this in Figure~\ref{fig:PML}.

\medskip

For what follows we adopt the following notational convention.
\begin{notation}\label{not:suppress}
Throughout the rest of the paper $A$ and $n$ (also $A_{\rm PML}$ and
$n_{\rm PML}$ -- introduced below) are random coefficients which depend on
$\bsx$, $\bsy$ (recall \eqref{eq:Axy}, \eqref{eq:nxy}). We suppress
dependence on $\bsx$, $\bsy$ from the notation, except where is needed to
aid understanding. Similarly the sesquilinear forms $a_k$ (from
\eqref{eq:agen}) and $a_{\rm PML}$ (introduced below) also depend on
$\bsy$ but we also suppress this from the notation.
\end{notation}

\begin{remark}[Choice of computational domain]\label{rem:poly}
To make the PML truncation as simple as possible we have assumed here that
the truncation boundary is $\partial B_{R_2}$. However this spherical
truncation boundary could be replaced by a polygon (or polyhedron), as
long as it contains $B_{R_1}$, and does not touch $\partial B_{R_1}$. We
explain in Remark~\ref{rem:curved} how the theory extends to this case.
\end{remark}

Under Assumption \ref{ass:1}, let $u$ be the solution of the EDP
\eqref{eq:vp}--\eqref{eq:agen}. Then the PML truncated approximation
$u_{\rm PML}\in H^1_0(D_{R_2})$ is defined to be the solution of the
variational problem
\begin{align}\label{eq:PMLvf}
 a_{\rm PML}(u_{\rm PML},v):=\int_{D_{R_2}}
 \Big(
 (A_{\rm PML} \nabla u_{\rm PML}) \cdot \overline{\nabla v} -k^2\, n_{\rm PML}\,
 u_{\rm PML}\,\overline{v}\Big)
 = \int_{D_{R_2}} f\,\overline{v} \quad \text{for all} \  {v} \in H^1_0(D_{R_2}).
\end{align}
Here ${A}_{\rm PML}$ and ${n}_{\rm PML}$ are complex modifications of
${A}$ and $ {n}$, which are defined as follows.

First, choose a univariate smooth cutoff function $\varphi_{\rm PML} \in C^3(\mathbb{R}_+)$, satisfying
\begin{align} \label{eq:sigma_prop}
\begin{cases}
\varphi_{\rm PML}(r)= 0
\text{ for }
 r\leq R_1, \\
\varphi_{\rm PML}(r) \text{ is increasing for } R_1 \leq r\leq R_2, \text{ and }\\
\varphi_{\rm PML}(r) = \varphi_{\rm PML, const} > 0 \text{ for } r\geq R_2,
\end{cases}
\end{align}
and then define
 \beq\label{eq:sigma}
 \sigma(r) := \big( r\, \varphi_{\rm PML}(r)\big)', \quad
 \alpha(r) := 1 + \ri\, \sigma(r), \quad \text{ and }\quad
 \beta(r) := 1 + \ri\, \varphi_{\rm PML}(r).
 \eeq

Then $A_{\rm PML}$ and $n_{\rm PML}$ are defined by:
\beq\label{eq:defAPML} {A}_{\rm PML} :=
\begin{cases}
A& \text{ for } r \leq R_1 \\
HKH^\top &\text{ for } r> R_1,
\end{cases}
\quad\text{ and }\quad
n_{\rm PML}:=
\begin{cases}
n& \text{ for } r \leq R_1 \\
\alpha(r) \beta(r)^{d-1} &\text{ for } r> R_1,
\end{cases}
\eeq
where, in polar coordinates $(r,\theta)$,
 \beqs
K =
\left(
\begin{array}{cc}
\beta(r)\alpha(r)^{-1} &0 \\
0 & \alpha(r) \beta(r)^{-1}
\end{array}
\right)
\quad\text{ and }\quad
H =
\left(
\begin{array}{cc}
\cos \theta & - \sin\theta \\
\sin \theta & \cos\theta
\end{array}
\right) \text{ for } d=2,
\eeqs
and, in spherical polar coordinates
$(r,\theta,\phi)$, \beqs K = \left(
\begin{array}{ccc}
\beta(r)^2\alpha(r)^{-1} &0 &0\\
0 & \alpha(r) &0 \\
0 & 0 &\alpha(r)
\end{array}
\right)
\,\,\text{ and }\,\,
H =
\left(
\begin{array}{ccc}
\sin \theta \cos\phi & \cos \theta \cos\phi & - \sin \phi \\
\sin \theta \sin\phi & \cos \theta \sin\phi & \cos \phi \\
\cos \theta & - \sin \theta & 0
\end{array}
\right) \text{ for } d=3. \eeqs

We note that $A=I$ and $n=1$ when $r=R_1$ and thus ${A}_{\rm PML}$ and
${n}_{\rm PML}$ are continuous at $r=R_1$.

\begin{remark}[The definition of the PML scaling function $\sigma$] \label{rem:pmlscf}
Here we have started from $\varphi_{\rm PML}$, and defined $\sigma$
in terms of $\varphi_{\rm PML}$, thus following, e.g.,
\cite[\S2]{BrPa:07}, \cite[\S1.2]{GLS2}. Alternatively, one can start from
a non-decreasing function $\sigma$ and define $\varphi_{\rm PML}$
such that the first equation in \eqref{eq:sigma} holds; see, e.g.,
\cite[\S3]{CoMo:98a}, \cite[\S2]{LaSo:98}, \cite[\S4]{HoScZs:03},
\cite[\S2]{ChGaNiTo:22}. The notation $\alpha$ and $\beta$ is also used by
\cite{LaSo:98, LiWu:19}.
\end{remark}

\begin{remark}[Relation to notation used in {\cite{GLS2, GLSW1}}] \label{rem:relation}
In \cite{GLS2,GLSW1} the PML is defined in terms of a scaling function
$f_\theta$ (see \cite[\S 1.2]{GLS2}, \cite[\S 1.3]{GLSW1}); this
corresponds to our function $\varphi_{\rm PML}$ above via the
relation $\varphi_{\rm PML}(r)= f_\theta(r)/r$.
The use in \cite{GLS2, GLSW1} of the subscript $\theta$ in $f_\theta$
originates in  the notation used in the method of complex scaling -- the
original name for PML when it was developed in the 1970s by analysts; see,
e.g., the references in \cite[\S4.7]{DyZw:19} (and see also comments given
in the proof of Lemma~\ref{lem:signPML}).
\end{remark}

\subsection{Properties of the PML sesquilinear form $\aPML$ and
solution $\uPML$}

In this subsection we collect some known properties of problem
\eqref{eq:PMLvf} that are needed in the proof of Theorem \ref{thm:FEM1}.
(Recall Notation~\ref{not:suppress}.)

\begin{lemma} 
\label{lem:signPML} Suppose $A,n$ satisfy Assumption~\ref{ass:3} and
$\varphi_{\rm PML}$ is chosen as in \eqref{eq:sigma_prop}. Then there
exists $C>0$ such that for all $\bsy \in U$ and $\bsx \in \mathbb{R}^d$,
\beq\label{eq:signPML}
 \Re \big( ({A}_{\rm PML}\bszeta) \cdot \bszeta \big)
\ \geq \ C |\bszeta|^2_2,  \quad \quad \text{for all} \quad  \bszeta \in
\mathbb{C}^d . \eeq
\end{lemma}

\begin{proof}
By Assumption \ref{ass:3} and Lemma \ref{lem:easy}, $A$ is real symmetric
and positive definite so \eqref{eq:signPML} holds with $A_{\rm PML}$
replaced by~$A$. Also \cite[Lemma 2.3]{GLSW1} proves \eqref{eq:signPML}
for $HKH^{\top}$. Thus the result follows from the definition
\eqref{eq:defAPML} of $A_{\rm PML}$.
(Note that \cite[Lemma 2.3]{GLSW1} actually considers a family of scaling
functions $f_\theta(r) : =  f(r) \tan \theta$, for a given $f$,
parametrized by $\theta \in (0,\pi/2)$, and proves the estimate of the
form \eqref{eq:signPML} uniformly with respect to $ \epsilon < \theta <
\pi/2 - \epsilon$ for any $\epsilon > 0$. However, we do not need that
detail here.)
\end{proof}

\begin{corollary}[Continuity and G\aa rding inequality for $\aPML$]\label{cor:PMLGarding}
Suppose $A$, $n$ satisfy Assumption \ref{ass:3} and $\varphi_{\rm PML}$ is as in \eqref{eq:sigma_prop}. Then there exists constants
$C_1,C_2,C_3>0$ such that, for all $\bsy \in U$, $k>0$, and $u,v \in
H^1(D_{R_2})$,
\begin{align} \label{eq:PML_cont}
 |a_{\rm PML}(u,v)|\ \leq\  C_1\,
 \N{u}_{H^1_k(D_{R_2})}\,\N{v}_{H^1_k(D_{R_2})} \quad \text{ (Continuity)} ,
\end{align}
and
\begin{align}  \label{eq:PML_Gaa}
 \Re\, a_{\rm PML}(v,v) \ \geq\  C_2\,
 \N{v}^2_{H^1_k(D_{R_2})} - C_3\, k^2\, \N{v}^2_{L^2(D_{R_2})} \quad \text{
 (G\aa rding inequality). }
\end{align}
\end{corollary}

\begin{proof}
The continuity property follows from the definitions of $a_{\rm
PML}(\cdot,\cdot)$ and $\|\cdot\|_{H^1_k(D_{R_2})}$, together with the
Cauchy--Schwarz inequality. For the G\aa rding inequality, we write
\[
 \Re\, a_{\rm PML}
 (v,v) = \int_{D_{R_2}} \Re
  \left(  (A_{\rm PML} \nabla v) \cdot   \nabla \overline{v} \right)
  - k^2\, n_{\rm PML}\, \vert v \vert^2 .
\]
Then we apply the reverse triangle inequality, using
Lemma~\ref{lem:signPML} and also the fact that $n_{\rm PML} $ can be
bounded above on $D_{R_2}$.
\end{proof}

The following theorem shows that the PML solution operator inherits the
bound satisfied by the non-truncated solution operator, up to a
$k$-independent constant, denoted $C$ below, uniformly in $\bsy\in U$.
(Compare \eqref{eq:inherit1} to \eqref{eq:EDP1}.)

\begin{theorem}[Bound on the solution of PML problem]\label{thm:inherit}
Let $D$, $A$, $n$, and $R_0$ satisfy Assumptions~\ref{ass:1},
\ref{ass:3}~and~\ref{ass:4}, and let $R_1> R_0$. Then, given $\epsilon>0$,
there exist $C>0$ and $k_0>0$ such that for all $R_{2} > (1+\epsilon)R_1$,
$k\geq k_0$, $\bsy \in U$, and $f\in L^2(D_{R_0})$, with
$\varphi_{\rm PML}$ as in \eqref{eq:sigma_prop}, the solution
$u_{\rm PML}(\cdot,\bsy)$ of \eqref{eq:PMLvf} exists, is unique, and
satisfies the bound
\begin{equation} \label{eq:inherit1}
 \N{u_{\rm PML}(\cdot,\bsy)}_{H^1_k (D_{R_2})} \leq C\, \Cstab\,
 R_2\, \|{f}\|_{L^2(D_{R_0})}\,,
\end{equation}
where $\Cstab$ is defined in \eqref{eq:Cstab}.
\end{theorem}

We highlight that, both here and elsewhere, the order of the quantifiers
in the statement of a result indicates on what quantities constants
depend. For example, the constant $C$ in Theorem~\ref{thm:inherit} depends
on $\epsilon$ but is independent of $k$ and $\bsy$.

\begin{proof}[Proof of Theorem \ref{thm:inherit}] This result follows from
\cite[Theorem 1.6]{GLS2}. To see this, note that (i) here we are taking
$R_{\rm tr}$ in \cite{GLS2} to be $R_2$; (ii) the norm $\|\chi
R_P(k)\chi\|_{\cH\to\cH}$ on the right-hand side of \cite[Equation
(1.22)]{GLS2} is the norm of the solution operator for the untruncated
problem, which in our case is estimated by \eqref{eq:EDP1}; (iii)  the
norm $\|\cdot\|_{\cH(\Omega_{\rm tr})}$ in \cite{GLS2} corresponds to $\|
\cdot \|_{ L^2(D_{R_2})}$ and the norm $\|\cdot\|_{\cD(\Omega_{\rm tr})}$
in \cite[Theorem 1.6]{GLS2} is the norm on the domain of the PDE operator
(see \cite[Equation 2.3]{GLS2}) and thus controls the norm
$\|\cdot\|_{H^1_k(D_{R_2})}$.
\end{proof}

\subsection{The accuracy of PML truncation for $k$ large} \label{sec:PML_error}

The following theorem shows that the error in approximating
$u(\cdot,\bsy)$ by $u_{\rm PML}(\cdot,\bsy)$, \emph{when measured on the
domain $D_{R_1}$} (i.e., the part of  $D_{R_2}$ that does not contain the
PML), decreases exponentially with respect to the wavenumber $k$ and
also with respect to the PML depth $R_2 - (1 + \epsilon)R_1$, for any
$\epsilon>0$. The exponential decay is uniform over all  $\bsy\in U$.

\begin{theorem}[Radial PMLs are exponentially accurate for $k$ large]\label{thm:GLS}
Let $D$, $A$, $n$, and $R_0$ satisfy Assumptions~\ref{ass:1},
\ref{ass:3}~and~\ref{ass:4}, and let $R_1> R_0$ and $k_0>0$. Then,
given $\epsilon
>0$, there exist $C_{\rm PML, 1}>0$, $C_{\rm PML,2}>0$ and $k_0>0$ such
that for all $R_{2}>(1 + \epsilon)R_1$, $k\geq k_0$, $\bsy \in U$, and
$f\in L^2(D_{R_0})$, with $\varphi_{\rm PML}$ as in
\eqref{eq:sigma_prop}, the solution $u_{\rm PML}(\cdot,\bsy)$ of
\eqref{eq:PMLvf} satisfies the bound
\beq\label{eq:GLS}
 \|u(\cdot,\bsy)-u_{\rm PML}(\cdot,\bsy)\|_{H^1_k(D_{R_1})}
 \ \leq \ C_{\rm PML, 1} \, \exp \big( - C_{\rm PML, 2}\,  k\big(R_2 -  (1 + \epsilon) R_1\big)
 \big)\,
 \|f\|_{L^2(D_{R_0})} .
\eeq
\end{theorem}

\bpf This result is a special case of \cite[Theorem 1.5]{GLS2}. To see
this note that (i) here we considering a nontrapping problem, so that in
the notation of \cite{GLS2}, $\Lambda(P,J)=0$ (see the discussion below
equation (1.18) in \cite{GLS2}); (ii) the role of $R_{\rm tr}$ in
\cite{GLS2} is played by $R_2$ here; (iii) the estimates in
\cite[Equations (1.11), (1.16), and (1.21)]{GLS2} imply that the
exponential on the right-hand side of \cite[Equation 1.21]{GLS2} is of the
form of the one on the right-hand side \eqref{eq:GLS}; (iv) as noted in
the proof of Theorem \ref{thm:inherit}, the norm
$\|\cdot\|_{\cD(\Omega_{\rm tr})}$ in \cite[Theorem 1.5]{GLS2} controls
the norm $\|\cdot\|_{H^1_k(D_{R_1})}$.

In fact \cite[Theorem 1.5]{GLS2} actually proves that the error in the PML
truncation decays exponentially not only in the parameters $k$ and $R_2 -
(1 + \epsilon) R_1$, but also in the the PML scaling parameter -- here
denoted $\varphi_{\rm PML,const}$ in \eqref{eq:sigma_prop}. However,
to reduce technical detail, here we have hidden the dependence on
$\varphi_{\rm PML,const}$ in the constant $C_{\rm PML, 2}$.
Increasing the scaling parameter provides more damping. \epf

\subsection{The accuracy of the $h$-FEM approximation of the PML solution}

In this section we give the theory of the approximation of problem
\eqref{eq:PMLvf} using a family of finite dimensional spaces
$(V_h)_{h>0}$.  The required approximation properties of $V_h$ as $h
\rightarrow 0$ are encapsulated in the following abstract assumption.

\begin{assumption}[Approximation property of $V_h$] \label{ass:spaces}
There exists approximation degree $m\in \mathbb{Z}^+$, and constants $C =
C(m)>0$ and $h_0>0$ such that, for all $R_2>R_0$ and for any $v \in H^1_0(D_{R_2})\cap H^{r'+1}(D_{R_2})$
with $0 \leq r' \leq m$, there exists an approximant $w_h\in V_{h}$
satisfying
 \beq\label{eq:pp_approx}
 \|v- w_h \|_{H^r(D_{R_2})} \leq C\, h^{r' + 1-r} \|v\|_{H^{r' +1}(D_{R_2})}
 \, \text{ for all }\, 0<h<h_0\, \text{ and } \, 0 \leq r \leq r' \leq
 m.
 \eeq
\end{assumption}

\begin{example}
When $D_{R_2}$ is a polygon/polyhedron, Assumption~\ref{ass:spaces} holds
when $(V_h)_{h>0}$ consists of continuous piecewise degree-$m$ Lagrange elements on shape-regular simplicial triangulations, indexed by the
meshwidth; see, e.g., \cite[Theorem 17.1]{Ci:91}, \cite[Proposition
3.3.17]{BrSc:08}. 
For a discussion of the geometric error incurred by using simplicial
elements on curved domains, see Remark \ref{rem:curved} below.
When $D_{R_2}$ is curved, Assumption~\ref{ass:spaces} holds for piecewise polynomials on curved meshes satisfying certain (natural) derivative bounds on the element maps; see 
\cite[Theorem~2]{CiRa:72}, \cite[Theorem~1]{Le:86}, \cite[Theorem~5.1]{Be:89}.
\end{example}

Now, using $V_h$, the solution of PML problem \eqref{eq:PMLvf} is defined
by the usual Galerkin formulation:
\begin{equation}\label{eq:FEM}
  \text{Find} \,\, u_{\rm PML}(\cdot,\bsy)_h \in V_h\,\, \text{ such that }\,\,
  a_{\rm PML}(u_{\rm PML}(\cdot,\bsy)_h ,{v}_h)
  = \int_{D_{R_2}}f\,\overline{v_h} \quad\text{ for all } v_h \in V_h.
\end{equation}

The convergence estimates in Theorem~\ref{thm:FEM1} below require the
following regularity assumptions on the expansion functions in
\eqref{eq:Axy}, \eqref{eq:nxy}.

\begin{assumption}[Additional smoothness of random coefficients]
\label{ass:smoothness} Following Assumption~\ref{ass:3}, suppose $A_0$,
$n_0$, $\Psi_j$, and $\psi_j$ are all in $W^{\tau,\infty}(D_{R_0})$ for
some smoothness degree $\tau \in \mathbb{Z}^+$, and that
 \beq\label{eq:Winftybound}
 \sum_{j=1}^\infty \N{\Psi_j}_{W^{\tau,\infty}(D_{R_0})} +\sum_{j=1}^\infty
 \N{\psi_j}_{W^{\tau,\infty}(D_{R_0})} < \infty,
 \eeq
where the norm for matrix-valued functions are defined analogously to
Notation~\ref{not:oneinf}.
\end{assumption}

Assumption~\ref{ass:smoothness} coincides with Assumption~\ref{ass:3} when
$\tau=1$, but when $\tau>1$ Assumption~\ref{ass:smoothness} is stronger.

The following theorem provides $k$-explicit error estimates for the FEM
approximation \eqref{eq:FEM} of the PML problem \eqref{eq:PMLvf}. These
error estimates hold on the whole computational domain $D_{R_2}$, which
includes the PML (recall \eqref{eq:defOmega}). For simplicity, we restrict
attention to the case where the data $f$ is a $k$-dependent oscillatory
function -- such an $f$ arises when we consider the plane-wave scattering
problem in the following section.

\begin{theorem}[$k$-explicit quasioptimality of FEM, uniformly in $\bsy \in U$]\label{thm:FEM1}
Let $D$, $A$, $n$, and $R_0$ satisfy Assumptions~\ref{ass:1},
\ref{ass:3}~and~\ref{ass:4}, and let $R_1> R_0$. Suppose $A$ and $n$ also
satisfy Assumption~\ref{ass:smoothness} with smoothness degree $\tau$. Let
$(V_h)_{h>0}$ satisfy Assumption~\ref{ass:spaces} with approximation
degree $m$. For
\[
 \ell :=\min(\tau,m),
\]
suppose that $\partial D_{R_2}$ is H\"older continuous $C^{\ell,1}$ (see, e.g., \cite[Page 90]{Mc:00}) and
that $\varphi_{\rm PML}$ as defined in \eqref{eq:sigma_prop} is in
$C^{\ell,1}(\bbR^+)$. 
Then, for all $\epsilon>0$, there exist $C_{\rm FEM,1}>0$ and $k_0>0$, and for all $\Coscil>0$ there exists $C_{\rm FEM,2}>0$ such that the following is true for all $R_{2}
> (1+\epsilon)R_1$, $k\geq k_0$, and $\bsy\in U$.

Suppose $f\in L^2(D_{R_0})$ is a
$k$-dependent oscillatory function in the sense that $f\in
H^{\ell-1}(D_{R_0})$ with 
 \beq     \label{eq:koscillatory}
 \|f\|_{H^{\ell-1}(D_{R_0})} \leq \Coscil\, k^{\ell-1}\, \|f\|_{L^2(D_{R_0})}
 \quad \tfa \  k\geq k_0.
 \eeq
Then, if $h = h(k)$ satisfies
\begin{gather}
  \label{eq:FEMthreshold}
(hk)^{2\ell}\, k\,R_2   \leq  C_{\rm FEM, 1},
\end{gather}
the solution $u_{\rm PML,h}(\cdot,\bsy)$ of \eqref{eq:FEM}
satisfies
\begin{align}
   \big\| u_{\rm PML}(\cdot,\bsy)
-u_{\rm PML, h}(\cdot,\bsy) \big\|_{H^1_k(D_{R_2})}
&\ \leq\  C_{\rm FEM, 2}\,
\big((hk)^\ell \;\;\; + (hk)^{2\ell}\, k\, R_2\big)\,
\|u_{\rm PML}(\cdot,\bsy)\|_{H^1_k(D_{R_2})}\,,
\label{eq:FEMrel1}
\\
k\, \big\| u_{\rm PML}(\cdot,\bsy)
-u_{\rm PML, h}(\cdot,\bsy) \big\|_{L^2(D_{R_2})}
&\ \leq\  C_{\rm FEM, 2}\,
\big( (hk)^{\ell+1} + (hk)^{2\ell}\, k\, R_2\big)\,
\|u_{\rm PML}(\cdot,\bsy)\|_{H^1_k(D_{R_2})}\,,
\label{eq:FEMrel2}
\end{align}
where $u_{\rm PML}(\cdot,\bsy)$ is the solution of \eqref{eq:PMLvf}.
\end{theorem}

Theorem \ref{thm:FEM1} implies that if $(hk)^{2\ell}\,k\,R_2$ is kept sufficiently small as $k \rightarrow \infty$,  then (by \eqref{eq:FEMrel1}) the relative $H^1_k(D_{R_2})$ error is
controllably small, uniformly in $k$, as $k\to \infty$. This is because
the condition \eqref{eq:FEMthreshold} implies $hk \le (C_{{\rm
FEM},1}/(k_0\, R_2))^{1/(2\ell)}$. The condition \eqref{eq:FEMthreshold} is
observed empirically to be sharp (going back to the work of Ihlenburg and
Babu\v{s}ka in 1D \cite{IhBa:95a, IhBa:97}); see the discussion in
\cite[\S1.2]{GS3} and the references therein.

\bpf 
This result for fixed $\bsy\in U$ (i.e., with the constants $C_{\rm FEM, 1}$ and $C_{\rm FEM,2}$ depending on $\bsy$)
 follows immediately from the Helmholtz $h$-FEM analysis in \cite{GS3}. This $h$-FEM analysis holds for the Galerkin method applied to sesquilinear forms that are continuous, satisfy a G\aa rding inequality, and satisfy the natural elliptic-regularity shift
results (i.e., if the data is in $H^{\ell-1}(D_{R_0})$ and the
coefficients and domain boundary are sufficiently smooth, then the
solution is in $H^{\ell+1}(D_{R_2})$; see, e.g.,  \cite[Theorem 4.18]{Mc:00}). Furthermore, \cite[Theorem 4.9]{GS3} shows that the Helmholtz PML problem we consider here falls into this category.
Finally, the assumption made on the FEM spaces in \cite{GS3}, namely \cite[Assumption 4.8]{GS3}, is equivalent to Assumption \ref{ass:spaces} (noting that, 
since the function $x \mapsto x^{r'-r+1}$ is increasing for $r'\geq r$, the bound in \eqref{eq:pp_approx} holding for $h$ sufficiently small
is equivalent to the bound holding for all $h>0$).

However, Theorem \ref{thm:FEM1} requires that the constants $C_{\rm FEM, 1}$ and $C_{\rm FEM,2}$ be independent of $\bsy$. 
For reasons of space, we do not reproduce here the arguments in \cite{GS3}, but instead explain why $C_{\rm FEM, 1}$ and $C_{\rm FEM,2}$ can be taken independent of $\bsy$.
The key point is that the constants $C_{\rm FEM, 1}$ and $C_{\rm FEM,2}$ depend only on 
 (i)
the constants in the continuity bound and G\aa rding inequality, which are independent of $\bsy$ by Corollary~\ref{cor:PMLGarding}, and 
 (ii) the constants in the elliptic-regularity shift results, which are independent of $\bsy$  by the $W^{\ell, \infty}$ bounds provided by
Assumption~\ref{ass:smoothness} (since $\ell\le\tau$). 
Thus the results of \cite[Theorem~4.9]{GS3} hold with $C_{\rm FEM, 1}$ and $C_{\rm FEM,2}$ independent of $\bsy$.

As a result, \eqref{eq:FEMrel1} follows from \cite[Eq.~4.19]{GS3}. The
bound \eqref{eq:FEMrel2} is not stated explicitly in \cite{GS3}, but
follows from the last displayed equation in \cite[\S2.2]{GS3} (one repeats
the arguments that obtain \eqref{eq:FEMrel1} from the quasi-optimality
bound \cite[Eq.~4.16]{GS3} but now one starts from the $L^2$ bound
\cite[Eq.~4.17]{GS3}).
\epf

\bre[Approximation of more general boundaries]\label{rem:curved} Theorem
\ref{thm:FEM1} requires that $\partial D_{R_2}$ is at least $C^{1,1}$, and
gives better convergence rates under additional smoothness of $\partial
D_{R_2}$. In these cases, using simplicial triangulations incurs a
non-conforming error coming from the approximation of $\partial D_{R_2}$.
The result of Theorem \ref{thm:FEM1} holds verbatim when this error is
taken into account, with the Galerkin errors on the left-hand sides of
\eqref{eq:FEMrel1} and \eqref{eq:FEMrel2} measured in the subset of
$D_{R_2}$ where both $u_{\rm PML}$ and $u_{{\rm PML},h}$ are defined. This
is thanks to the recent result of \cite{ChSp:24}, which extend the
arguments of \cite{GS3} (used in the proof of Theorem \ref{thm:FEM1}) to
the non-conforming setting. The key point is that, for nontrapping
problems solved with simplicial elements, the geometric error is of order
$hk$, which is smaller than the pollution error $(hk)^{2\ell} k R_2$
under the mesh condition \eqref{eq:FEMthreshold} (see \cite[Main Result 1.1 and 
Theorem 4.11]{ChSp:24}).
\ere

We now apply the results above to the estimation of the quantity $\EFEM$
defined in \eqref{eq:tentative1}. In this connection we make the following
remark.
\begin{remark}[Dimension truncation]\label{rem:alls}
Now recall from \eqref{eq:tentative1} that the analysis of $\EFEM$
requires estimating the error in approximating the dimension truncated
solution $u_s$ (as defined in \eqref{eq:defus}) by its combined PML-FEM
approximation $u_{s,{\rm PML},h}$. However, because the estimates in these
theorems are uniform in $\bsy \in U$, they hold without modification for
$\bsy$ replaced by $\bsy_s$.
\end{remark}

Using this we obtain the following Corollary to Theorems~\ref{thm:GLS}
and~\ref{thm:FEM1}.

\begin{corollary} \label{cor:functional}
Assume that the conditions of Theorems~\ref{thm:GLS} and~\ref{thm:FEM1}
all hold and assume $G \in H^r(D_R)'$ for $r \in [0,1]$ and $R_0\leq R
\leq R_1$. Given $\epsilon>0$ and $k_0>0$ there exists $C$ such that, for
all $R_2 > (1 + \epsilon) R_1$, $k \geq k_0$, $s$, $\bsy_s\in U_s$, and $f
\in L^2(D_{R_0})$, with $\varphi_{\rm PML}$ as in
\eqref{eq:sigma_prop}, if $h = h(k)$ satisfies \eqref{eq:FEMthreshold}
then
\begin{align} \label{eq:estG}
 &\vert G(u_s - u_{s,{\rm PML},h})(\cdot, \bsy_s) \vert \nonumber  \\
 &\leq  C\,  k^{r-1}\,  K_r(f,G)\,
 \Big[\exp
 \big(-C_{{\rm PML},2}\, k\,\big(R_2 - (1+\epsilon) R_1\big) \big)
 +  \big( (hk)^{\ell+1-r} + (hk)^{2\ell}\, k\, R_2\big)\, R_2\, \Big], 
\end{align}
with $K_r(f,G)$ defined in \eqref{eq:def-K} and $C_{{\rm PML},2}$ as in
Theorem~\ref{thm:GLS}.
\end{corollary}

\begin{proof}
We prove \eqref{eq:estG} for $r\in \{0,1\}$. The result then follows for
all $r \in [0,1]$ by interpolation (see, e.g., \cite[\S 14]{BrSc:08}). We
write
\begin{align} \label{eq:estG1}
 \vert G(u_s - u_{s,{\rm PML},h})(\cdot, \bsy_s) \vert
 \,\le\, \Vert u_s(\cdot, \bsy_s) - u_{s,{\rm PML},h}(\cdot , \bsy_s) \Vert _{H^r(D_R)}\,
  \Vert G \Vert_{H^r(D_R)'}\,.
\end{align}
Then, we use the triangle inequality
\begin{align}
& \Vert u_s(\cdot, \bsy_s) - u_{s,{\rm PML},h}(\cdot , \bsy_s) \Vert _{H^r(D_R)}\nonumber  \\
& \leq \ \Vert u_s(\cdot, \bsy_s) - u_{s,{\rm PML}}(\cdot , \bsy_s) \Vert _{H^r(D_R)}
+  \Vert u_{s, {\rm PML}}(\cdot, \bsy_s) - u_{s,{\rm PML},  h}(\cdot , \bsy_s) \Vert _{H^r(D_R)}\,.  \label{eq:triangle}
\end{align}
Using the norm bound $\|v\|_{L^2(D_R)} \le k^{-1}\,\|v\|_{H^1_k(D_R)}$
when $r=0$, and $\|v\|_{H^1(D_R)} \le (1 + k_0^{-2})^{1/2}\,
\|v\|_{H^1_k(D_R)}$ for $k\ge k_0$ when $r=1$, the first term in
\eqref{eq:triangle} is estimated by \eqref{eq:GLS} noting Remark
\ref{rem:alls}. To estimate the second term in \eqref{eq:triangle}, we use
\eqref{eq:FEMrel1} when $r=0$, and \eqref{eq:FEMrel2} when $r=1$, and then
use \eqref{eq:inherit1} to estimate $\Vert u_{s,{\rm PML}}(\cdot, \bsy_s)
\Vert_{H^1_k(D_R)}$.
\end{proof}

\section{Sound-soft plane-wave scattering}
\label{sec:SSSP}

Returning to the plane wave scattering problem introduced in \S
\ref{subsec:SSSP}, we now show that this can be reformulated as a
particular case of the EDP \eqref{eq:intro_tedp}--\eqref{eq:src}. Then we
provide a computable formula for the far-field pattern $u_\infty^{S}$
defined in \eqref{eq:FFP},  as well as an error estimate for computing
this via PML truncation and FEM approximation. The error estimate is
uniform in the random parameter $\bsy$ and therefore leads to an analogous
estimate for approximation of the expected value of $u^{S}_\infty$,
allowing us to estimate the quantity  $\EFEM$ (discussed in \S
\ref{sec:intro}) in the particular case when $G$ denotes the far-field
pattern operator.

\subsection{The plane-wave scattering problem formulated as an EDP}
\label{subsec:plane_wave_EDP}

While both $u^T$ (the total field) and $u^S$ (the scattered field) solve
problem \eqref{eq:intro_tedp}, neither satisfy the whole EDP problem
\eqref{eq:intro_tedp}--\eqref{eq:src}; ~$u^T$ satisfies
\eqref{eq:intro_dbc}, but not \eqref{eq:src}, while  $u^S$ satisfies
\eqref{eq:src} but not \eqref{eq:intro_dbc}. However, by introducing a
cutoff function $\varphi_{\rm alt}$ with value $1$ in a neighbourhood
of $D$ and vanishing in the far field, and defining
 \beq\label{eq:tilde_u}
 \ualt : = u^S + \varphi_{\rm alt}\, u^I = u^T- (1-\varphi_{\rm alt})\,u^I,
 \eeq
we see that $\ualt$ satisfies both \eqref{eq:intro_dbc} and \eqref{eq:src}
and thus it satisfies an EDP with a special source term $\falt$ (see Lemma
\ref{lem:L2}). Then one can immediately apply the bounds in Lemma
\ref{lem:boundsDtN}, and also approximate $\ualt$ via PML truncation and
FEM as in \S \ref{sec:FEM}.

We construct $\varphi_{\rm alt}$ so that it has value $1$ on the set
$\supp(A-I) \cup \supp(n-1) \cup \overline{D}$: this then ensures that
$\falt$ has a simple form and also leads to a simple formula for the
quantity $u^{S}_\infty$ (see Theorem~\ref{thm:Monk}).

In more detail, by definition of $R_0$ in Assumption~\ref{ass:1}, there
exists $\eta>0$ such that
\begin{align} \label{eq:exist_eta}
 \supp(A-I) \cup \supp(n-1) \cup \overline{D}\, \subseteq \,  B_{R_0-\eta}.
\end{align}
Accordingly, we introduce a univariate smooth cutoff function
$\varphi_{\rm alt} \in C^2(\mathbb{R}_+)$ satisfying
 \beq\label{eq:psi_pw1}
 \varphi_{\rm alt}(r) = 1 \quad \text{for} \quad 0 \leq r \leq R_0-\eta
 \quad \text{and} \quad \varphi_{\rm alt}(r) = 0 \quad \text{for} \quad r
\geq R_0 - \eta/2. \eeq Then (with an obvious abuse of notation) we extend
this to a radial function $\varphi_{\rm alt} \in C^{2}(\mathbb{R}^d)$
with values in $[0,1]$,  by setting
\begin{align} \label{eq:psi_pw}
 \varphi_{\rm alt}(\bsx) =
  \varphi_{\rm alt}(|\bsx|_2).
\end{align}

\begin{figure}[t]
\begin{center}
\begin{tikzpicture}[scale=0.33]
 \draw[ultra thick,fill=white] (0,0) circle (8);   
 \shade[even odd rule,ring shading={from blue!30!white at 5 to white at 6.5}] (0,0) circle (5) circle (6.5);
 \draw[ultra thick,red,fill=blue!30!white] (0,0) circle (5); 
 \draw[ultra thick,red,dashed] (0,0) circle (6.5); 
 \draw[ultra thick, pattern={Lines[angle=45,distance=8pt]}, blue, fill=white]
   plot[smooth cycle, tension=.7]
   coordinates {(-2,-2) (-1,2) (3,2.1) (3.5,-1) (1.5,-1.5)};
 \node[right] at (0,0) {$D$};
 \node[blue] at (2,1.5) {$\partial D$};
 \node[left] at (10.6,0) {$B_{R_0}$};
 \draw[red,-latex] (8,-6) -- (6.6,0); \node[red,below] at (11,-6) {radius $R_0-\eta/2$};
 \draw[red,-latex] (6,-7.4) -- (5.1,0); \node[red,below] at (8,-7.4) {radius $R_0-\eta$};
\end{tikzpicture}
\end{center}
\caption{When formulating the plane-wave sound-soft scattering problem as an EDP,
the blue gradient shading illustrates the radial cutoff function $\varphi_{\rm alt}\in C^2(\bbR^d)$
satisfying \eqref{eq:psi_pw1}. The support of the data $f^{\rm alt}$ as defined in \eqref{eq:pw_f}
is in the annulus between radii $R_0-\eta$ (solid red circle) and $R_0-\eta/2$ (dashed red circle).
}\label{fig:SSSP}
\end{figure}
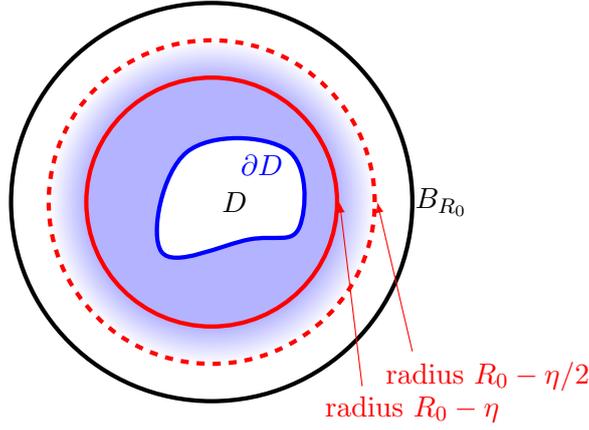

Note that, without loss of generality, we may assume that $R_0$ in
Definition \ref{ass:1} is chosen sufficiently far from the boundary of the
set on the left-hand side of \eqref{eq:exist_eta} so that
$\varphi_{\rm alt}$ does not have an unpleasantly-large gradient.
This is illustrated in Figure~\ref{fig:SSSP}.

\begin{lemma}[Plane-wave sound-soft scattering rewritten as an EDP with $L^2$ data]\label{lem:L2}
Let $u^T$ be the solution of the plane-wave sound-soft scattering problem
of Definition \ref{def:SSSP}. Then $\ualt$ defined by \eqref{eq:tilde_u}
and \eqref{eq:psi_pw} satisfies the EDP: \beq\label{eq:pw_f} \nabla \cdot(
A\nabla\ualt) +k^2 n \ualt = \falt :=
  2 \nabla \varphi_{\rm alt}\cdot \nabla u^I + u^I \Delta \varphi_{\rm alt} \quad \text{in } \quad D_+ \, .
\eeq together with \eqref{eq:intro_dbc} and   \eqref{eq:src}. Moreover,
with $R_0>0$ as in Assumption~\ref{ass:1}, $\ualt \equiv u^T $ on
$D_{R_0-\eta}$ and $\ualt \equiv u^S$ on $\mathbb{R}^{d}\setminus
B_{R_0-\eta/2}$, while $\falt$ is only nonzero on the annulus
$B_{R_0-\eta/2}\setminus B_{R_0-\eta}$.
\end{lemma}

\bpf Given $\varphi_{\rm alt}$ as in \eqref{eq:psi_pw}, the function
$\ualt =  u^T- (1-\varphi_{\rm alt})u^I$ satisfies both the
Sommerfeld radiation condition \eqref{eq:src} and
\begin{align*}
(\nabla\cdot(A\nabla) +k^2 n)
\ualt
&= -(\nabla\cdot(A\nabla) +k^2 n) (1-\varphi_{\rm alt}) u^I \\
&=-(\Delta+k^2) (1-\varphi_{\rm alt}) u^I
= 2\nabla \varphi_{\rm alt}\cdot \nabla u^I + u^I \Delta \varphi_{\rm alt},
\nonumber
\end{align*}
where we used the fact that $A\equiv I$ and $n\equiv 1$ on
$\supp(1-\varphi_{\rm alt})$ and also the fact that $(\Delta +k^2)u^I
= 0$.
\epf

\subsection{The far-field pattern of $u^S$}

The goal of this section is to express the far-field pattern of $u^S$,
\eqref{eq:FFP},  as a functional of $\ualt$, with the result being
\eqref{eq:Monk2} below. The numerical experiments in \S
\ref{sec:Numerical} then compute this quantity of interest.

We first show in Theorem~\ref{thm:Monk} how the far-field pattern of any
solution $v$ of $(\Delta +k^2)v=0$ satisfying \eqref{eq:src} can be
obtained from knowledge of $v$ in an annulus. Then in
Corollary~\ref{cor:FFPuS} we apply this technique to obtain a formula for
the far-field pattern of $u^S$. For convenience, the result is stated for
the deterministic case, under Assumptions~\ref{ass:1} and~\ref{ass:2},
although the random case is entirely analogous.

\begin{theorem}[Expression for far-field pattern as integral over an annulus]\label{thm:Monk}
Suppose that for some $\rho_0>0$, $v\in C^2(\mathbb{R}^d\setminus
B_{\rho_0})$ satisfies the Helmholtz equation $(\Delta +k^2)v=0$ in
$\mathbb{R}^d\setminus \overline{B_{\rho_0}}$,  together with the
Sommerfeld radiation condition \eqref{eq:src}. Let $\rho_1 > \rho_0$ and
suppose $\varphi_{\rm ffp} \in C^2(\mathbb{R}^d)$ is a radial
cutoff function with $\varphi_{\rm ffp} \equiv 0$ in a
neighbourhood of $B_{\rho_0}$ and $\varphi_{\rm ffp} \equiv 1$ on
$\mathbb{R^d} \backslash B_{\rho_1}$. Then
  \beq\label{eq:Monk}
  v_\infty(\widehat{\bsx}) \,=\, c(d,k)
  \int_{\supp \nabla \varphi_{\rm ffp} } v(\bsx')\, \big( \Delta \varphi_{\rm ffp} (\bsx') - 2\, \ri\, k\,
  \widehat{\bsx}\cdot\nabla\varphi_{\rm ffp} (\bsx') \big)\,
  \exp\big(-\ri\, k\, \bsx' \cdot \widehat{\bsx}\big) \, {\rm
    d} \bsx', \quad \widehat{\bsx} \in \mathbb{S}^{d-1}, \eeq
    where
    \beq\label{eq:Cdk} c(2,k) =\frac{1}{\sqrt{k} } \frac{e^{\ri
        \pi/4}}{2\sqrt{2\pi}} \quad\text{ and }\quad c(3,k) = \frac{1}{4\pi}.
       \eeq
\end{theorem}

Theorem~\ref{thm:Monk} is essentially \cite[Theorem 2.2]{Mo:95}, but is
written in a slightly-more general way here.

\bpf For $\bsx, \bsx' \in \mathbb{R}^d$,
let \beq \label{eq:fund}
 G(\bsx,\bsx'):= \frac{\ri}{4} H_0^{(1)}(k\,
 |\bsx-\bsx'|_{2}) \,\,\text{ for } d=2, \quad \text{ and }\quad
 G(\bsx,\bsx')
:= \frac{\exp(\ri\, k\, |\bsx-\bsx'|_{2})}{4\pi\,
  |\bsx-\bsx'|_{2}} \,\,\text{ for } d=3, \eeq
where $H^{(1)}_0$ denotes the Hankel function of the first kind of order zero, see e.g.,
\cite[(3.83)]{CoKr:13}.
Recall that $G$ is the
fundamental solution for the operator $-(\Delta + k^2)$,  i.e.,
$$(\Delta_{\bsx'} + k^2)\, G(\bsx,\bsx') \,=\, - \delta (\bsx-\bsx') \quad
\text{for all} \quad \bsx , \bsx' \in \mathbb{R}^d,$$ with $\delta$
denoting the Dirac delta. Using this, for any $R > 0$, any $w\in
C^2(\overline{B_R})$ can be represented (using Green's integral
representation theorem -- see, e.g., \cite[Theorem 2.1]{CoKr:13}) as:
 \beq \label{eq:Green}
 w(\bsx) = \int_{\partial B_R}\left( \frac{\partial w(\bsx')}{\partial
 \nu(\bsx')} \,
G(\bsx,\bsx') - w(\bsx') \frac{\partial G(\bsx,\bsx')}{\partial
\nu(\bsx')}\right) \rd S(\bsx') \,-\, \int_{B_R} G(\bsx,\bsx')\, ((\Delta
+ k^2)w)(\bsx')\, \rd \bsx',  \eeq for  $\bsx \in B_R$,  where $\partial/\partial\nu(\bsx')$
denotes the outward normal derivative with respect to $\bsx'$.

Now, with $v$ and $\varphi_{\rm ffp}$ as given in the hypotheses, let
$R > \rho_1$, consider the function $w$ defined by $w\equiv
v\,\varphi_{\rm ffp}$ on $\mathbb{R}^d\backslash
\overline{B_{\rho_0}}$ and $w \equiv 0$ on $B_{\rho_0}$. Then (since
$\varphi_{\rm ffp} \equiv 0$ on a neighbourhood of $B_{\rho_0}$) we
have $w\in C^2(\overline{B_R})$. With this choice of $w$, as $R\to
\infty$, the integral over $\partial B_R$ in \eqref{eq:Green} tends to
zero. This is because both $v$ and $G(\bsx,\cdot)$ (for $\bsx$ fixed)
satisfy the Sommerfeld radiation condition and both $v$ and $G(\bsx,
\cdot)$ decay with $\mathcal{O}(r^{(1-d)/2})$ as $r \rightarrow \infty$.
(Similar arguments are used in, e.g., \cite[Last equation in the proof of
Theorem 3.3]{CoKr:83}). Moreover,
 \beq\label{eq:Helmw} (\Delta +k^2) w
\,=\, (\Delta +k^2)(\varphi_{\rm ffp}  v) \,=\, 2
\nabla\varphi_{\rm ffp} \cdot\nabla v + v \Delta \varphi_{\rm ffp}  + \varphi_{\rm ffp} (\Delta +k^2)v \eeq Therefore inserting
\eqref{eq:Helmw} into \eqref{eq:Green} and letting $R \rightarrow \infty$,
we obtain \beqs v(\bsx) \,=\, -\int_{\supp \nabla\varphi_{\rm ffp} }
 G(\bsx,\bsx')\, \big(2 \nabla\varphi_{\rm ffp}  \cdot\nabla v + v \Delta
 \varphi_{\rm ffp}\big)(\bsx')\, 
 \rd \bsx', \quad \bsx \in \mathbb{R}^d \backslash \overline{B_{\rho_1}}. 
 \eeqs
where we used the fact that $\nabla \varphi_{\rm ffp} (\bsx') =
\Delta \varphi_{\rm ffp}  (\bsx') =0 $ for $\bsx' \not \in
\mathrm{supp} \nabla \varphi_{\rm ffp}$ and $(\Delta + k^2 )v(\bsx')
= 0 $ for $\vert \bsx' \vert_{2} > \rho_0$.

Now, using the definition \eqref{eq:FFP} of the far-field pattern and
(when $d=2$) using the large-argument asymptotics of the Hankel function $H_0^{(1)}$
(see, e.g., \cite[Equation 10.17.5]{Di:22}) or (when $d=3$) via an
elementary calculation, we find that \beqs v_\infty(\widehat{\bsx}) \,=\,
-c(d,k) \int_{\supp \nabla\varphi_{\rm ffp} } \big(2
\nabla\varphi_{\rm ffp} \cdot\nabla v + v \Delta \varphi_{\rm ffp} \big)(\bsx')\,
\exp\big(-\ri\, k\, \bsx' \cdot \widehat{\bsx}\big)
\,\rd \bsx'. \eeqs The result \eqref{eq:Monk} then follows by integrating
by parts (via the divergence theorem) the term involving
$\nabla\varphi_{\rm ffp} \cdot\nabla v$ (moving the derivative from
$v$ onto $\varphi_{\rm ffp}$). \epf

Recall that our goal is to  find the far-field pattern of $u^S$ from
knowledge of $u^{\rm alt}$ defined by \eqref{eq:tilde_u}. Recall also that
$u^{\rm alt} \equiv u^S$ on $\mathbb{R}^d\setminus B_{R_0-\eta/2}$, and
that we will ultimately compute an approximation to $u^{\rm alt}$ on
$D_{R_0}$ via PML truncation. We therefore choose the annulus $\supp
\nabla\varphi_{\rm ffp}$ in Theorem \ref{thm:Monk} to be
$B_{R_0}\setminus B_{R_0-\eta/2}$, since in this region we will have an
approximation to $u^{\rm alt} \equiv u^S$.

\begin{corollary}[Far-field pattern for $u^S$ as an integral involving $u^{\rm alt}$]\label{cor:FFPuS}
Given $\eta>0$ as in \eqref{eq:exist_eta}, let $\varphi_{\rm ffp} \in
C^2(\mathbb{R}^d)$ be such that \beq\label{eq:part_psi}
\varphi_{\rm ffp} \equiv 1 \text{ on } \mathbb{R}^d \setminus B_{R_0}
\quad\tand\quad \varphi_{\rm ffp} \equiv 0 \text{ on }
B_{R_0-\eta/2}. \eeq Then, with $u^{\rm alt}$ defined by
\eqref{eq:tilde_u} and \eqref{eq:psi_pw},
  \beq\label{eq:Monk2a}
  u^S_\infty(\widehat{\bsx}) \,=\, c(d,k)
  \int_{\supp \nabla\varphi_{\rm ffp}} u^{\rm alt}(\bsx')\, \big( \Delta \varphi_{\rm ffp}(\bsx') - 2\, \ri\,
  k\,
  \widehat{\bsx}\cdot\nabla\varphi_{\rm ffp}(\bsx') \big)\,
  \exp\big(-\ri\, k\, \bsx' \cdot \widehat{\bsx}\big) \, \rd \bsx'. \eeq
\end{corollary}

\bpf The conditions in \eqref{eq:part_psi} imply that $\supp \nabla
\varphi_{\rm ffp}\subset B_{R_0} \setminus B_{R_0 - \eta/2}$. That
is, $\nabla \varphi_{\rm ffp}$ is only supported where $u^{\rm alt}
\equiv u^S$ by Lemma \ref{lem:L2}. The result then follows from Theorem
\ref{thm:Monk}. \epf

\subsection{PML-FEM approximation}
\label{subsec:PMLFEM}

In this subsection we obtain an error estimate for the approximation of
$u^S_\infty $, as given by the formula \eqref{eq:Monk2a}. We apply
dimension truncation (replacing $\bsy$ by $\bsy_s$, as in
\S\ref{sec:dimtrunc}) and then combined PML-FEM approximation (as in
\S\ref{sec:FEM}), to the EDP \eqref{eq:pw_f}, thus obtaining  an
approximation  $\ualt_{s,{\rm PML}, h}$, This is then substituted into
\eqref{eq:Monk2a} in place of $\ualt$, yielding the far-field
approximation which we call $u^S_{\infty,s,{\rm PML}, h}$.  We note that
source  $\falt$ in \eqref{eq:pw_f} is oscillatory (in the sense of
\eqref{eq:koscillatory}), allowing  us to use the estimates from Theorem
\ref{thm:FEM1}, in particular the $L^2$ bound in \eqref{eq:FEMrel2}.

In the following theorem we give the PML-FEM error estimate for the far-field pattern without truncation of the variable $\bsy$ (i.e., $s = \infty$). However the case of finite $s$ is identical 
as explained in Remark~\ref{rem:alls}.

\begin{theorem}[FEM approximation of the far-field pattern for plane-wave
scattering] \label{thm:FEM-ffp}
Let $D$, $A$, $n$, and $R_0$ satisfy Assumptions~\ref{ass:1},
\ref{ass:3}~and~\ref{ass:4}, and let $R_1> R_0$. Suppose $A$ and $n$ also
satisfy Assumption~\ref{ass:smoothness} with smoothness degree $\tau$. Let
$(V_h)_{h>0}$ satisfy Assumption~\ref{ass:spaces} with approximation
degree $m$. For $\ell :=\min(\tau,m)$
suppose that $\partial D_{R_2}$ is H\"older continuous $C^{\ell,1}$ and
that $\varphi_{\rm PML}$ as defined in \eqref{eq:sigma_prop} is in
$C^{\ell,1}(\bbR^+)$. 
Let $u^{\rm alt}_{\rm PML}$ be the PML approximation of
$u^{\rm alt}$, i.e., the solution of the PML variational problem
\eqref{eq:PMLvf} with $f$ replaced by $f^{\rm alt}$ defined by
\eqref{eq:pw_f}. Let $u^{\rm alt}_{{\rm PML}, h}$ be the FEM approximation
to $u^{\rm alt}_{\rm PML}$, i.e., the solution of the variational problem
\eqref{eq:FEM} with $f$ replaced by $f^{\rm alt}$. Let
\begin{align} \label{eq:Monk2}
  &u^S_{\infty,{\rm PML},h}(\widehat{\bsx}) \\ 
  &\,:=\, c(d,k)
  \int_{\supp \nabla\varphi_{\rm ffp}} u^{\rm alt}_{{\rm PML}, h}(\bsx')\,
  \big( \Delta \varphi_{\rm ffp}(\bsx') - 2\,\ri\, k\,
  \widehat{\bsx}\cdot\nabla\varphi_{\rm ffp}(\bsx') \big)\,
  \exp\big(-\ri\, k\, \bsx' \cdot \widehat{\bsx}\big) \,\rd \bsx',
  \quad \widehat{\bsx} \in \mathbb{S}^{d-1}, \nonumber
\end{align}
with $\varphi_{\rm ffp}$ satisfying \eqref{eq:part_psi}.
Given $\epsilon >0$, there exist constants $C>0$
and $k_0>0$ such that for all $R_{2}>(1 + \epsilon)R_1$, $k\geq k_0$,
$\bsy \in U$, and for all incident directions $\bsalpha\in\bbR^d$ with
$\vert \bsalpha \vert_{2} = 1$, if $h = h(k)$ satisfies
\eqref{eq:FEMthreshold} then
\begin{align} \label{eq:FFFEM}
 & \big\|u^S_\infty-u^S_{\infty,{\rm PML},h}\big\|_{L^\infty(\mathbb{S}^{d-1})} \nonumber \\
 & \leq c(d,k)\,C \Big[
 \exp \big(-C_{{\rm PML},2}\, k\,(R_2 - (1+\epsilon) R_1) \big)
+ \big( (hk)^{\ell+1} + (hk)^{2\ell}\, k\, R_2 \big) \,R_2 \Big]
\|f^{\rm alt}\|_{L^2(D_{R_0})},
\end{align}
where $c(d,k)$ is given by \eqref{eq:Cdk}, $C_{{\rm PML},2}$ is as in
Theorem~\ref{thm:GLS}, and $\|f^{\rm alt} \|_{L^2(D_{R_0})} \lesssim k$.
\end{theorem}

The bound \eqref{eq:FFFEM} implies that if
\[
  c(d,k)\,(hk)^{2\ell}\,k^2\,R_2^2
\]
is a sufficiently small constant, then the error in computing the
far-field pattern via \eqref{eq:Monk2} is bounded independently of $k$.
From \eqref{eq:Cdk} we see that the condition for $d=3$ is
$(hk)^{2\ell}\,k^2\,R_2^2$ being sufficiently small, and for $d=2$ it is
(slightly weaker) $(hk)^{2\ell}\,k^{3/2}\,R_2^2$ being sufficiently small.

\bpf By the definitions of $u^S_\infty$ and $u^S_{\infty,{\rm PML},h}$
together with Cauchy--Schwarz inequality, there exists $C'>0$ (independent
of $\epsilon, \bsy, k, \bsalpha$, and $h$) such that
\begin{align}\nonumber
 &\big\|u^S_\infty-u^S_{\infty,{\rm PML},h}\big\|_{L^\infty(\mathbb{S}^{d-1})}
 \leq c(d,k)\,C'\,k\, \big\| u^{\rm alt} - u^{\rm alt}_{{\rm PML}, h}\big\|_{L^2(\supp \nabla \varphi_{\rm ffp})}\\
 &\leq c(d,k)\,C'\Big( k\, \big\| u^{\rm alt} - u^{\rm alt}_{{\rm PML}}\big\|_{L^2(\supp \nabla \varphi_{\rm ffp})}
 + k\,\big\| u^{\rm alt}_{\rm PML} - u^{\rm alt}_{{\rm PML}, h}\big\|_{L^2(\supp \nabla \varphi_{\rm ffp})}\Big).\label{eq:temp0}
\end{align}
By the exponential accuracy of the PML approximation \eqref{eq:GLS},
\begin{align} \label{eq:temp1}
 k\, \big\| u^{\rm alt} - u^{\rm alt}_{{\rm PML}}\big\|_{L^2(\supp \nabla \varphi_{\rm ffp})}
 \leq C_{{\rm PML},1}\, \exp \big(-C_{{\rm PML},2}\, k\,(R_2 - (1+\epsilon) R_1) \big)\,
 \|f^{\rm alt}\|_{L^2(D_{R_0})}.
\end{align}
By the FEM error estimate \eqref{eq:FEMrel2}, the fact that $f^{\rm alt}$
is $k$-oscillating in the sense of \eqref{eq:koscillatory}, and the bound
\eqref{eq:inherit1},
\begin{align} \label{eq:temp2}
 k\, \big\| u^{\rm alt}_{\rm PML} - u^{\rm alt}_{{\rm PML}, h}\big\|_{L^2(\supp \nabla \varphi_{\rm ffp})}
 \,\leq\, C''\,\big( (hk)^{\ell+1} + (hk)^{2\ell}\, k\, R_2\big)\,
  \, R_2\,\|f^{\rm alt}\|_{L^2(D_{R_0})}.
\end{align}
The result \eqref{eq:FFFEM} then follows from inputting \eqref{eq:temp1}
and \eqref{eq:temp2} into \eqref{eq:temp0}. Note that we have $\|f^{\rm
alt} \|_{L^2(D_{R_0})} \lesssim k$ by the definitions of $f^{\rm alt}$
\eqref{eq:pw_f} and $u^I$ \eqref{eq:incident}. \epf

\section{Numerical Experiments}
\label{sec:Numerical}

In this section we illustrate the performance of our algorithm for
computing the far-field pattern of the scattered field in the case of
$2$-dimensional plane-wave scattering by a random medium with an
impenetrable obstacle.
To verify our implementation, we computed the solution for a Dirac impulse at the origin with no obstacle, for which the exact solution is known. We also computed the far-field pattern with a triangular
obstacle with a homogeneous field  ($A = I$, $n= 1$) and obtained good agreement with the results for
this case from the authors of \cite{HLM13}. 
The paper \cite{HLM13} uses a much more efficient and accurate boundary element method, but one suitable only for the homogeneous case, see also \cite{GibLan24}.

\subsection{Problem specification}

For $d=2$ we consider a butterfly-shaped obstacle defined by the smooth radial function
$(0.3 + \sin^2(\theta))(1.5+1.4\cos(2\theta))$ for $\theta\in [0,2\pi]$.
The boundary of the obstacle is contained in the annulus between radii $0.130$ and
$1.249$. For our domain we consider the annuli given by
\begin{align*}
  \mbox{radius} &= 1.25 && (\mbox{contains a butterfly-shaped obstacle}), \\
  R_0 - 3\eta/2 &= 3 && (\mbox{random field fluctuation starts to diminish}), \\
  R_0 - \eta &= 3.5 && (\mbox{random field fluctuation ends}; \mbox{support of $f^{\rm alt}$ starts}), \\
  R_0 - \eta/2 &= 4 && (\mbox{support of $f^{\rm alt}$ ends}; \mbox{far-field pattern calculation region starts}), \\
  R_0 &= 4.5  && (\mbox{far-field pattern calculation region ends}), \\
  R_1 &= 4.52 && (\mbox{PML starts}), \\
  R_2 &= 5 && (\mbox{PML ends}).
\end{align*}
The incident wave comes from above, i.e., at angle $90^\circ$ and $\bsalpha = (0,-1)$ in \eqref{eq:incident}.
We need a number of cutoff functions with varying smoothness requirements
as we specify below.

The PML takes place between radii $R_1=4.52$ and $R_2=5$. The function
$\varphi_{\rm PML}\in C^3(\bbR_+)$ satisfying \eqref{eq:sigma_prop} is
taken as (with a scaling parameter $\varphi_{\rm PML,const}:=3$)
\[
 \varphi_{\rm PML}(r) :=
 3\, \varphi_{{\rm cut}, 3}\Big(\frac{r-R_1}{R_2-R_1}\Big),
 \quad\mbox{with}\quad
 \varphi_{{\rm cut}, 3}(r) :=
 \begin{cases}
 0 & \mbox{for } r\le 0,  \\
 35r^4-84r^5+70r^6-20r^7 & \mbox{for } 0\le r\le 1, \\
 1 & \mbox{for } r\ge 1.
 \end{cases}
\]

The active region of $f^{\rm alt}$ defined in \eqref{eq:pw_f} for
plane-wave scattering is between radii $R_0 - \eta = 3.5$ and $R_0 -
\eta/2 = 4$ (thus $\eta = 1$), while the far-field pattern is calculated
between radii $R_0-\eta/2 = 4$ and $R_0 = 4.5$. The functions
$\varphi_{\rm alt}\in C^2(\bbR^d)$ satisfying
\eqref{eq:psi_pw1}--\eqref{eq:psi_pw} and $\varphi_{\rm ffp}\in
C^2(\bbR^d)$ satisfying \eqref{eq:part_psi} are taken respectively as
\[
 \varphi_{\rm alt}(\bsx) := 1 - \varphi_{{\rm cut,2}}\Big(\frac{|\bsx|_2-(R_0 - \eta)}{\eta/2}\Big)
 \quad\mbox{and}\quad
 \varphi_{\rm ffp}(\bsx) := \varphi_{{\rm cut,2}}\Big(\frac{|\bsx|_2-(R_0 - \eta/2)}{\eta/2}\Big),
\]
with
\[
 \varphi_{{\rm cut}, 2}(r) :=
 \begin{cases}
 0 & \mbox{for } r\le 0,  \\
 10r^3 - 15r^4 + 6r^5 & \mbox{for } 0\le r\le 1, \\
 1 & \mbox{for } r\ge 1.
 \end{cases}
\]
We compute the far-field pattern \eqref{eq:Monk2} at directions $\widehat{\bsx}$ corresponding to integer degrees $1^\circ, 2^\circ, \ldots, 360^\circ$. 

We restrict here to the case of \eqref{eq:PWSS} with fixed homogeneous
matrix $A(\bsx,\bsy)=I$ (therefore $A_0(\bsx) = I$, $\xi_A = 0$,
$\Psi_j(\bsx) = 0$ in \eqref{eq:Axy}, and we may take $\mu_A = 1$ in
\eqref{eq:hh-Acond}), and we take the coefficient $n(\bsx,\bsy)$ as in
\eqref{eq:nxy}, with
\begin{align} \label{eq:nspec}
  n_0(\bsx) = 1,
  \quad
  \xi_n = 0.8319, 
  \quad
  \psi_j(\bsx) = j^{-q} \sin(j \pi x_1)\, \sin(j\pi x_2)\, \varphi_{\rm fluc}(\bsx),
  \quad
  q = 3,
\end{align}
where
\[
 \varphi_{\rm fluc}(\bsx) := 1 - \varphi_{{\rm cut},1}\Big(\frac{|\bsx|_2-(R_0 - 3\eta/2)}{\eta/2}\Big),
 \quad\mbox{with}\quad
 \varphi_{{\rm cut}, 1}(r) :=
 \begin{cases}
 0 & \mbox{for } r\le 0,  \\
 3r^2-2r^3 & \mbox{for } 0\le r\le 1 \\
 1 & \mbox{for } r\ge 1.
 \end{cases}
\]
The cutoff function $\varphi_{\rm fluc} \in C^1(\bbR^d)$ is chosen to
ensure that $\supp (\psi_j) \subseteq B_{R_0-\eta}$ as required in
\eqref{eq:exist_eta}, as well as $\|\nabla \psi_j\|_{L^\infty(D_{R_0})} <
\infty$ as required in \eqref{eq:npsimeaspos}. Consequently, Assumption~\ref{ass:smoothness} holds with smoothness degree for the random coefficients being $\tau = 1$.

The value of $\xi_n$ is chosen so that the condition \eqref{eq:npsimeaspos} is satisfied. The condition \eqref{eq:npsimeaspos} is sufficient to ensure positivity of
$n(\bsx,\bsy)$. There is another condition \eqref{eq:nontrap} which is sufficient to
ensure nontrapping, but in our experiment we will not impose this latter condition, see discussion in \S\ref{sec:trapping}.
Indeed, with $n_0(\bsx) = 1$ we have $n_{0,\min} =
1$, $n_{\min} = 1/2$, $n_{\max} = 3/2$, and we may take $\mu_n = 1$ in
\eqref{eq:hh-Acond}. The definition of $\psi_j$ in \eqref{eq:nspec} yields
\begin{align*}
  \sum_{j=1}^\infty \|\psi_j\|_{L^\infty(D_{R_0})}
  &\le \sum_{j=1}^\infty j^{-q} \le 1.2021,
  \\
  \sum_{j=1}^\infty \esssup_{\bsx \in D_{R_0}} \big|\psi_j(\bsx) + \bsx\cdot\nabla\psi_j(\bsx)\big|
  &\le \sum_{j=1}^\infty \esssup_{\bsx \in D_{R_0}}
  j^{-q} \bigg(1 + \sum_{i=1}^d |x_i|
  \bigg(j\pi + \bigg|\frac{\partial\varphi_{\rm fluc}(\bsx)}{\partial x_i}\bigg|\bigg) \bigg) \\
  &\le \sum_{j=1}^\infty j^{-q} \Big(1 + 2 (R_0 - \eta)(j\pi + 3) \Big)
  = \sum_{j=1}^\infty j^{-q} (7j\pi + 22)
  \le 62.6193,
\end{align*}
where we used $|x_i|\le R_0-\eta$ since $\supp (\psi_j) \subseteq
B_{R_0-\eta}$, and
\[
  \bigg|\frac{\partial\varphi_{\rm fluc}(\bsx)}{\partial x_i}\bigg|
  = \bigg|(-6r+6r^2)\Big|_{r = \frac{|\bsx|_2 - (R_0-3\eta/2)}{\eta/2}}
  \cdot\frac{x_i}{|\bsx|_2}\cdot\frac{1}{\eta/2}\bigg|
 \le \frac{3}{2}\times 1\times \frac{1}{\eta/2} = 3.
\]
Thus, to satisfy \eqref{eq:npsimeaspos} it suffices to take 
$\xi_n\le n_{0,\min}/1.2021 \approx 0.8319$, 
while \eqref{eq:nontrap} would ask for the much more restrictive
$\xi_n\le \mu_n / 62.6193 \approx 0.01597$. As indicated above, we
take $\xi_n = 0.8319$.

Below we will consider $k\in \{12, 24, 48\}$ so we may take $k_0
= 12$. This yields from~\eqref{eq:Cstab}
\begin{align*}
  \Cstab
  =
  \frac{1}{\min\{1,1/2\}} \bigg( \frac{1}{12\times4.5} +
    4 \sqrt{ \frac{1}{1/2} + \frac{1}{1/2}
    \bigg( 1 + \frac{2-1}{2\times 12\times4.5}\bigg)^2}
    \frac{3/2}{\sqrt{\min\{1/2, 1/2\}}}
    \bigg)
    \le 48.26.
\end{align*}
The wavelength is $2\pi/k \le 2\pi/k_0 = 2\pi/12 \le 0.5236$.
The radii listed at the start of this subsection were chosen to be roughly $0.5$ apart so that we have around one or more waves in each annulus.

\subsection{PML-FEM discretization}

We use the Julia packages Gridap and GridapEmbedded \cite{Badia2020,Verdugo2022,BNV22} to discretize our domain and solve the PDE using an aggregated unfitted finite element method (AgFEM) \cite{BNV22}. Specifically, we divide the square $[-R_2,R_2]^2 = [-5,5]^2$ into $800^2$ square cells each divided into two triangles, resulting in more than a million triangles in the background mesh, with
\[
  h = \frac{2\,R_2}{800} = \frac{10}{800} = 0.0125.
\]
A circle of radius $R_2$ is extracted from this background mesh and then the butterfly-shaped obstacle is cut out from the centre. Both boundaries are ``refined'' by triangles and then the AgFEM method is used to solve the PDE. The Dirichlet boundary conditions are enforced in a weak form and we use a ``Nitsche factor'' (see \cite{BNV22}) of $\gamma / h$ with $\gamma = 100$.

We set the approximation degree of AgFEM to be $m=2$ (cf.~Assumption~\ref{ass:spaces}). Together with the smoothness degree of the random field being $\tau=1$ (cf.~Assumption~\ref{ass:smoothness}), we see from Theorem~\ref{thm:FEM1} that, with $\ell = \min(\tau,m)=1$, if $h = h(k)$ satisfies the condition \eqref{eq:FEMthreshold}, i.e.,
\[
  (hk)^{2\ell}\, k\, R_2 \le C_{\rm FEM, 1}
  \implies (hk)^{2\ell}\, k \lesssim 1
  \implies h \lesssim k^{-1-1/(2\ell)} = k^{-3/2},
\]
then the relative error for the FEM solution of the PML problem is controllably small. Moreover, we see from Theorem~\ref{thm:FEM-ffp} that (with $d=2$) if 
\[
  (hk)^{2\ell}\, k^{3/2}\, R_2^2 \le C
  \implies (hk)^{2\ell}\, k^{3/2} \lesssim 1
  \implies h \lesssim k^{-1-3/(4\ell)} = k^{-7/4},
\]
then the error in computing the far-field pattern is bounded independently of $k$.
Hence our fixed value of $h$ above is sufficient for all of the values of $k\in\{12,24,48\}$.

\ifdefined\journalstyle 
  \newcommand{\widthseven}{6cm}
  \newcommand{\widtheight}{7cm}
\else
  \newcommand{\widthseven}{7cm}
  \newcommand{\widtheight}{8cm}
\fi

\begin{figure} [t] 
 \centering
 \scriptsize{$|u^T|$, $\bsy=\bszero$ \hspace{5.5cm} $|u^S|$, $\bsy=\bszero$} \\
 \includegraphics[width=\widthseven]{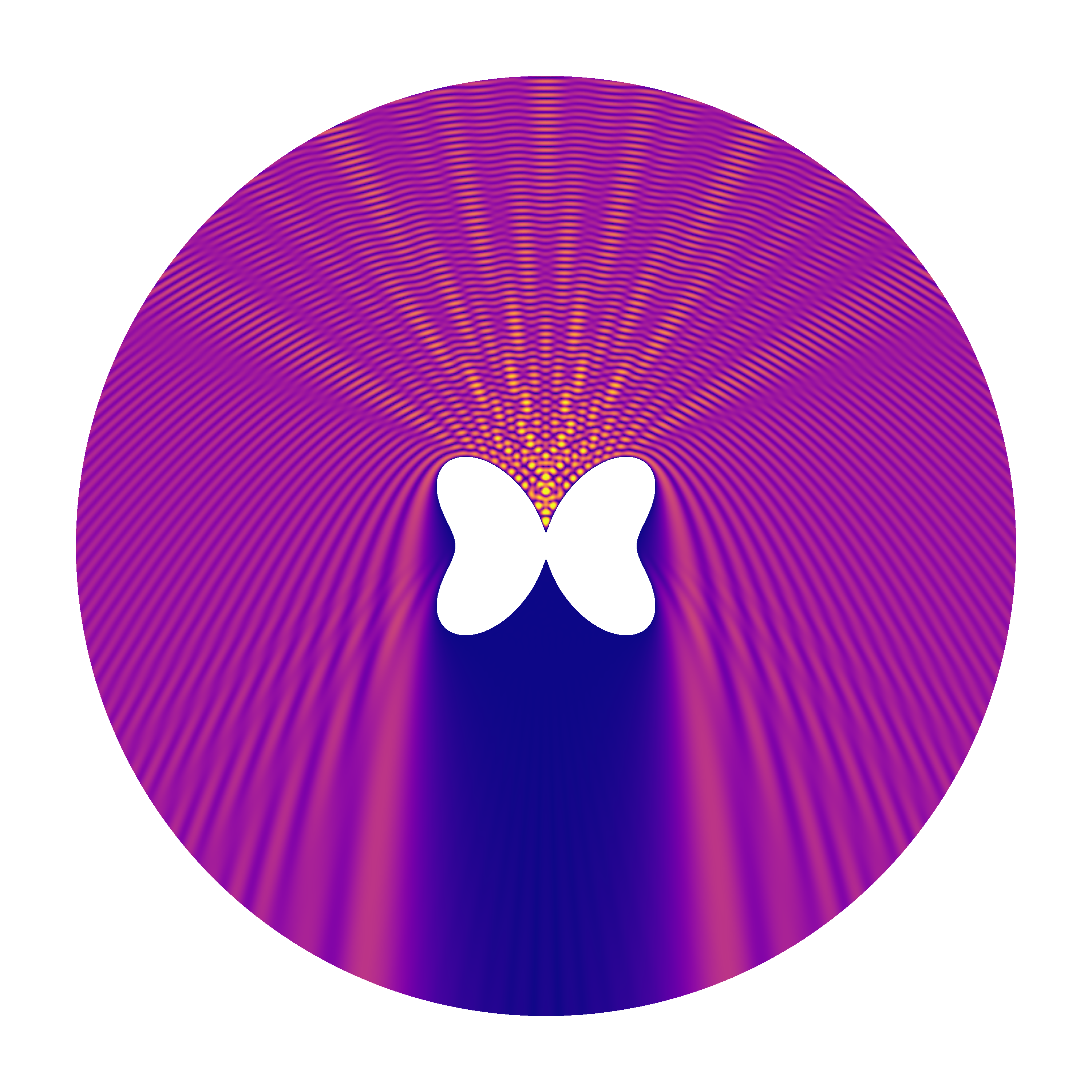} \quad
 \includegraphics[width=\widthseven]{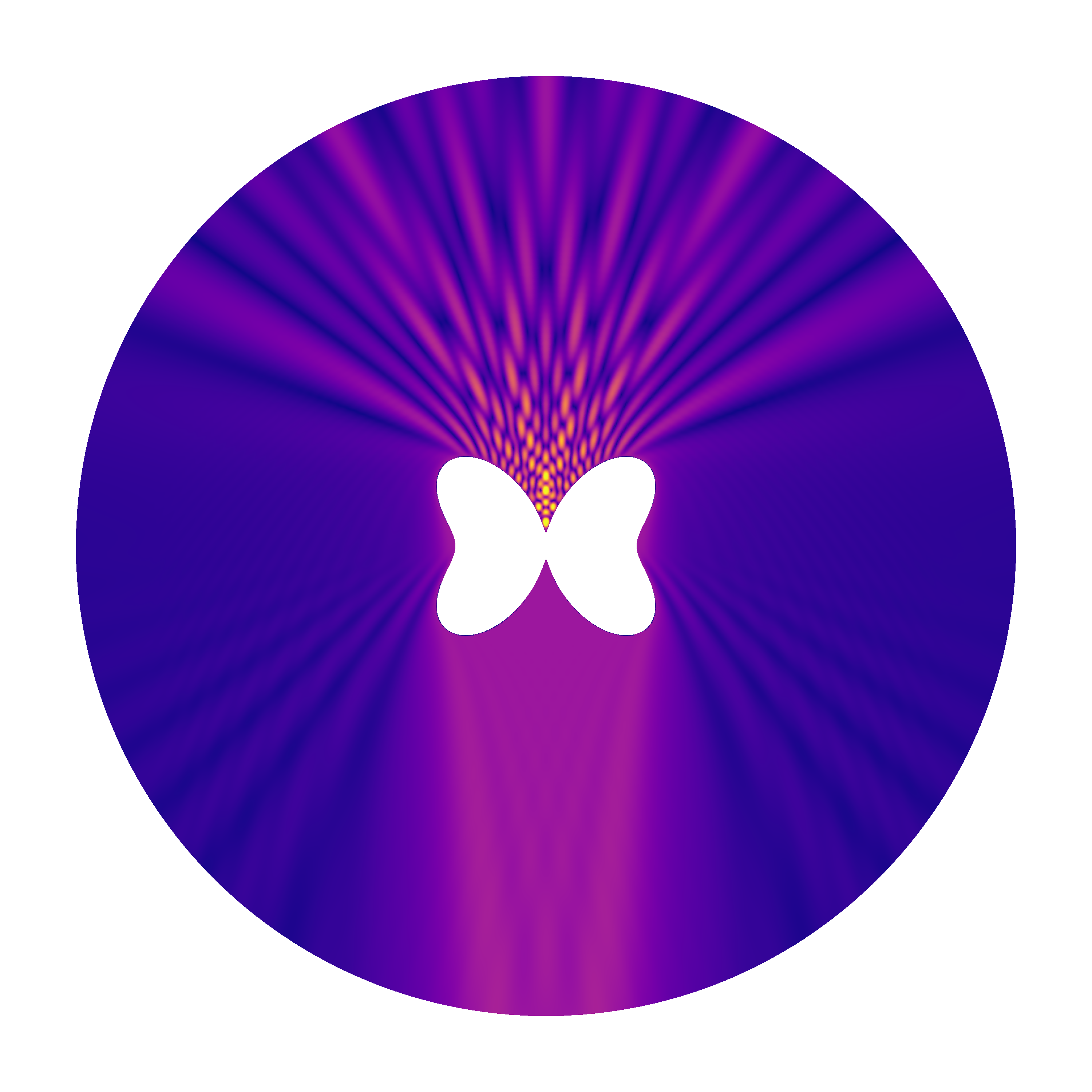} 
 \medskip \\
 \scriptsize{$|u^T|$, $\bsy=(\frac{1}{2},\frac{1}{2},\frac{1}{2},\bszero)$ \hspace{4.5cm} 
 $|u^S|$, $\bsy=(\frac{1}{2},\frac{1}{2},\frac{1}{2},\bszero)$} \\
 \includegraphics[width=\widthseven]{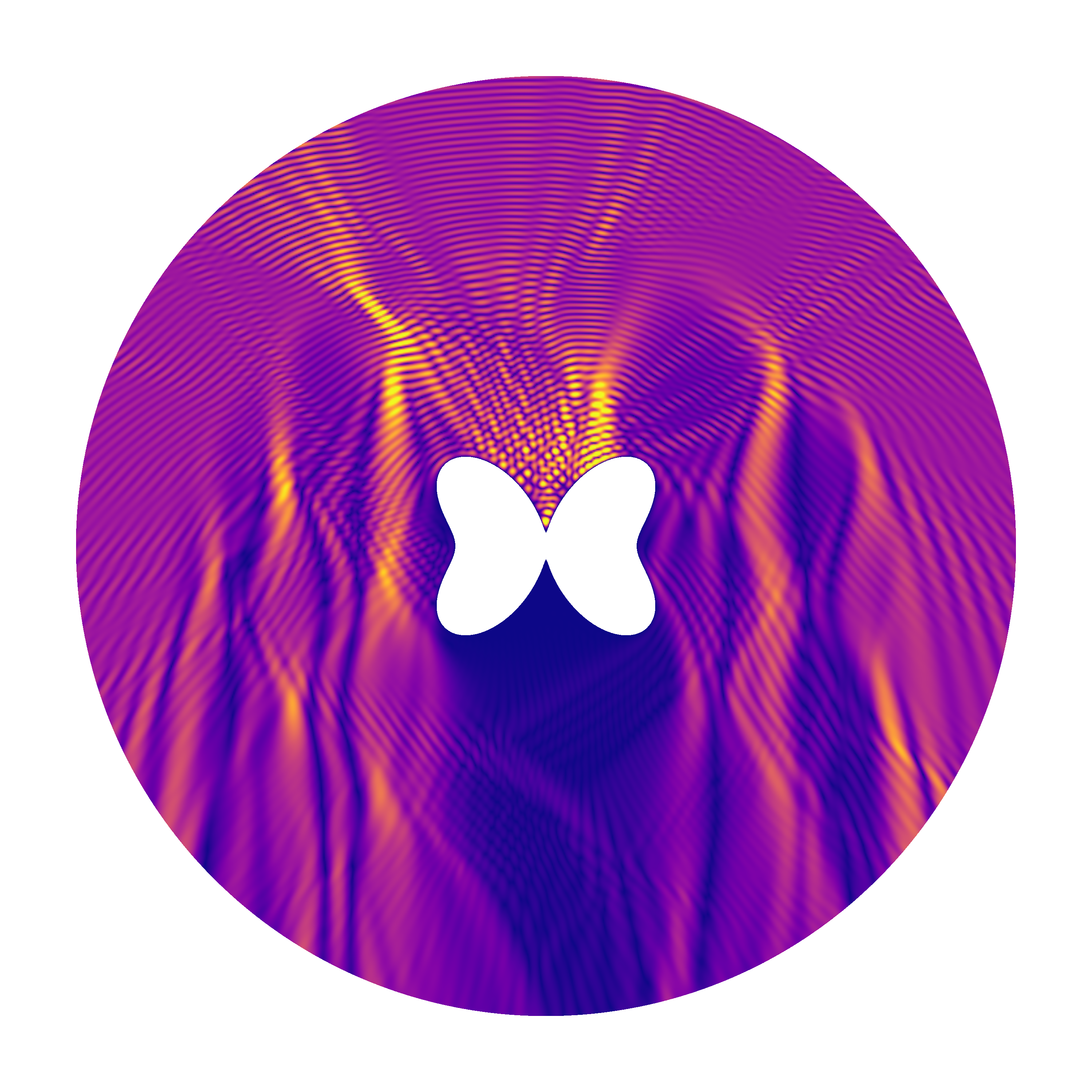} \quad 
 \includegraphics[width=\widthseven]{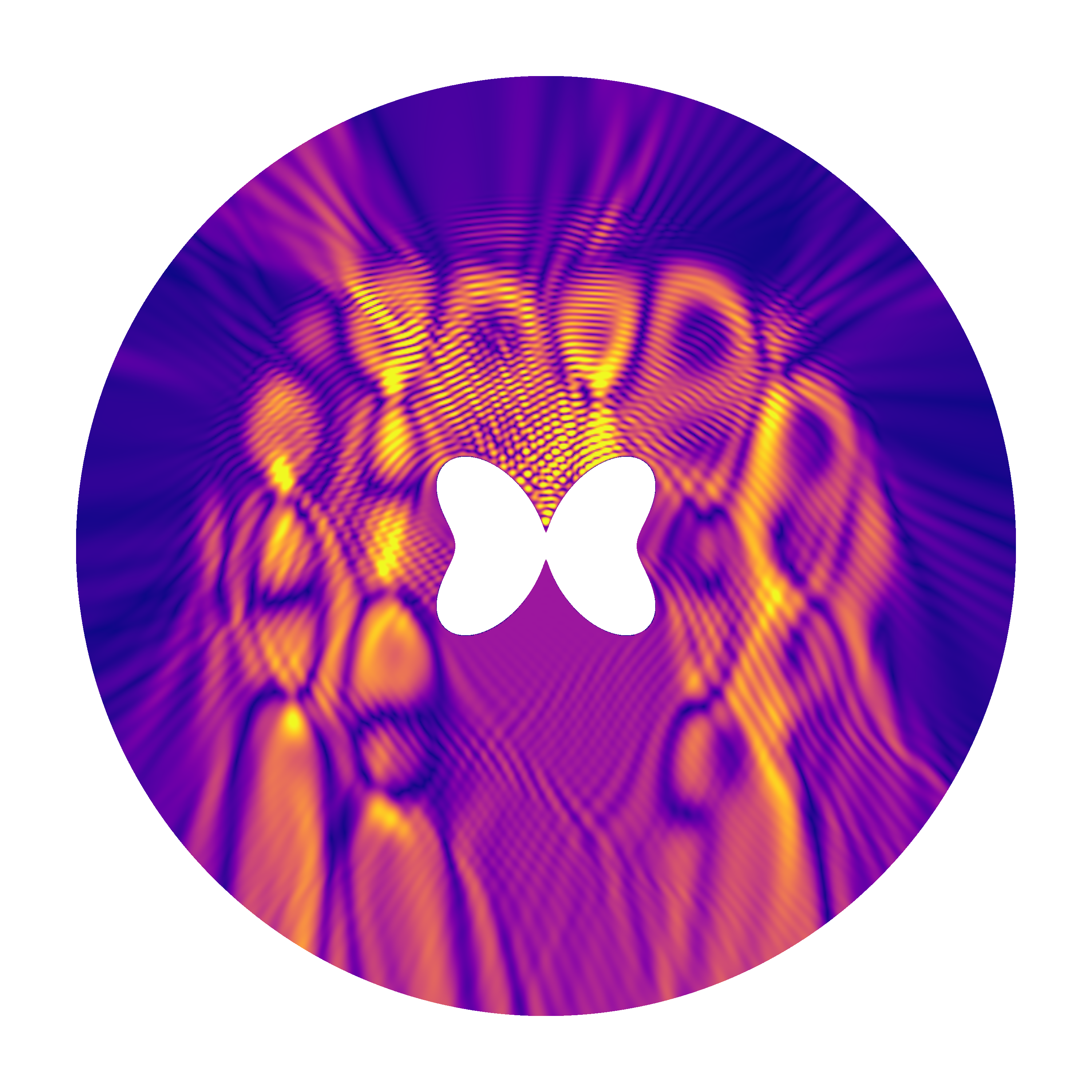} 
 \caption{Computed total field $|u^T|$ and scattered field $|u^S|$ within radius $R_0=4.5$ for $k=48$. Top row: homogeneous field $n(\bsx,\bszero) = 1$. Bottom row: heterogeneous field $n(\bsx,\bsy)$ with $\bsy=(\frac{1}{2},\frac{1}{2},\frac{1}{2},\bszero)$.} \label{fig:soln}
\end{figure} 

Figure~\ref{fig:soln} plots the modulus of the computed total field~$|u^T|$ and scattered field~$|u^S|$ in the circle with radius $R_0=4.5$ for $k=48$. The top two images correspond to the homogeneous field $n(\bsx,\bsy) = n(\bsx,\bszero) = n_0(\bsx)=1$. The bottom two images correspond to one single realization of the heterogeneous field $n(\bsx,\bsy)$ where $\bsy = (\frac{1}{2},\frac{1}{2},\frac{1}{2},\bszero)$, which is equivalent to taking $\bsy = \bshalf := (\frac{1}{2}, \frac{1}{2},\ldots)$ and then truncating the series \eqref{eq:nxy} to $s=3$ terms. Recall that the incident wave comes from the top so we observe some shadow below the butterfly in the plots of $u^T$. We clearly see the effect of the field in the heterogeneous region (circle with radius $R_0-\eta = 3.5$) in the bottom images compared to the top images. Notice that the solutions in the bottom images are not symmetric around the vertical axis of the butterfly. This is because $\sin(j\pi x_1)$ in the definition \eqref{eq:nspec} is not symmetric along the $x_1=1/2$ line. However, since $\sin(j\pi x_1)$ is actually anti-symmetric in $x_1$, later when we average over many realizations with $y_j\sim {\rm Uniform}([-\frac{1}{2},\frac{1}{2}])$, we should observe symmetry in the expected solution.
With $k=48$ the wavelength is $2\pi/k \approx 0.13$ so in each circle with radius $R_0 = 4.5$ we expect to see about $2R_0/0.13\approx 69$ ripples.

\begin{figure} [t] 
 \centering
 \includegraphics[width=\widtheight]{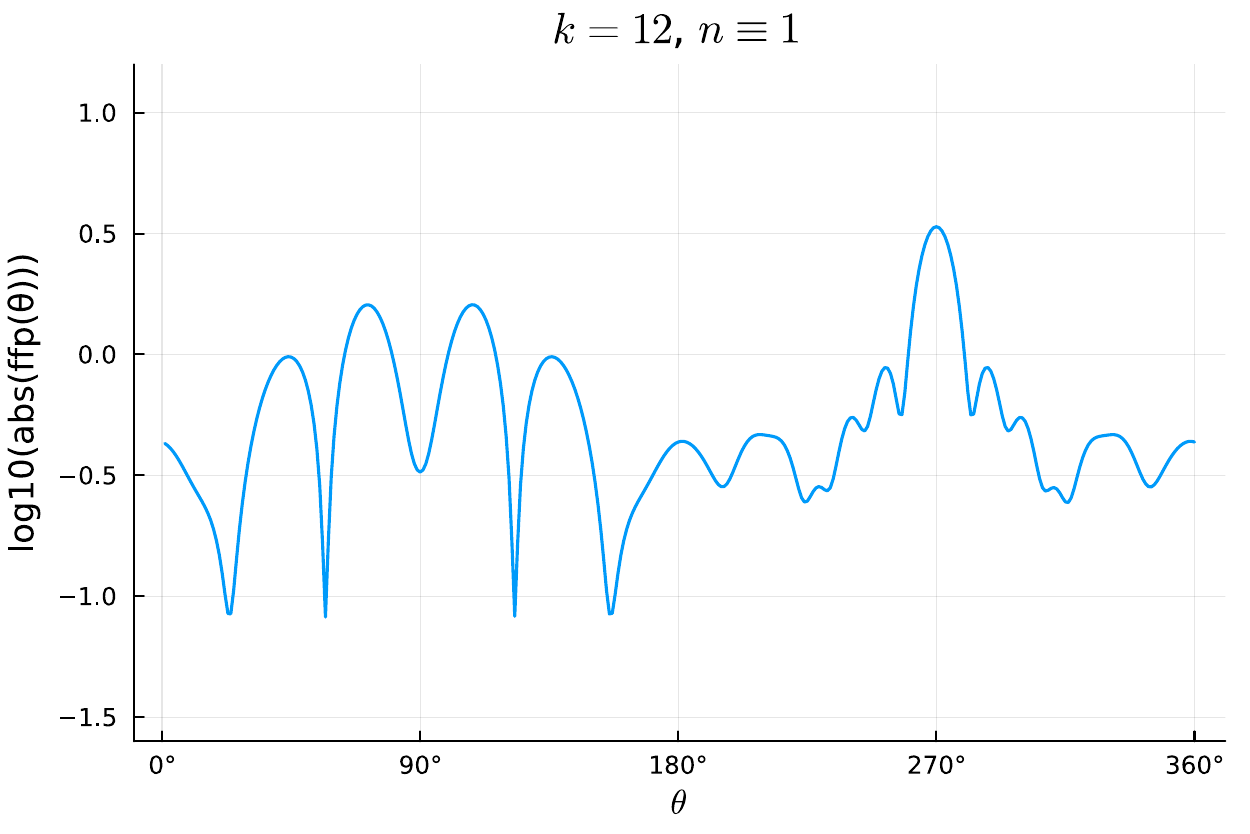} \quad
 \includegraphics[width=\widtheight]{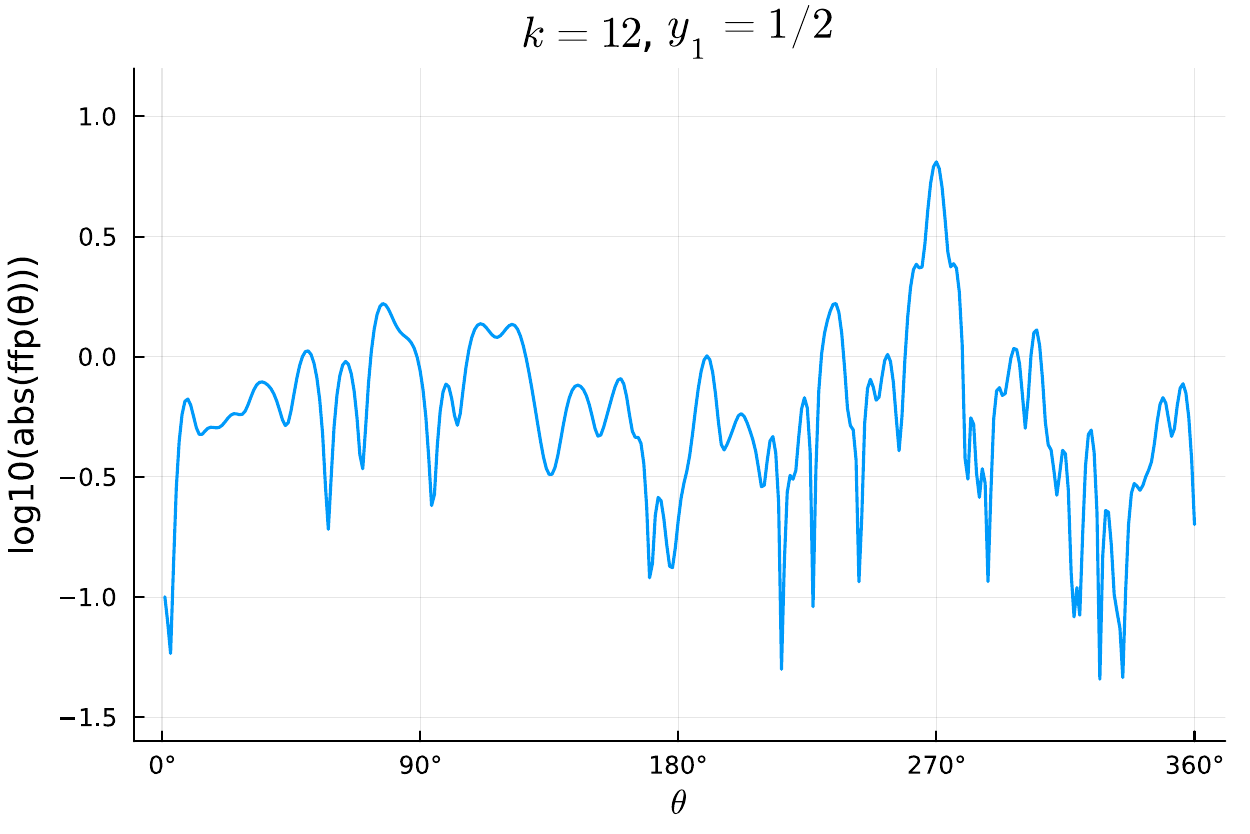} \\
 \includegraphics[width=\widtheight]{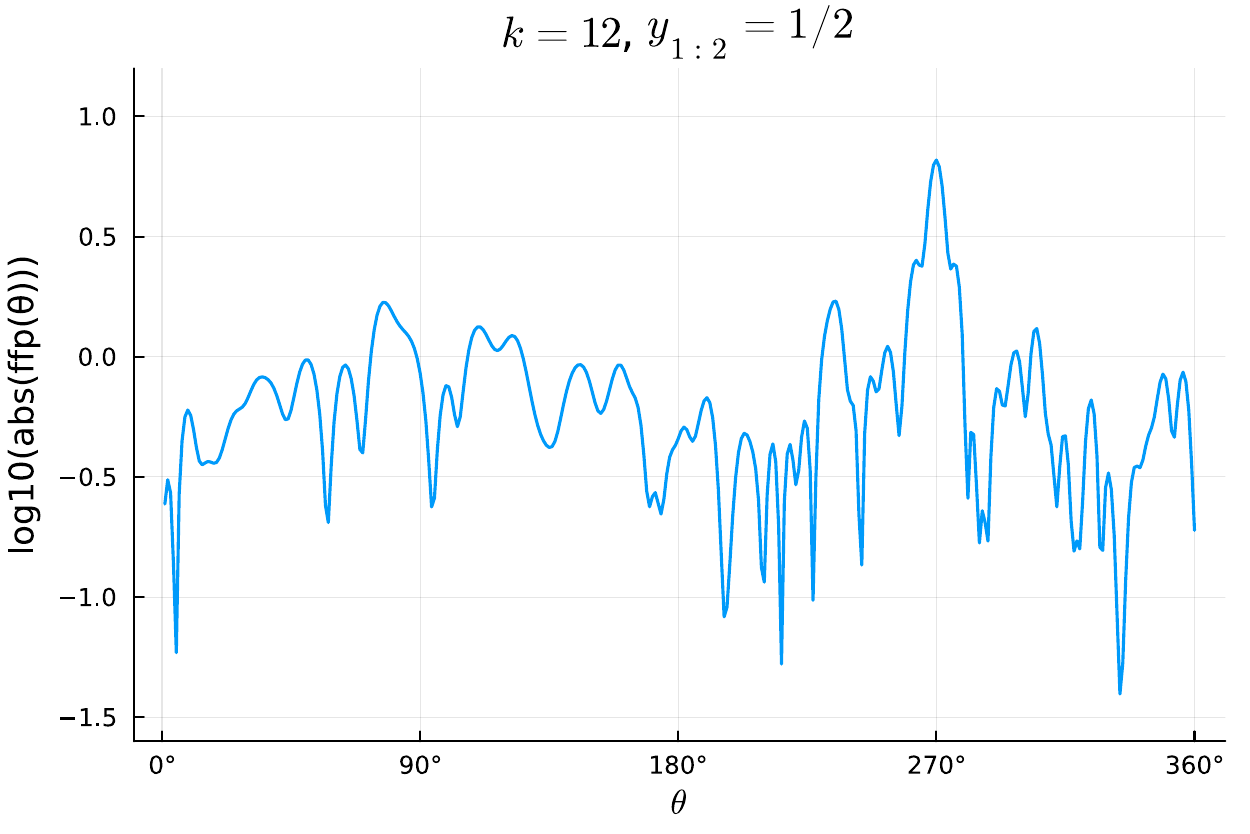} \quad
 \includegraphics[width=\widtheight]{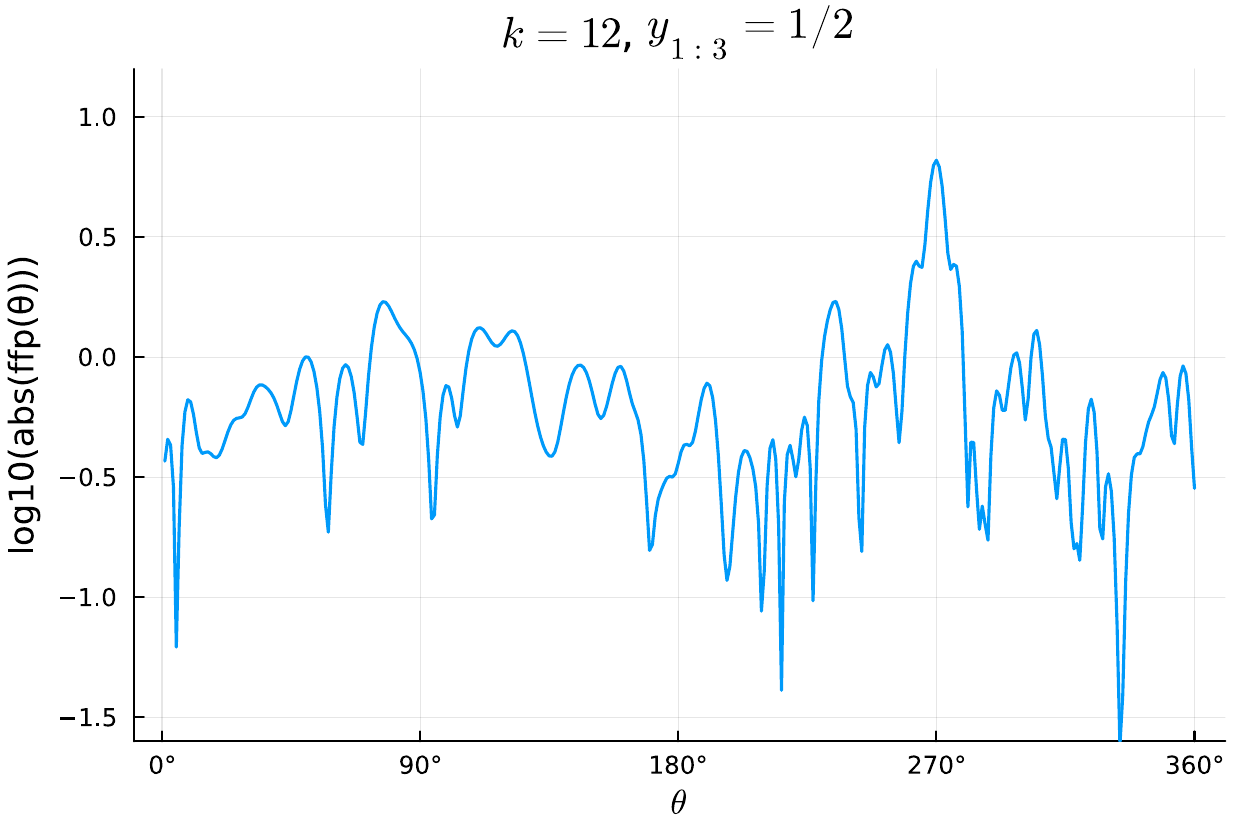} \\
 \includegraphics[width=\widtheight]{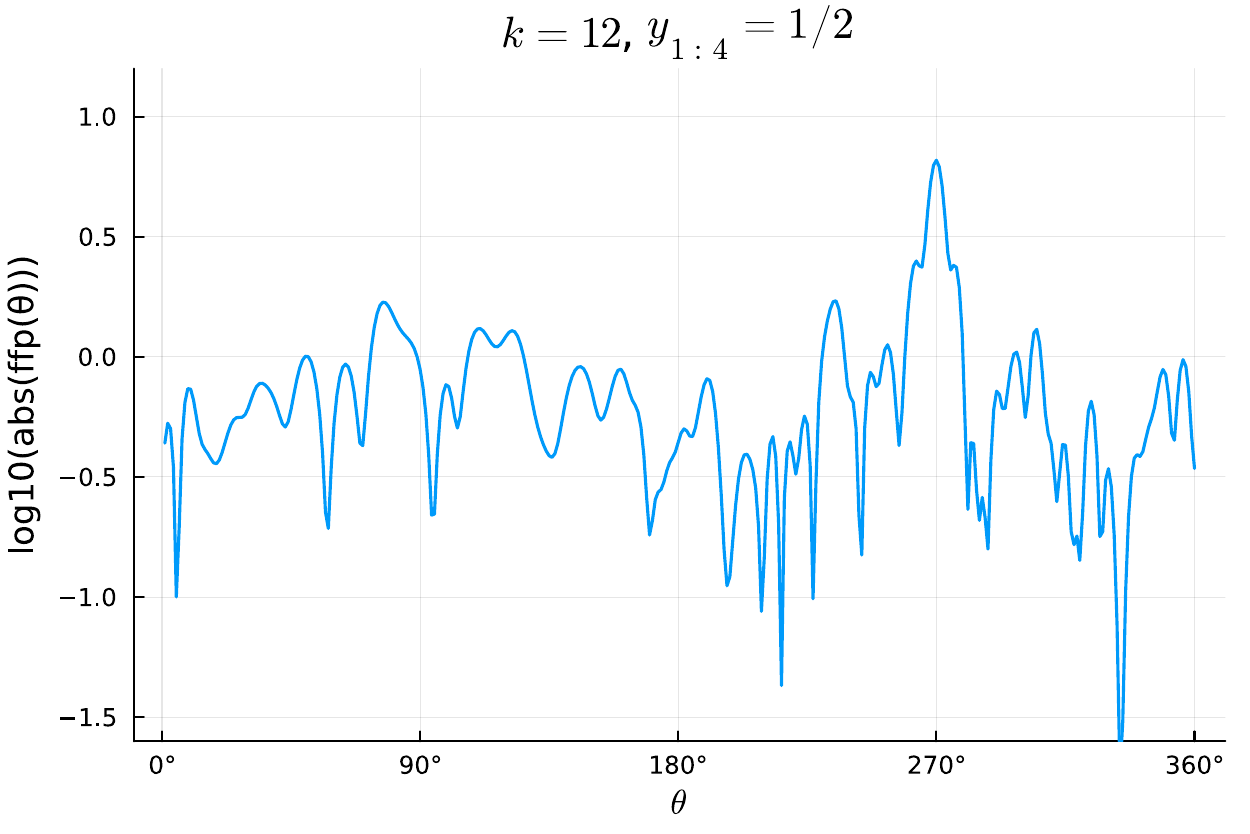} \quad
 \includegraphics[width=\widtheight]{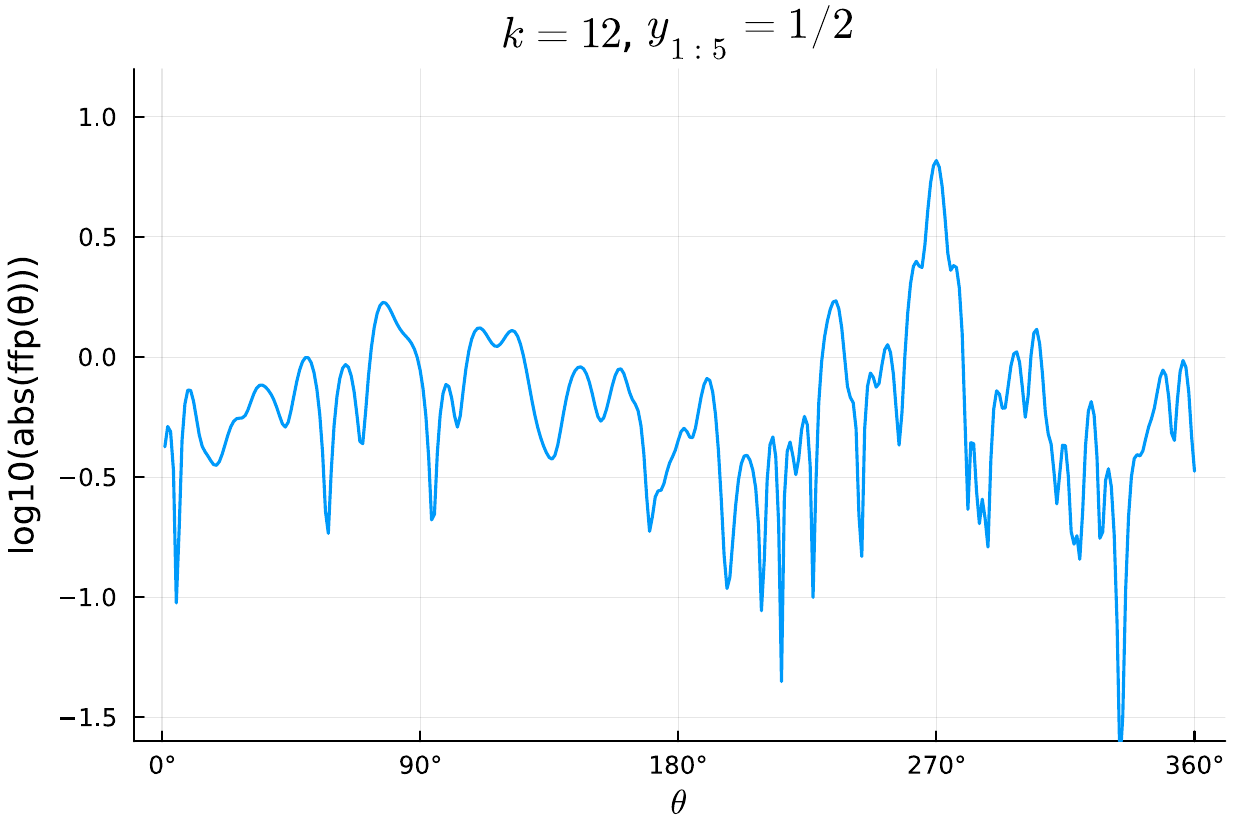} \\
\caption{Far-field pattern for $k=12$ computed with the homogeneous field $n(\bsx,\bszero)=1$ (top left) and with realizations of the heterogeneous field $n(\bsx,\bshalf)$ truncated to $s=1,2,3,4,5$ terms.} \label{fig:ffp}
\end{figure}

Figure~\ref{fig:ffp} plots (the $\log_{10}$ of the modulus) of the far-field pattern for $k=12$ evaluated at angles $1^\circ,2^\circ,\ldots,360^\circ$ for the solutions with homogeneous field $n(\bsx,\bszero)=1$ and with realizations of the heterogeneous field $n(\bsx,\bsy)$ where $\bsy$ is taken successively to be $(\frac{1}{2},\bszero)$, $(\frac{1}{2},\frac{1}{2},\bszero)$, $(\frac{1}{2},\frac{1}{2},\frac{1}{2},\bszero)$, $(\frac{1}{2},\frac{1}{2},\frac{1}{2},\frac{1}{2},\bszero)$, $(\frac{1}{2},\frac{1}{2},\frac{1}{2},\frac{1}{2},\frac{1}{2},\bszero)$, which is equivalent to taking $\bsy = \bshalf$ and truncating the series \eqref{eq:nxy} to $s=1,2,3,4,5$ terms, respectively. We see that the homogeneous case (top left) looks very different to the other five heterogeneous cases, while the five heterogeneous cases look quite similar. From further plots for even higher dimensions we can spot a visible difference up to dimension $s=6$, and in the corresponding plots of the solutions (similar to Figure~\ref{fig:soln} but with $k=12$), we can spot a visible difference only up to dimension $s=4$. These indicate that our problem has a fairly low ``effective dimension''. Similar observations hold for $k=24, 48$, but the far field patterns become rougher as $k$ increases.

\subsection{QMC quadrature}

Since the far-field pattern is calculated using the numerical approximation to $u^{\rm alt}$ in the annulus between radii $4$ and $4.5$, we may take $R = 4.5$ in Theorem~\ref{thm:qmc}. We construct one generating vector up to dimension $s=64$ for an ``embedded'' (see \cite{CKN06}) randomly shifted lattice rule for $N \in \{ 128, 256, 512, 1024\}$ based on POD weights given by \eqref{eq:pod}, taking (since $p\approx 1/3$)
\[
  \delta = 0.1, \quad
  \lambda = 1/1.8, \quad
  \vartheta(\lambda) = \tfrac{2\zeta(2\lambda)}{(2\pi^2)^\lambda} \approx 3.666, \quad
  \beta_j := b_j = 1/j^3.
\]
We do not follow the definition \eqref{eq:Gu-norm}
\[
  \beta_j := \Cstab\,k\,R\,\xi\,b_j = 48.26 \times 48 \times 4.5 \times 0.8319 
  \times 1/j^3
  = 8672/j^3,
\]
since the constant is too big (due to loose theoretical estimates) for practical use. Although this means theoretically that the error bound could scale with $(\Cstab\,k\,R\,\xi)^s$, we do not expect to observe this in practice. In our component-by-component construction we avoid picking a repeated choice of component to improve the quality of the lower-dimensional projections of the lattice rule. This one generating vector is used for all our computations with different values of $k\in \{12, 24, 48\}$, $s\in\{16,32,64\}$ and $N \in \{ 128, 256, 512, 1024\}$. 
Our QMC quadrature is the average $\overline{Q}_{s,N,L}(\Theta) := (1/L) \sum_{\ell=1}^L Q_{s,N,\bsDelta_\ell}(\Theta)$ over $L=10$ independent random shifts
$\bsDelta_1,\ldots,\bsDelta_L \in [0,1]^s$, with $Q_{s,N,\bsDelta_\ell}(\Theta)$ as defined in \eqref{eq:lattice},
and the standard error is estimated via
\[
  \sqrt{\frac{1}{L(L-1)} \sum_{\ell=1}^L \big(
  Q_{s,N,\bsDelta_\ell}(\Theta) - \overline{Q}_{s,N,L}(\Theta) \big)^2}.
\]

\begin{figure} [t] 
 \centering
 \includegraphics[width=\widtheight]{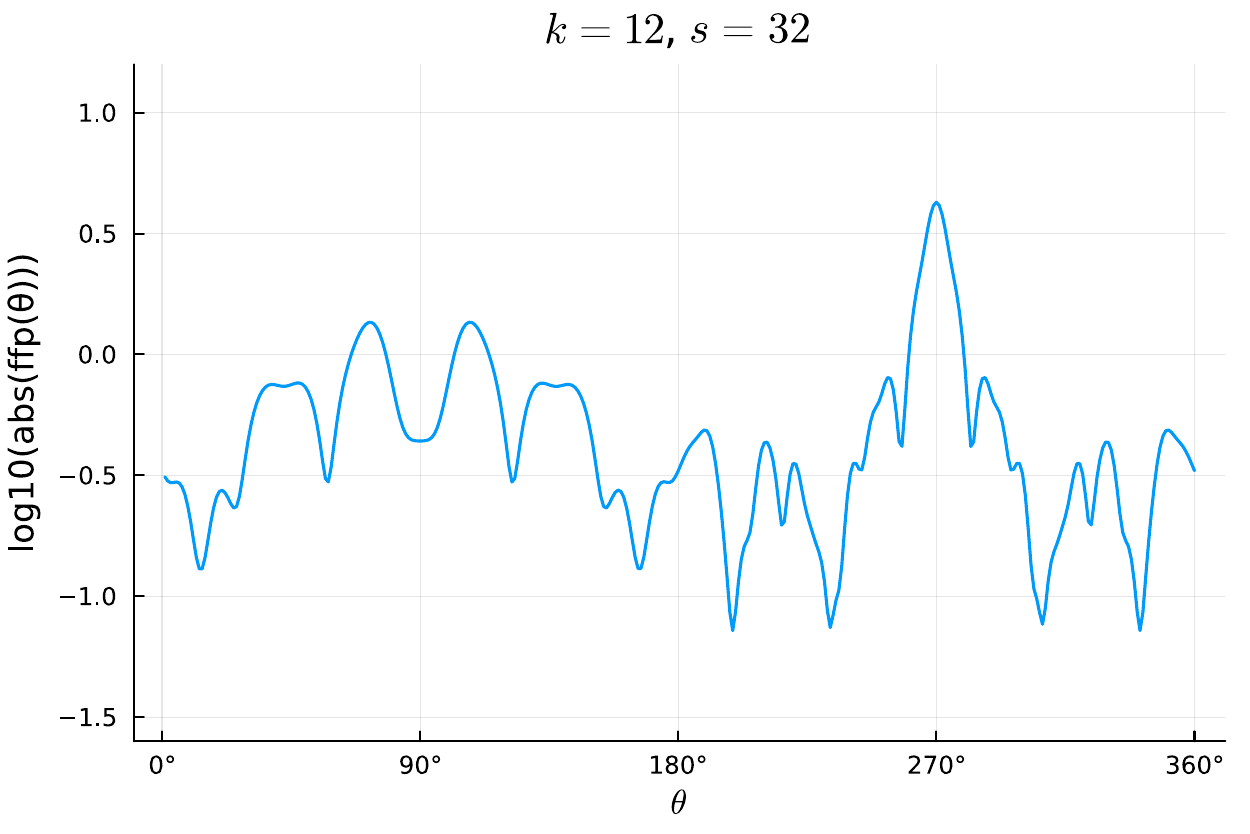} \quad 
 \includegraphics[width=\widtheight]{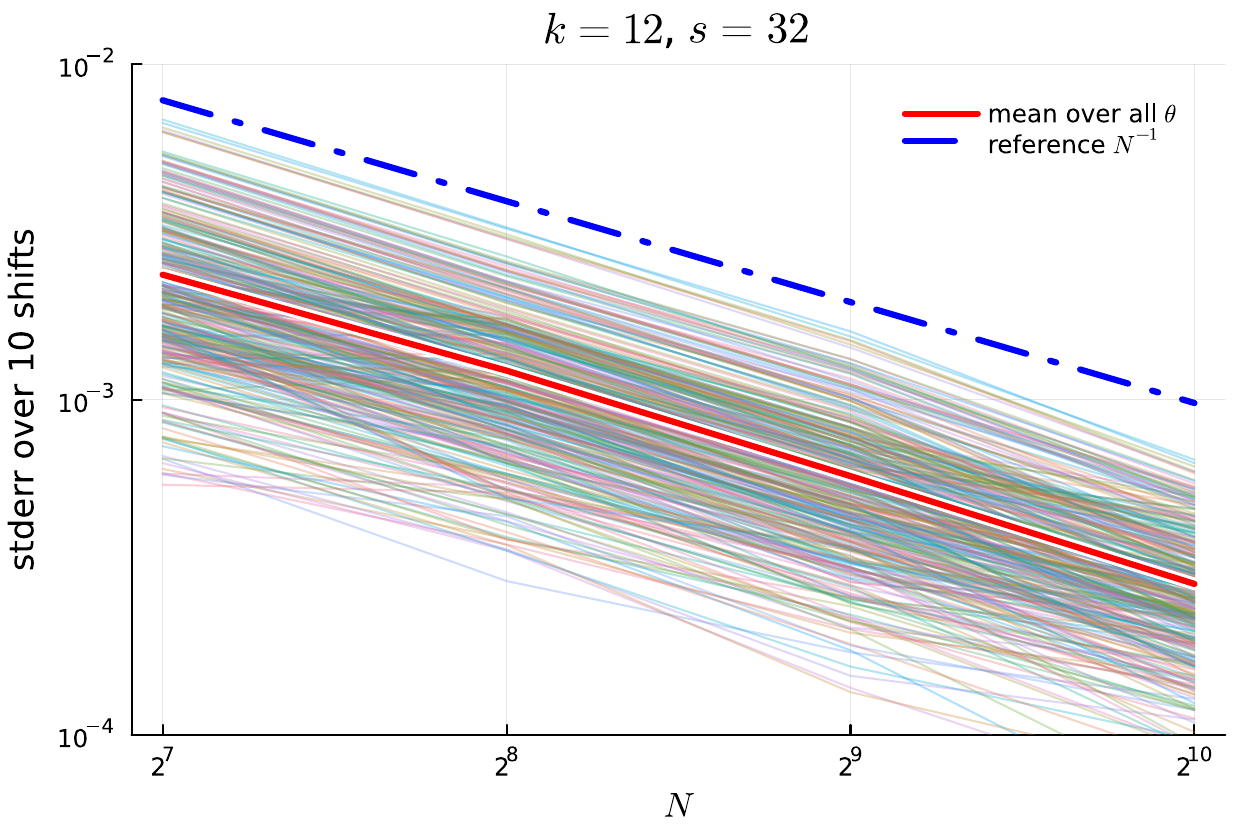} \\
 \includegraphics[width=\widtheight]{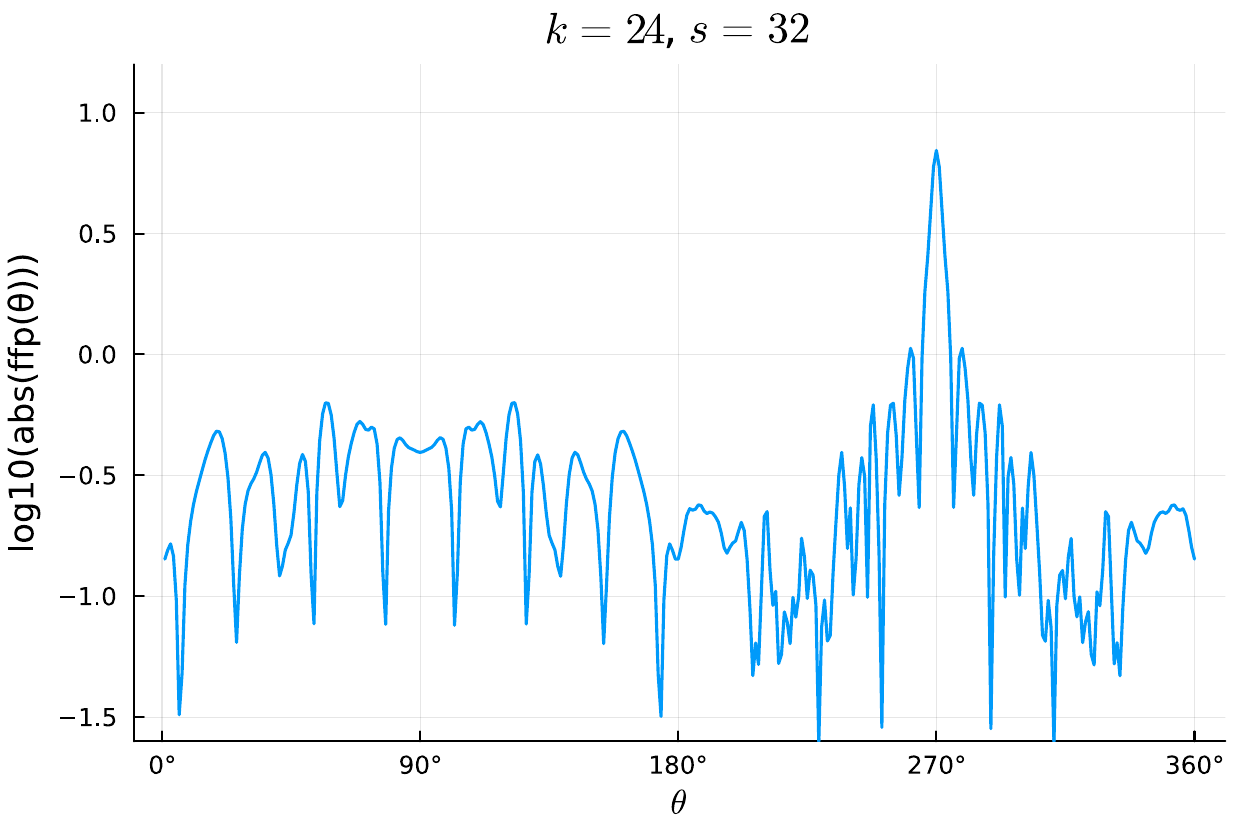} \quad 
 \includegraphics[width=\widtheight]{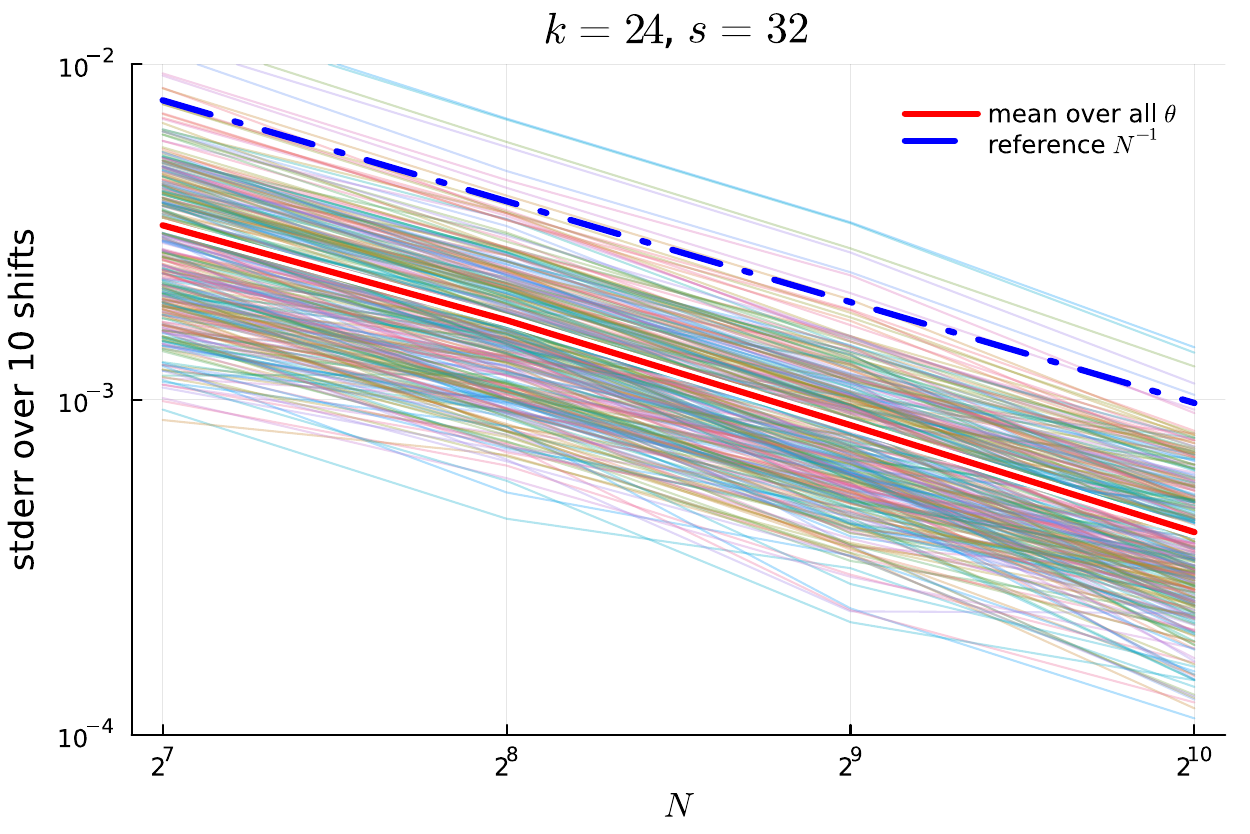} \\
 \includegraphics[width=\widtheight]{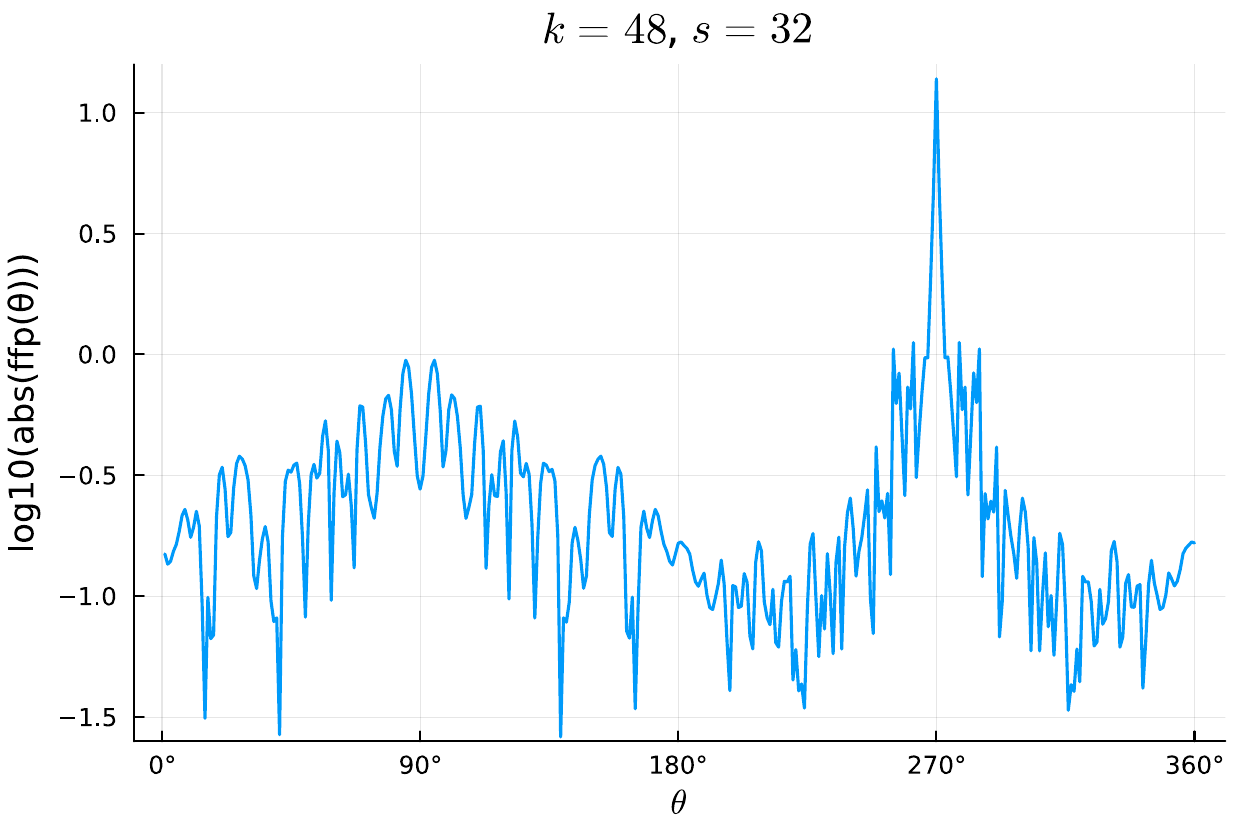} \quad 
 \includegraphics[width=\widtheight]{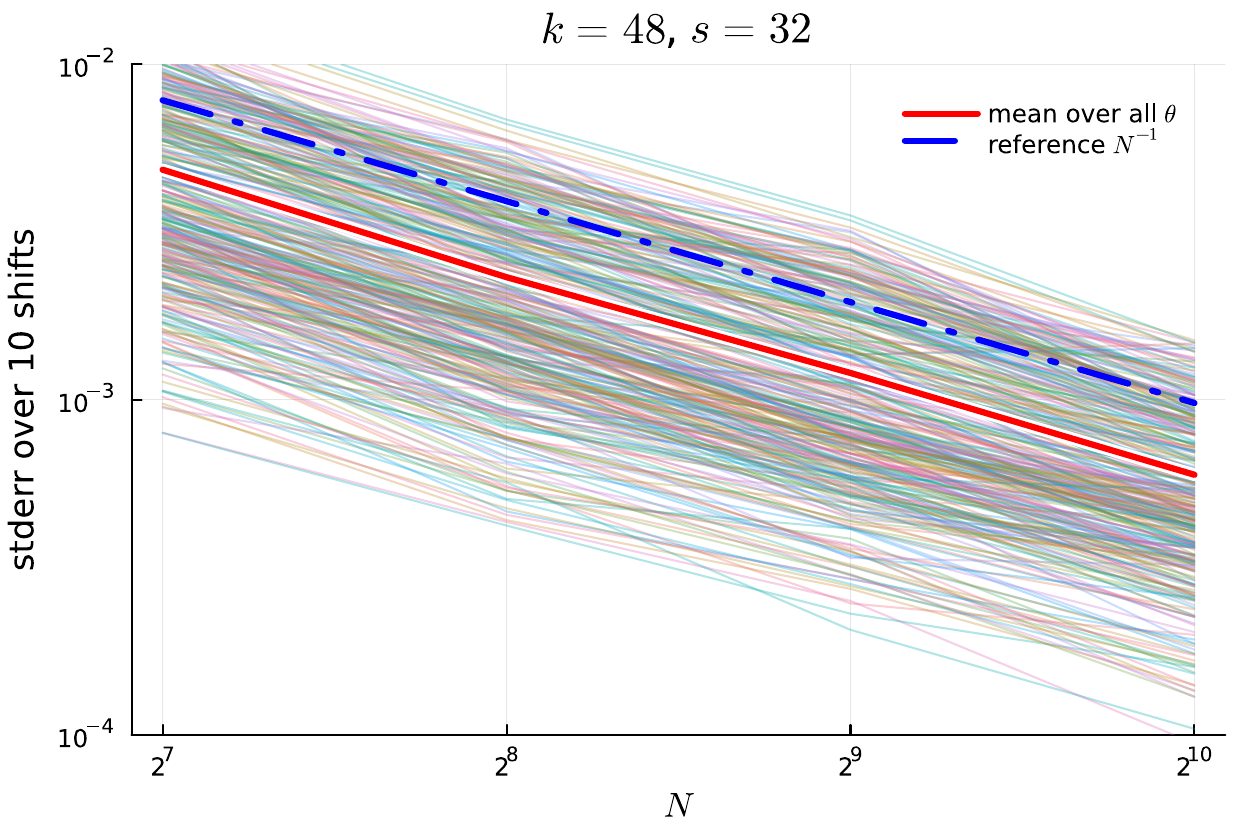} \\
 \caption{Expected value of the far-field pattern and estimated standard errors for the random field with $s=32$ terms and varying $k \in \{12, 24, 48\}$ from top down.} \label{fig:vary-k}
\end{figure}

Figure~\ref{fig:vary-k} plots the expected value of the far-field pattern and the estimated standard errors computed using $L=10$ random shifts of $N\in\{ 128, 256, 512, 1024\}$ realizations of the random field with $s=32$ terms, varying $k \in \{12,24,48\}$ from top down. 
Individual lines on the right plots correspond to the estimated errors at angles $1^\circ, 2^\circ,\ldots, 360^\circ$. The red solid lines are the means and the blue dashed lines are reference lines for order $N^{-1}$, demonstrating that we are achieving the order $N^{-1}$ convergence rate in all cases. The error is $k$ dependent, as predicted by Theorem \ref{thm:qmc}. The far-field pattern becomes significantly rougher as $k$ increases, and the standard error lines shift noticeably higher as $k$ increases, all as we expected. We observe symmetry in the expected value of the far-field pattern along $90^\circ$ and $270^\circ$ as we predicted.
We also made similar plots for $k=48$ and varying $s\in \{16,32,64\}$, but we could not spot any visible difference between the three different value of $s$. This indicates that for $k=48$ the effective dimension of the far-field pattern is lower than $16$. (Note that we use the same set of random shifts for varying $N$ and $s$. In particular, when $s$ varies the initial components are embedded.)

\begin{figure} [t] 
 \centering
 \includegraphics[width=\widtheight]{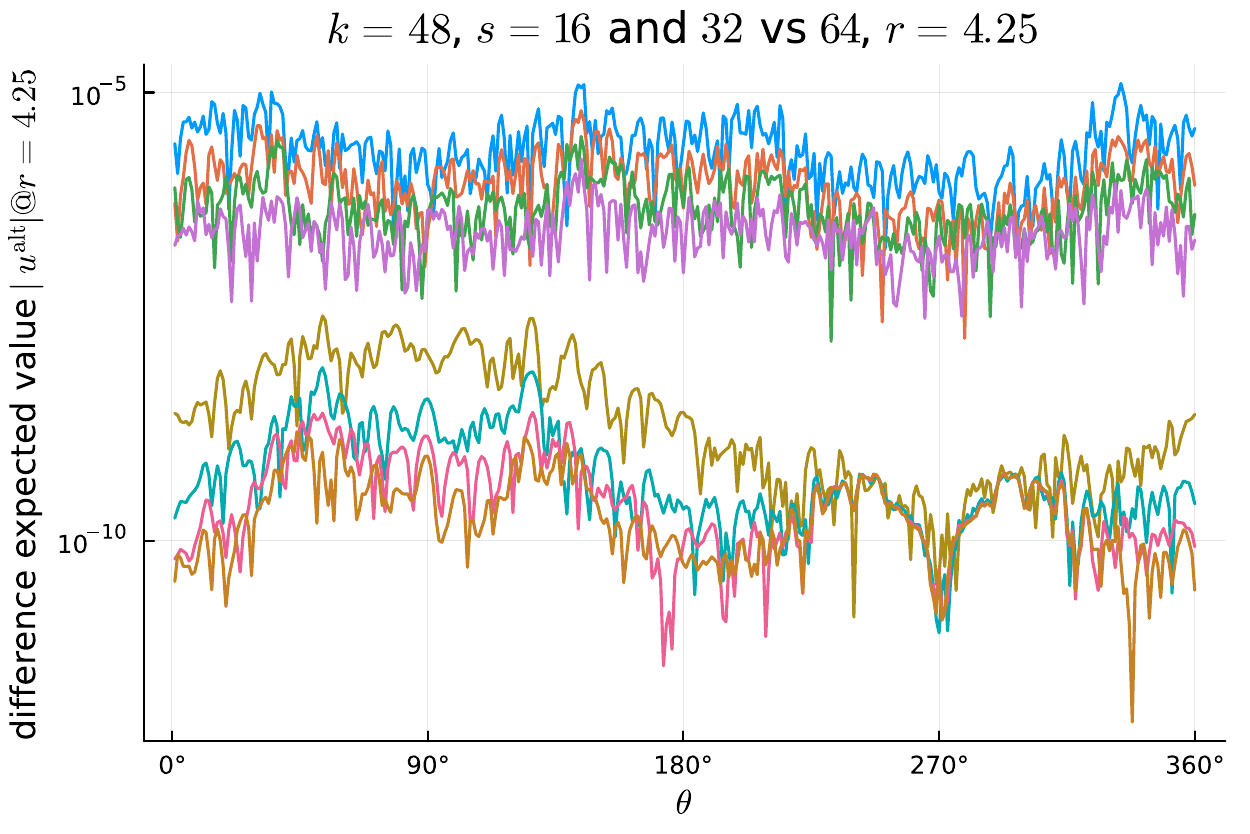} \quad 
 \includegraphics[width=\widtheight]{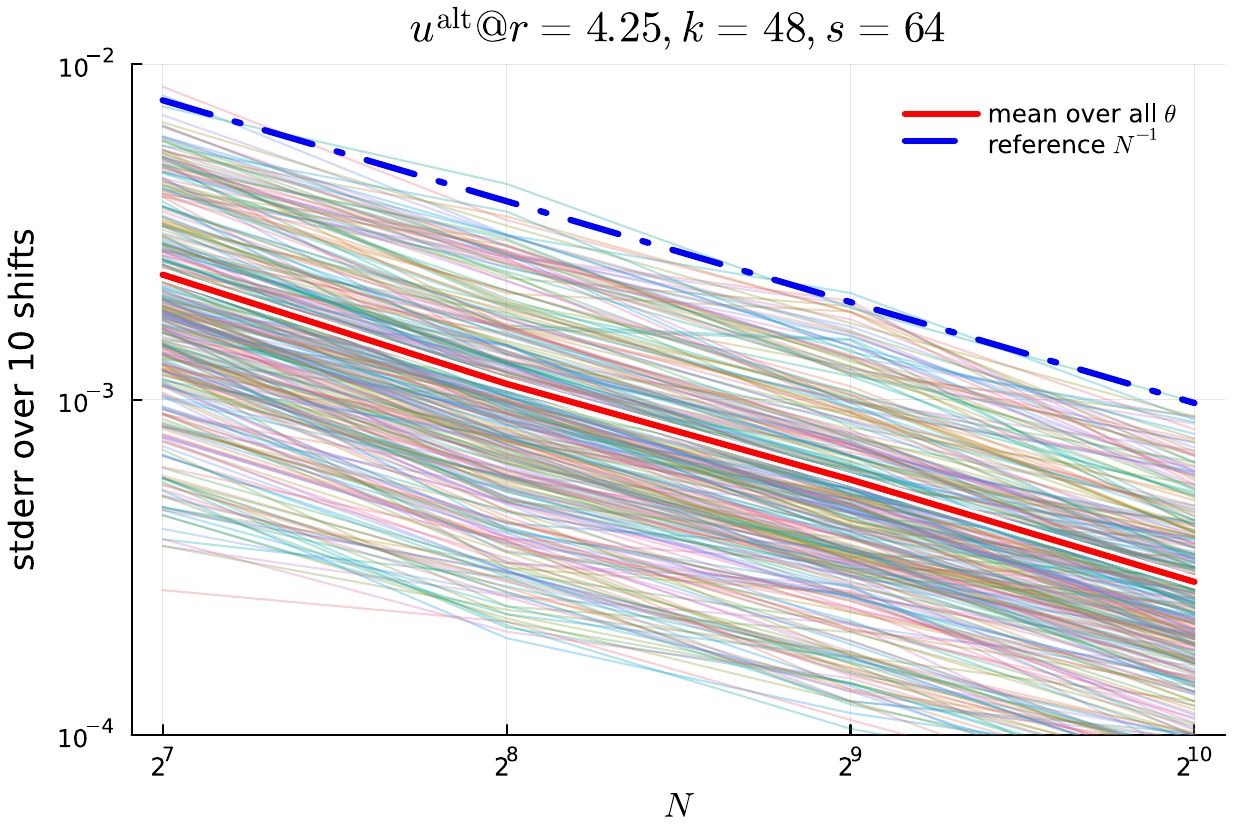} 
 \caption{Left: 
   Differences in $|u^{\rm alt}(\bsx)|$ between $s=16$ and $s=64$ terms (top $4$ curves, 
   $4$ different values of $N$), and 
   between $s=32$ and $s=64$ terms (bottom $4$ curves). 
   Right: Estimated standard errors in computing $|u^{\rm alt}(\bsx)|$ for $k=48$ and $s = 64$.} \label{fig:uScircle}
\end{figure}

Figure~\ref{fig:uScircle} shows that the solution $u^{\rm alt}$ has a higher effective dimension than the corresponding far-field pattern. For $k=48$, we compute the expected value of $|u^{\rm alt}(\bsx)|$ with $\bsx$ on the circle of radius $4.25$ at angles $1^\circ, 2^\circ,\ldots, 360^\circ$, and as before we
use $L=10$ random shifts of $N\in\{128,256,512,1024\}$ realizations of the random field, but now varying the number of terms $s\in \{16,32,64\}$ used in \eqref{eq:nxy}. 
The top $4$ curves (one for each value of $N$) on the left of Figure~\ref{fig:uScircle} plot in log-scale the absolute difference between the estimates for $s=16$ and $s=64$ against the polar angle of $\bsx$. These curves are about $10^{-7}$ which is a rough guide to the dimension truncation error for $s=16$.
The bottom $4$ curves plot the difference between $s=32$ and $s=64$, 
and are about $10^{-9}$, thus the effective dimension of $|u^{\rm alt}(\bsx)|$ is more than $32$. 
On the right of Figure~\ref{fig:uScircle} we plot the estimated standard errors (obtained with $L=10$ random shifts) for $k=48$ and $s=64$, where individual lines again correspond to $\bsx$ on the same circle at different angles. These indicate an optimal QMC convergence rate of $\calO(N^{-1})$. The corresponding plots for $s = 16,32$ (not included here) turned out to be almost identical, indicating that the $\mathcal{O}(N^{-1})$ convergence is dimension independent (as predicted by Theorem \ref{thm:qmc}, since $k$ is fixed here).

\ifdefined\journalstyle {}
\else
  \clearpage
\fi 

\section{Conclusion}

In this paper we combined QMC quadrature with PML truncation and FEM approximation to solve the
Helmholtz equation formulated on the infinite propagation domain, after
scattering by the heterogeneity as well as a bounded impenetrable scatterer. 
One important special case is the plane-wave sound-soft scattering problem 
where our quantity of interest is the far-field pattern of the scattered field. We developed theory for 
(a) dimension truncation in parameter space; (b) QMC quadrature for computing expectation;
and (c) PML truncation and FEM approximation error. Our error estimates are explicit in $s$
(the dimension truncation parameter), $N$ (the number of QMC points), $h$
(the FEM grid size) and most importantly $k$~(the Helmholtz wavenumber). The
method is also exponentially accurate with respect to the PML truncation
radius. We give numerical experiments computing the expected value of the far-field patterns for domains containing about $70$ wavelengths. We observe close to $\calO(N^{-1})$ convergence which is optimal for randomly shifted lattice rules.

\section*{Acknowledgments} I.~G.~Graham and E.~A.~Spence acknowledge support from
EPSRC grant EP/S003975/1, and E.~A.~Spence acknowledges support from EPSRC grant EP/R005591/1.
F.~Y.~Kuo and I.~H. Sloan
acknowledge the support from the Australian Research Council Discovery Project DP210100831.
We are grateful to Andrew Gibbs and Stephen Langdon for sharing their data for some computed far-field patterns.
This research includes computations using the computational cluster Katana
supported by Research Technology Services at UNSW Sydney~\cite{katana}.

\end{document}